\RequirePackage{fix-cm}
\documentclass{svjour3}                     
\smartqed  
\usepackage{tabularx}
\usepackage{amsmath}
\usepackage{amssymb}
\usepackage{dsfont}
\usepackage{rotating}
\usepackage{multirow}
\usepackage{epstopdf}
\usepackage{newtxmath}
\usepackage{mathtools}

\usepackage{multirow} 
\usepackage{longtable,tabu} 
\usepackage{hhline}
\usepackage[linesnumbered,ruled,vlined]{algorithm2e}

\usepackage{hyperref}
\hypersetup{
	colorlinks=true,
	linkcolor=black, 
	citecolor=black, 
	urlcolor=blue  } 

\newcommand{\R}{\mathds{R}}

\renewcommand{\vv}{\mathbf{v}}
\newcommand{\xx}{\mathbf{x}}
\newcommand{\nn}{\mathbf{n}}
\newcommand{\q}{\mathbf{q}}

\newcommand{\dt}{{\Delta t}}
\newcommand{\dx}{{\Delta x}}
\newcommand{\dy}{{\Delta y}}
\newcommand{\dz}{{\Delta z}}
\newcommand{\Q}{\mathbf{q}}

\newcommand{\CFL}{\text{CFL}}

\newcommand{\B}{\mathbf{B}}
\newcommand{\A}{\mathbf{A}}
\newcommand{\f}{\mathbf{f}}

\renewcommand{\vv}{\mathbf{v}}

\newcommand{\pdd}[2]{\frac{\partial #1}{\partial #2}}
\newcommand{\refv}{\text{ref}}

\newcommand{\Mc}{{M}_c}  
\newcommand{\Ma}{{M}_b}  

\newcommand{\ijk}{i  j  k}
\newcommand{\ihp}{i+\frac{1}{2} j k}
\newcommand{\ihm}{i+\frac{1}{2} j k}
\newcommand{\jhp}{i j+\frac{1}{2} k}
\newcommand{\jhm}{i j-\frac{1}{2} k}
\newcommand{\khp}{i j k+\frac{1}{2}}
\newcommand{\khm}{i j k-\frac{1}{2}}
\newcommand{\ip}{i+1 j k}
\newcommand{\im}{i-1 j k}
\newcommand{\jp}{i j+1 k}
\newcommand{\jm}{i j-1 k}
\newcommand{\kp}{i j k+1}
\newcommand{\km}{i j k-1}

\newcommand{\Fop}{\vmathbb{F}} 
\newcommand{\Bop}{\vmathbb{B}} 
\newcommand{\Gop}{\vmathbb{G}} 
\newcommand{\Cop}{\vmathbb{C}} 
\newcommand{\Dop}{\vmathbb{D}} 
\newcommand{\Lop}{\vmathbb{L}} 
\newcommand{\Hop}{\vmathbb{H}} 

\newcommand{\scheme}{SIFV-EB~}

\newfont{\numerikEleven}{ecrm1000}
\newfont{\numerikTen}{cmss10}
\newfont{\numerikNine}{cmss9}
\newfont{\numerikEight}{cmss8}

\journalname{Journal of Scientific Computing}
\begin{document}

\title{A structure-preserving semi-implicit IMEX finite volume scheme for ideal magnetohydrodynamics at all Mach and Alfv{\'e}n numbers}

\titlerunning{A semi-implicit finite volume scheme for MHD at all Mach and Alfv{\'e}n numbers}        

\author{Walter Boscheri \and
		Andrea Thomann }


\institute{W. Boscheri \at
	    Laboratoire de Mathématiques UMR 5127 CNRS, Universit{\'e} Savoie Mont Blanc, 73376 Le Bourget du Lac, France \\
		\email{walter.boscheri@univ-smb.fr}
		\and
		A. Thomann \at
        Universit{\'e} de Strasbourg, CNRS, Inria, IRMA, F-67000 Strasbourg, France \\
		\email{andrea.thomann@inria.fr}
}

\date{Received: date / Accepted: date}

\maketitle

\begin{abstract}
We present a divergence-free semi-implicit finite volume scheme for the simulation of the ideal magnetohydrodynamics (MHD) equations which is stable for large time steps controlled by the local transport speed at all Mach and Alfv\'en numbers. An operator splitting technique allows to treat the convective terms explicitly while the hydrodynamic pressure and the magnetic field contributions are integrated implicitly, yielding two decoupled linear implicit systems. The linearity of the implicit part is achieved by means of a semi-implicit time linearization. This structure is favorable as second-order accuracy in time can be achieved relying on the class of semi-implicit IMplicit-EXplicit Runge-Kutta (IMEX-RK) methods.
In space, implicit cell-centered finite difference operators are designed to discretely preserve the divergence-free property of the magnetic field on three-dimensional Cartesian meshes. The new scheme is also particularly well suited for low Mach number flows and for the incompressible limit of the MHD equations, since no explicit numerical dissipation is added to the implicit contribution and the time step is scale independent. Likewise, highly magnetized flows can benefit from the implicit treatment of the magnetic fluxes, hence improving the computational efficiency of the novel method. The convective terms undergo a shock-capturing second order finite volume discretization to guarantee the effectiveness of the proposed method even for high Mach number flows.
The new scheme is benchmarked against a series of test cases for the ideal MHD equations addressing different acoustic and Alfv\'en Mach number regimes where the performance and the stability of the new scheme is assessed.
\subclass{MSC 65 \and MSC 68}
\end{abstract}


\section{Introduction} \label{sec.intro}

Compressible magnetized plasma flows are mathematically modeled by the equations of magnetohydrodynamics (MHD), that constitute a first order hyperbolic system of nonlinear partial differential equations (PDE). Potential practical applications concern astrophysics to describe plasma flows in the magnetosphere of neutron stars, or inertial and magnetic confinement fusion, in order to calibrate or design fusion reactors for civil energy production.

The dynamics of plasma flows can involve different time scales, that correspond to the characteristics of the governing PDE system. The multi-scale nature of the MHD equations can be analyzed by considering the Mach number $\Mc$, which is the ratio between the fluid velocity and the sound speed, and the Alfv\'en number $\Ma$, that represents the ratio between the fluid velocity and the speed of the Alfv\'en wave related to the magnetic field. Different values of $\Mc$ and $\Ma$ yield different fluid regimes, and, consequently, the nature of the governing equations can change passing for instance from hyperbolic to elliptic. A classical example is given by low Mach number flows, which are dominated by acoustic waves, hence recovering an incompressible fluid regime in the stiff relaxation limit for $\Mc \to 0$. On the other hand, the material speed is prevalent in high Mach number flows, thus leading to compressible fluids with shock waves. Another source of stiffness arises for highly magnetized flows, i.e. for $\Ma \to 0$, that are typically encountered in Tokamak scenarios for magnetic confinement fusion applications.

The numerical resolution of the MHD equations is therefore very challenging if the numerical method is required to deal with the entire spectrum of the different time scales. Explicit Godunov-type finite volume schemes \cite{Godunov1959,HLL1983,Balsara2004} exhibit excellent shock-capturing properties that can perfectly handle high Mach number flows. However, if different regimes coexist in the flow, two main problems are generated by the different time scales that can be assumed by the model: i) the computational time step obeys a stability condition that is constraint by the magnitude of the fastest characteristic speeds, hence leading to very small time steps in the stiff relaxation limits of the model which make the simulation extremely slow, and, ultimately, even impossible to be accomplished; ii) the numerical viscosity does not have the correct scaling but is of order $\mathcal{O}(1/\Mc)$ which is especially severe for low Mach number flows and therefore the accuracy of the solution is compromised \cite{Dellacherie1,GH}.

To tackle these challenges, a first strategy consists in rescaling the dissipation term in the numerical flux so that the numerical viscosity becomes independent of the Mach number \cite{bouchut2007multiwave,barsukow2021,barsukow2023,Birke2023}. A fully implicit time integration has then to be adopted to circumvent the restrictive time step condition \cite{BK2022,Viallet2011}, which leads to a highly nonlinear system to be solved with a large number of unknowns, or local time stepping techniques can be foreseen \cite{altmann2009local,LagrangeLTS}.
Another possibility is instead given by the implicit treatment of only one part of the system to be solved while keeping the remaining part explicit \cite{Degond2,Dumbser_Casulli16,BosPar2021,Chen2022,SIMHD_Dumbser2019,Zampa2023,birke2024,LukPesTho2024}. The terms responsible to a possible stiffness of the governing model undergo an implicit discretization in which no numerical diffusion is needed, since stability follows from implicitness, so that the above mentioned problems of explicit solvers can be simultaneously solved. This technique makes often use of the class of IMplicit-EXplicit (IMEX) time integrators forwarded e.g. in \cite{AscRuuSpi,BP2017,BosRus,PR_IMEX,BPR2017,Hofer}, that have been developed with the specific purpose of dealing with multiple time scales, typically a slow and a fast time scale. Furthermore, these IMEX schemes are proven to be Asymptotic Preserving (AP) and, some even Asymptotically Accurate (AA) meaning that they recover a consistent discretization of the limit model also with high order accuracy.

Besides the multi-scale character, the MHD equations embed another constraint on the divergence of the magnetic field, which must remain zero all times because no magnetic monopoles exist in nature. This is an involution that has to be carefully taken into account by the numerical scheme, since in principle it is not satisfied at the discrete level. In the literature, there are three main approaches to satisfy the solenoidal property of magnetic field in a numerical scheme: i) the use of staggered grids is typically adopted to cancel the divergence errors, which requires the knowledge of the electric field \cite{yee1966,BalsaraSpicer1999}; ii) divergence cleaning methods \cite{munz2000,dedner2002} add one hyperbolic equation to the MHD system in order to transport away the divergence errors, which require the determination of the cleaning speed as a new parameter; iii) constrained transport schemes on collocated grids that also make use of the electric field \cite{fey2003,helzel2011,helzel2013}.

The aim of this paper is to design a numerical scheme for the ideal MHD equations that can handle fast time scales which might be induced by both the stiff limits $\Mc \to 0$ and $\Ma \to 0$, while retaining shock-capturing properties in the case of high Mach number flows. To face the multi-scale nature of the model, we rely on the implicit-explicit approach, whereas the divergence-free involution is preserved by means of mimetic finite difference operators \cite{HymanShashkov1997,Margolin2000,Lipnikov2014} on collocated meshes. For the low Mach number regime, we follow the schemes forwarded in \cite{BosPar2021} for compressible viscous flows. This approach has been extended to magnetohydrodynamics in \cite{birke2024,Chen2022}, thus devising asymptotic preserving methods for the hydrodynamic sub-system of the MHD model up to third order space-time accuracy. On staggered meshes, semi-implicit finite volume schemes have been developed in \cite{SIMHD_Dumbser2019}, and hybrid finite volume/finite element methods are presented in \cite{Zampa2023}. These schemes are second order accurate in space and first order accurate in time. For the low Alfv\'en regime, the only semi-implicit scheme currently available in the literature is presented in \cite{3splitMHD} with first order in time accuracy, where staggered grids are used, and the implicit sub-system results to be non-linearly coupled, hence requiring to resort to a fixed point algorithm to linearize the problem combined with a nested Newton algorithm which is very challenging in practise.

In this work, we propose a novel semi-implicit finite volume scheme that is second order accurate in both space and time. The semi-implicit IMEX strategy of \cite{BosFil2016} is adopted to devise an overall second-order all Mach and Alfv\'en number scheme that only involves linear implicit systems for the fast scales. Consequently, neither Newton solvers nor Picard iteration techniques will be needed. A standard explicit TVD finite volume method is employed for the non-linear convective terms. To respect the involution on the magnetic field, simple collocated grids are used onto which we design a structure-preserving div-curl operator that is retrieved as the second order version of the div-curl discontinuous Galerkin operators recently forwarded in \cite{divcurlB}. To this end, the MHD model is reformulated in terms of the magnetic potential which generates the magnetic field by its curl.

The paper is structured as follows. In Section \ref{sec.pde}, we briefly recall the model of ideal MHD equations and its non-dimensional formulation in terms of Mach and Alfv\'en number.
Then, the numerical scheme is introduced in Section \ref{sec.numscheme} where we first focus on the operator splitting in a one-dimensional framework where we subsequently describe the multi-dimensional time semi-discrete scheme based on the magnetic potential instead of evolving the magnetic field. In the latter we focus on the cell-centered operators that are designed to preserve exactly the div-curl involution constraint thus yielding a divergence-free fully discrete multi-dimensional finite volume scheme.
In Section \ref{sec.numtest}, a series of test problems in up to three space dimensions is presented to numerically verify the performance of the numerical schemes at various Mach and Alfv\'en numbers. A section providing a summary and an outlook to future developments concludes the paper.

\section{Governing equations} \label{sec.pde}
We consider the ideal MHD equations on a three-dimensional domain $\Omega \subset \mathbb{R}^3$ where $\xx = (x,y,z) \in \Omega$ is the vector of spatial coordinates and $t \in \mathbb{R}^+_0$ denotes the time coordinate. They write
\begin{equation}
    \label{sys:MHD}
    \pdd{}{t} \left( \begin{array}{c}
    	\displaystyle \rho \\[2mm] \displaystyle \rho \vv \\[2mm] \displaystyle \rho E \\[2mm] \displaystyle \B
    \end{array} \right) + \nabla \cdot \left( \begin{array}{c}
    \displaystyle \rho \vv \\[2mm] \displaystyle \rho \vv \otimes \vv + \left(p + \frac{\|\B\|^2}{8 \pi}\right) \mathbb{I} - \frac{1}{4\pi} \B \otimes \B \\[2mm] \displaystyle \left(\rho E + p + \frac{\|\B\|^2}{8 \pi}\right) \vv - \frac{1}{4 \pi} \B (\vv \cdot \B) \\[2mm] \displaystyle \vv \otimes \B - \B \otimes \vv
\end{array} \right) = \mathbf{0}.
\end{equation}
Therein, the density is denoted by $\rho$, the pressure by $p$, the velocity field by $\vv = (u,v,w)$ and the magnetic field by $\B = (B_x, B_y, B_z)$ which obeys
\begin{equation}
	\nabla \cdot \B = 0,
	\label{eqn.divBFree}
\end{equation}
hence imposing that the magnetic field lines are always closed. The identity matrix is denoted by $\mathbb{I}$. The total energy density is composed of the internal energy $\rho e$, the kinetic energy $\rho k$ and the magnetic energy density $ \rho m$ as follows:
\begin{equation}
    \rho E = \rho e + \rho k + \rho m \qquad \text{with} \qquad k = \frac{1}{2}\|\vv\|^2, \quad m = \frac{\|\B\|^2}{8 \pi}.
\end{equation}
The system is closed by the ideal gas equation of state (EOS)
\begin{equation}
    \rho e = \frac{p}{\gamma - 1},
    \label{eqn.EOS}
\end{equation}
with $\gamma$ representing the ratio of specific heats.
Thus, the pressure is a linear function of the total energy density $\rho E$ and can be computed as $p = (\gamma - 1)(\rho E - \rho k - \rho m)$.
Further, the sound speed, the so-called total and directional Alfv\'en speeds are defined respectively by
\begin{equation}
    c = \sqrt{\frac{\gamma p}{\rho}}, \quad b = \sqrt{\frac{\|\B\|^2}{4 \pi \rho}}, \quad b_{n} = \sqrt{\frac{(\B \cdot \nn)^2}{4 \pi \rho}},
\end{equation}
where $\nn$ denotes a unit normal vector pointing to some arbitrary direction.
The ideal MHD system \eqref{sys:MHD} exhibits 8 characteristic speeds $\lambda_{\{1\ldots8\}}$ in normal direction $\nn$ given by
\begin{equation}
    \lambda_{1,8} = \vv \cdot \nn \mp c_f, \quad \lambda_{2,7} = \vv \cdot \nn \mp b_n, \quad \lambda_{3,6} = \vv \cdot \nn \mp c_s, \quad \lambda_{4,5} = \vv \cdot \nn,
\end{equation}
where $c_{s,f}$ denote the slow and fast magneto-sonic wave speeds which read
\begin{equation}
    c_{s,f} = \sqrt{\frac{1}{2}\left(c^2 + b^2 \mp \sqrt{(c^2 + b^2)^2 - 4 c^2 b_n^2}\right)}.
\end{equation}
For a detailed discussion of the ideal MHD equations \eqref{sys:MHD} and its structure see e.g. \cite{BalsaraRoe}.

\subsection{Scaling of the ideal MHD equations}
The MHD system \eqref{sys:MHD} is characterized by multiple time scales, which depend on the physical phenomena under consideration. A sub-system of the ideal MHD equations is indeed given by the Euler equations of gas dynamics, in which a stiffness might be originated in the momentum equations by the pressure waves. In this case, the sound speed is much faster than the local fluid velocity, thus weakly compressible flows are retrieved. The scaling parameter is the Mach number, defined as $\Mc = \|\vv\|/c$, and these flows are typically affected by a low Mach regime. Another interesting case is given by large magnetic fields, where the Alfv\'en speeds $b$ are significantly faster than the local velocity or even faster than the sound speed $c$. Analogously to the Mach number, we can define a parameter $\Ma = \|\vv\|/b$ associated to the Alfv\'en speed $b$, which is referred to as Alfv\'en Mach number.

To analyze the influence of those parameters on the dynamics of the model, we consider the non-dimensional formulation of the ideal MHD equations \eqref{sys:MHD}.
To that end, we employ a decomposition of a dimensional variable $q$ in the associated non-dimensional variable $\tilde{q}$ and a reference value $q_\refv$ that contains the units of measurements as well.
Thus we can write $q = \tilde{q} q_\refv$.
The convective scale is given by $v_\refv = x_\refv/t_\refv$, where $x_\refv$ and $t_\refv$ denote the reference length and time respectively.
Further, we have a reference density $\rho_\refv$, sound speed $c_\refv$, pressure $p_\refv = \rho_\refv c_\refv^2$, specific internal energy $e_\refv = c_\refv^2$ and a reference Alfv\'en speed $b_\refv$ associated to the magnetic field.
Inserting this decomposition into the governing equations, the rescaled non-dimensional MHD system writes
\begin{equation}
	\label{sys:MHD_nd}
	\pdd{}{t} \left( \begin{array}{c}
		\displaystyle \tilde{\rho} \\[2mm] \displaystyle \tilde{\rho} \tilde{\vv} \\[2mm] \displaystyle \tilde{\rho} \tilde{E} \\[2mm] \displaystyle \tilde{\B}
	\end{array} \right) + \nabla \cdot \left( \begin{array}{c}
		\displaystyle \tilde{\rho} \tilde{\vv} \\[2mm] \displaystyle \tilde \rho \tilde \vv \otimes \tilde \vv + \left(\frac{\tilde p}{M_c^2} + \frac{1}{M_b^2}\frac{\|\tilde\B\|^2}{2}\right) \mathbb{I} - \frac{1}{M_b^2} \tilde\B \otimes \tilde\B \\[2mm] \displaystyle \left(\tilde\rho \tilde E + \frac{\tilde p}{M_c^2} + \frac{1}{M_b^2}\frac{\|\tilde \B\|^2}{2}\right) \tilde \vv - \frac{1}{M_b^2} \tilde\B (\tilde\vv \cdot \tilde\B) \\[2mm] \displaystyle \tilde\vv \otimes \tilde\B - \tilde\B \otimes \tilde\vv
	\end{array} \right) = \mathbf{0},
\end{equation}
where the total energy is given by
\begin{equation}
    \tilde E = \frac{1}{M_c^2}\tilde e + \frac{1}{M_b^2}\tilde m + \tilde k.
\end{equation}
Introducing the plasma-beta parameter $\beta = \tilde{p}/\tilde{m}$, which is given by ratio between (reference) pressure and (reference) magnetic pressure, the non-dimensional formulation \eqref{sys:MHD_nd} can be rewritten as
\begin{equation}
	\label{sys:MHD_nd_beta}
	\pdd{}{t} \left( \begin{array}{c}
		\displaystyle \tilde{\rho} \\[2mm] \displaystyle \tilde{\rho} \tilde{\vv} \\[2mm] \displaystyle \tilde{\rho} \tilde{\tilde{E}} \\[2mm] \displaystyle \tilde{\B}
	\end{array} \right) + \nabla \cdot \left( \begin{array}{c}
		\displaystyle \tilde{\rho} \tilde{\vv} \\[2mm] \displaystyle \tilde \rho \tilde \vv \otimes \tilde \vv + \left(\frac{\tilde p}{M_c^2} + \frac{1}{M_c^2\beta}\frac{\|\tilde\B\|^2}{2}\right) \mathbb{I} - \frac{1}{M_c^2\beta}\tilde\B \otimes \tilde\B \\[2mm] \displaystyle \left(\tilde\rho \tilde{\tilde E} + \tilde p + \frac{1}{\beta}\frac{\|\tilde \B\|^2}{2}\right) \tilde \vv - \frac{1}{\beta} \tilde\B (\tilde\vv \cdot \tilde\B) \\[2mm] \displaystyle \tilde\vv \otimes \tilde\B - \tilde\B \otimes \tilde\vv
	\end{array} \right) = \mathbf{0},
\end{equation}
where the total energy is given by
\begin{equation}
    \tilde{\tilde E} = \tilde e + \frac{1}{\beta}\tilde m + M_c^2\tilde k,
\end{equation}
which is closer to the formulation known from the Euler equations with $\Mc$ in front of the kinetic energy contribution \cite{BosPar2021}.

%

\subsection{Novel 3-split form of the ideal MHD equations}
Let us consider only one space dimension, so that the MHD system \eqref{sys:MHD} reduces to
\begin{equation}
	\frac{\partial \Q}{\partial t} + \frac{\partial \f(\Q)}{\partial x} = \mathbf{0},
	\label{eqn.MHD1D}
\end{equation}
with the definitions
\begin{equation}
	\Q = \left( \begin{array}{c}
		\displaystyle \rho \\[2mm] \displaystyle \rho u \\[2mm] \displaystyle \rho v \\[2mm] \displaystyle \rho w \\[2mm] \displaystyle \rho E \\[2mm] \displaystyle B_x \\[2mm] \displaystyle B_y \\[2mm] \displaystyle B_z \\[2mm]
	\end{array} \right), \qquad \f(\Q) = \left( \begin{array}{c} \displaystyle \rho u \\[2mm] \displaystyle \rho u^2 + p + \frac{\|\B\|^2}{8 \pi} - \frac{1}{4\pi} B_x^2 \\[2mm] \displaystyle \rho u v  - \frac{1}{4\pi} B_x B_y \\[2mm] \displaystyle \rho u w  - \frac{1}{4\pi} B_x B_z \\[2mm] \displaystyle \left(\rho E + p + \frac{\|\B\|^2}{8 \pi}\right) u - \frac{1}{4 \pi} B_x (\vv \cdot \B) \\[2mm] \displaystyle 0  \\[2mm] \displaystyle u B_y - v B_x \\[2mm] \displaystyle u B_z - w B_x
	\end{array} \right).
\label{eqn.MHD1Dterms}
\end{equation}
To handle multiple time scales, we resort to a flux splitting technique that permits to design different time discretizations according to the stiffness of each resulting sub-system. Following \cite{3splitMHD}, the fluxes of the MHD system \eqref{sys:MHD} are split into three contributions: i) convective fluxes $\f^c$ related to $u$; ii) pressure fluxes $\f^p$ related to $p$; iii) magnetic fluxes $\f^b$ related to $\B$.
Here, we propose the following definitions
\begin{equation}
	\f^c(\q) = \left( \begin{array}{c} \displaystyle \rho u \\[2mm] \displaystyle \rho u^2 \\[2mm] \displaystyle \rho u v  \\[2mm] \displaystyle \rho u w  \\[2mm] 0 \\[2mm] \displaystyle 0  \\[2mm] 0 \\[2mm] \displaystyle 0
	\end{array} \right), \qquad \f^p(\q) = \left( \begin{array}{c} \displaystyle 0 \\[2mm] \displaystyle p \\[2mm] \displaystyle 0 \\[2mm] \displaystyle 0 \\[2mm] \displaystyle \left(\rho E + p \right) u \\[2mm] \displaystyle 0  \\[2mm] 0 \\[2mm] 0
\end{array} \right), \qquad \f^b(\q) = \left( \begin{array}{c} \displaystyle 0 \\[2mm] \displaystyle \frac{\|\B\|^2}{8 \pi} - \frac{1}{4\pi} B_x^2 \\[2mm] \displaystyle  - \frac{1}{4\pi} B_x B_y \\[2mm] \displaystyle - \frac{1}{4\pi} B_x B_z \\[2mm] \displaystyle \frac{\|\B\|^2}{8 \pi}\, u - \frac{1}{4 \pi} B_x (\vv \cdot \B) \\[2mm] \displaystyle 0  \\[2mm] \displaystyle u B_y - v B_x \\[2mm] \displaystyle u B_z - w B_x
\end{array} \right).
	\label{eqn.fluxsplit}
\end{equation}
We remark that the above splitting is different from the formulation proposed in \cite{Toro_Vazquez12} for the Euler equations which was applied in \cite{SIMHD_Dumbser2019,3splitMHD} on the MHD system.
The difference here is that in \eqref{eqn.fluxsplit} the hydrodynamic component of the energy flux is not split into a convective and a pressure part, but is entirely assigned to the pressure fluxes $\f^p$.
In particular, the magnetic energy contribution $\rho m$ will also play a role in the fluxes $\f^p$. The novel 3-split form of the one-dimensional system \eqref{eqn.MHD1D} is then given in compact notation by
\begin{equation}
	\frac{\partial \Q}{\partial t} + \frac{\partial \f^c(\q)}{\partial x} + \frac{\partial \f^p(\q)}{\partial x} + \frac{\partial \f^b(\q)}{\partial x}= \mathbf{0}.
	\label{eqn.MHD1D_3split}
\end{equation}
In the following we detail for each sub-system the associated set of eigenvalues.
\begin{enumerate}
	\item For the convective sub-system
	\begin{subequations}
		\begin{align}
			\pdd{\Q}{t} + \pdd{\f^c(\q)}{x} = \mathbf{0} \label{eqn.fc_sub}
		\end{align}
    we find the following real eigenvalues
    \begin{equation}
        \lambda^c_{1,2,3,4} = 0, \quad \lambda^c_{5,6,7,8} = u \label{eqn.fc_eig}.
    \end{equation}
    Thus sub-system \eqref{eqn.fc_sub} is hyperbolic with characteristic speeds which only depend on the fluid flow.
	\end{subequations}

    \item For the pressure sub-system
    \begin{subequations}
    	\begin{align}
    		\pdd{\Q}{t} + \pdd{\f^p(\q)}{x} = \mathbf{0} \label{eqn.fp_sub}
    	\end{align}
        we find the following eigenvalues
        \begin{equation}
            \lambda^p_{1,2,3,4,5,6} = 0, \quad
            \lambda^p_{7,8} = \frac{1}{2}\left( u \mp \sqrt{u^2 + 4 (c^2 - (\gamma-1)(m+k+u^2))} \right). \label{eqn.fp_eig}
        \end{equation}
    \end{subequations}
    Compared to the eigenvalues of the pressure sub-system derived in \cite{SIMHD_Dumbser2019,3splitMHD}, we notice the additional contribution of the specific magnetic and kinetic energy. We remark that the eigenvalues $\lambda^p_{7,8}$ might become complex in the case of highly magnetized fluids, thus the pressure sub-system may result to be elliptic.

    \item For the magnetic sub-system
     \begin{subequations}
    	\begin{align}
    		\pdd{\Q}{t} + \pdd{\f^b(\q)}{x} = \mathbf{0}, \qquad \lambda^b_{1,2,3,4} &= 0 \label{eqn.fb_sub}
    	\end{align}
        we find the following real eigenvalues
        \begin{equation}
            \lambda^b_{1,2,3,4} = 0, \quad
            \lambda^b_{5,6} = \frac{1}{2}\left( u \mp \sqrt{u^2 + 4 \left( \frac{B_x}{\sqrt{4\pi\rho}} \right)^2} \right), \quad
            \lambda^b_{7,8} = \frac{1}{2}\left( u \mp \sqrt{u^2 + 4 \left( \frac{\|\B\|^2}{\sqrt{4\pi\rho}} \right)^2} \right). \label{eqn.fb_eig}
        \end{equation}
    \end{subequations}
    Thus the magnetic sub-system is hyperbolic.
\end{enumerate}

Looking at the eigenvalues of the three sub-systems, it is clear that for low Mach number flows, the acoustic sound speed $c$ in \eqref{eqn.fp_eig} becomes dominant compared to the material speed $u$ of the flow, hence needing an implicit discretization of the terms related to the pressure. This is also confirmed by analyzing the non-dimensional system \eqref{sys:MHD_nd} for $\Mc \to 0$. Likewise, if the magnitude of the magnetic field $\|\B\|$ increases, i.e. when $\beta \to 0$ in \eqref{sys:MHD_nd_beta}, the MHD model is dominated by the magnetic sub-system \eqref{eqn.fb_sub} with eigenvalues $\lambda^b$. As such, an implicit treatment of the terms related to the magnetic field would allow for a larger stability region of the resulting numerical method. Therefore, the possible elliptic nature of the pressure sub-system \eqref{eqn.fp_sub} does not pose any stability problem.
The numerical scheme which will be designed in this work starts from these considerations.

\section{Numerical scheme}
\label{sec.numscheme}

Our objective is to construct a numerical scheme for the ideal MHD equations \eqref{sys:MHD} which is stable independently of the parameters $\Mc$ and $\beta$ which are associated to the scaling of the pressure (and internal energy) and the strength of the magnetic field, respectively. To achieve this goal, we rely on the novel 3-split formulation \eqref{eqn.MHD1D_3split}, dividing the fluxes into explicitly and implicitly treated terms. In particular, the pressure and magnetic sub-systems with $\f^p$ and $\f^b$ are discretized implicitly, while we retain an explicit scheme for the convective sub-system with $\f^c$.

This approach simultaneously yield two main advantages: i) no degradation of accuracy occurs which would be caused for explicit schemes by an excessive due to scale dependent numerical viscosity in case of low Mach number flows, see e.g. \cite{Dellacherie1,Klein}; ii) improvement of the computational efficiency of the scheme in the low Mach number or low plasma-beta regime, since the stability of the scheme is guaranteed by a time step imposed only by the velocity of the local fluid flow $\vv$.

However, the resulting implicit system is in general nonlinear due to the non-linearity of the implicitly treated flux terms. Moreover, this implies a huge computational overhead since nonlinear solvers would be required to solve large coupled systems. To overcome this problem, the non-linearity can be treated by adopting a fixed point iteration algorithm, that has been successfully applied to low Mach number flows for the Euler and MHD systems \cite{Dumbser_Casulli16,SIMHD_Dumbser2019}. If the magnetic sub-system is discretized implicitly as done in \cite{3splitMHD}, the nonlinear coupling is even stronger, thus a fixed point method combined with a nested Newton solver \cite{NestedNewton} has to be employed, which makes the resulting method quite demanding in terms of implementation and computational efficiency.

Therefore, our goal is to give a numerical scheme which only includes linear systems to be solved reducing the computational costs without modifying the model equations \eqref{sys:MHD}. Furthermore, a 3-split form as given in \eqref{eqn.MHD1D_3split} and proposed in \cite{3splitMHD} is not suitable for the class of IMEX Runge-Kutta methods \cite{PR_IMEX,Avgerinos2019} that allows only the splitting according to two time scales, namely a slow and a fast scale.
Hence dealing with a 2-split form of the governing equations is favourable and we will detail how to resort to a semi-implicit IMEX discretization of the fluxes, following the time integrators proposed in \cite{BosFil2016} yielding an overall second-order time discretization.

Besides the treatment in time, the MHD equations require to be compliant with the divergence-free property of the magnetic field \eqref{eqn.divBFree}. A widespread technique consists in the use of staggered grids \cite{yee1966}, since it has been shown that an electric field defined on a mesh that is dual with respect to the mesh of the magnetic field allows the divergence-free condition to be preserved up to machine accuracy.
Nevertheless, the staggering of the meshes can become cumbersome, especially in higher dimensions.
Consequently, we prefer to adopt a collocated approach utilizing a reformulation of the governing equations in terms of the magnetic potential $\A$ \cite{helzel2011,helzel2013}. Since the magnetic field $\B$ is the curl of $\A$, a mimetic finite difference div-curl operator will then be designed following the general approach introduced in \cite{divcurlB}.

Summarizing, we will first derive the time semi-discrete numerical scheme in one space dimension where the focus lies on the semi-implicit formulation of the flux splitting. Since the divergence-free constraint reduces to $\partial_x(B_x) = 0$, we solve for the magnetic field directly. Secondly, the time semi-discrete numerical scheme in multiple space dimensions is considered, where the magnetic field $\B$ is replaced by the curl of the magnetic potential, i.e. $\B = \nabla \times \A$. Thirdly, the extension to second order time accuracy is carried out relying on semi-implicit IMEX time marching schemes.
Then, the details of the spatial discretization are provided, including the structure-preserving div-curl operator needed to respect the solenoidal property of the magnetic field at the discrete level upon which the section is concluded with the description of the fully discrete numerical scheme.

\paragraph{Space and time computational domain.} The three-dimensional computational domain is defined in terms of the spatial intervals as $\Omega =[x_{\min};x_{\max}] \times [y_{\min};y_{\max}] \times [z_{\min};z_{\max}]$, and it is discretized using a Cartesian grid composed of $N_x \times N_y \times N_z$ cells. The characteristic mesh sizes are then given by
\begin{equation}
	\dx=\frac{x_{\max}-x_{\min}}{N_x}, \qquad \dy=\frac{y_{\max}-y_{\min}}{N_y}, \qquad  \dz=\frac{z_{\max}-z_{\min}}{N_z}.
\end{equation}

The time interval $T=[0; t_f]$ is discretized by time steps $\dt = t^{n+1}-t^n$ which are always subject to a CFL-type stability condition that is only based on the characteristic speeds of the convective sub-system \eqref{eqn.fc_sub}.
Thus, with $|\lambda^{c}|^n$ denoting the maximum eigenvalue of the convective sub-system \eqref{eqn.fc_eig} at the time level $t^n$, the CFL condition is given by
\begin{equation}
	\label{eq:CFL}
	\Delta t \leq \CFL \frac{(\dx \dy \dz)^{1/3}}{\max \limits_{\Omega}|\lambda^{c}|^n}.
\end{equation}

\subsection{Time semi-discrete scheme in one space dimension}\label{sec.numscheme.1d}
We consider the one-dimensional MHD system \eqref{eqn.MHD1D}-\eqref{eqn.MHD1Dterms}. The divergence free condition on the magnetic field in one dimension reduces to $\partial_x(B_x) = 0$, thus the first component of the magnetic field needs to be constant, that is $B_x=const$.

Applying a forward Euler method on the interval $[t^{n};t^{n+1}]$ to the convective sub-system \eqref{eqn.fc_sub}, we obtain the intermediate states $\q^\star$ for density and momentum as
\begin{equation}
    \q^\star = \q^n - \Delta t \, \pdd{\f^{c}(\q^{n})}{x}.
    \label{eqn.qstar1D}
\end{equation}
The remaining flux terms will be considered implicitly, leading to a coupling of the magnetic field with the momentum and energy equation. The two implicit sub-systems will be solved sequentially.
\begin{enumerate}
	\item The coupling of the magnetic field with the momentum equation provides the magnetic field at the new time level $\B^{n+1}$.
	\item The coupling of the energy with the momentum equation yields the total energy density at the new time level $(\rho E)^{n+1}$.
\end{enumerate}

The time discretization of the implicit magnetic subsystem reads
\begin{subequations}
    \label{sys:MHD_1D_Bfield}
    \begin{align}
        \displaystyle (\rho u)^{n+1} &= \rho u^\star - \Delta t \pdd{}{x} \left(p^n + \frac{\B^n\cdot \B^{n+1}}{8 \pi} - \frac{1}{4\pi} B_x^2\right), \label{eq:MomentumX_implicit}\\
        \displaystyle (\rho v)^{n+1} &=  \rho v^\star - \Delta t \pdd{}{x} \left( - \frac{1}{4\pi} B_x B_y^{n+1}\right),  \label{eq:MomentumY_implicit}\\
        \displaystyle (\rho w)^{n+1} &= \rho w^\star - \Delta t \pdd{}{x} \left( - \frac{1}{4\pi} B_x B_z^{n+1}\right),  \label{eq:MomentumZ_implicit}\\
        B_y^{n+1} &= B_y^n - \Delta t\pdd{}{x} \left(\frac{(\rho u)^{n+1}}{\rho^{n+1}} B_y^n - \frac{(\rho v)^{n+1}}{\rho^{n+1}} B_x \right),  \label{eq:By_implicit} \\
        B_z^{n+1} &= B_z^n - \Delta t \pdd{}{x} \left(\frac{(\rho u)^{n+1}}{\rho^{n+1}} B_z^n - \frac{(\rho w)^{n+1}}{\rho^{n+1}} B_x \right). \label{eq:Bz_implicit}
    \end{align}
\end{subequations}
Note that the density at the new time $\rho^{n+1}$ is known after the explicit step \eqref{eqn.qstar1D}, and $B_x$ is constant.
In addition, the pressure contribution in the momentum equation is treated explicitly inspired by the approach forwarded in \cite{BT2023} for an implicit treatment of the viscous terms in the compressible Navier-Stokes equations. We remark that we could even not consider the pressure contribution at this stage, since here the momentum equation is needed only to provide a flux for the magnetic sub-system, and we are not solving for the momentum at the new time level.
Formal substitution of the evolution equations  \eqref{eq:MomentumX_implicit}-\eqref{eq:MomentumZ_implicit} into the equation for the magnetic field components \eqref{eq:By_implicit}-\eqref{eq:Bz_implicit} yields the following sub-system for the unknowns $B_y^{n+1}$ and $B_z^{n+1}$:
\begin{subequations}
    \label{sys:MHD_1D_Bfield_elliptic}
    \begin{eqnarray}
        && B_y^{n+1} = B_y^\star + \Delta t^2\pdd{}{x} \left(\frac{B_y^n}{\rho^{n+1}}\pdd{}{x} \left(\frac{B_y^n B_y^{n+1} + B_z^n B_z^{n+1}}{8 \pi}\right) + \frac{B_x}{\rho^{n+1}} \pdd{}{x} \left(\frac{B_x B_y^{n+1}}{4\pi} \right)\right), \\
        && B_z^{n+1} = B_z^\star + \Delta t^2 \pdd{}{x} \left(\frac{B_z^n}{\rho^{n+1}}\pdd{}{x} \left(\frac{B_y^n B_y^{n+1} + B_z^n B_z^{n+1}}{8 \pi}\right) + \frac{B_x}{\rho^{n+1}} \pdd{}{x} \left(\frac{B_x B_z^{n+1}}{4\pi} \right)\right),
    \end{eqnarray}
\end{subequations}
where
\begin{eqnarray}
    && B_y^\star = B_y^n - \Delta t\pdd{}{x} \left(\frac{B_y^n}{\rho^{n+1}}\left(\rho u^\star - \Delta t\pdd{}{x} \left(p^n - \frac{B_x^2}{8\pi} \right)\right) + \frac{B_x}{\rho^{n+1}} \rho v^\star\right),\\
    && B_z^\star = B_z^n - \Delta t\pdd{}{x} \left(\frac{B_z^n}{\rho^{n+1}}\left(\rho u^\star - \Delta t\pdd{}{x} \left(p^n - \frac{B_x^2}{8\pi} \right)\right) + \frac{B_x}{\rho^{n+1}} \rho w^\star\right).
\end{eqnarray}
Thanks to the semi-implicit discretization of the term $\|\B \|^2 \approx \B^n\cdot \B^{n+1}$ in the $x-$momentum equation \eqref{eq:MomentumX_implicit}, the above system is linear.

Finally, we can write the implicit system for the energy. Since the magnetic field at the new time level $\B^{n+1}$ is known from the solution of \eqref{sys:MHD_1D_Bfield_elliptic}, the approaches introduced in implicit-explicit schemes for the Euler equations can be readily applied. Following \cite{Avgerinos2019}, the energy sub-system reads
\begin{subequations}
    \label{sys:MHD_1D_Energy_subsystem}
    \begin{eqnarray}
        && \displaystyle (\rho u)^{n+1} = (\rho u)^\star - \Delta t \pdd{}{x} \left(p^{n+1} + \frac{\|\B^{n+1}\|^2}{8 \pi} - \frac{1}{4\pi} B_x^2\right) \label{eq:MomentumX_implicit_E} \\
        && \displaystyle (\rho E)^{n+1} = (\rho E)^n - \Delta t \pdd{}{x} \left(\left(\rho E^n + p^n\right) \frac{(\rho u)^{n+1}}{\rho^{n+1}} + \frac{\|\B^{n+1}\|^2}{8 \pi} u^n - \frac{1}{4 \pi} B_x (\vv^n \cdot \B^{n+1})\right). \label{eq:E_implicit_E}
    \end{eqnarray}
\end{subequations}
Let us notice that a semi-implicit scheme has again been applied to the split fluxes \eqref{eqn.fluxsplit}, in which the chosen discretization for the terms $\left(\rho E^n + p^n\right) u^{n+1}$ in $\f^p$ and $\frac{\|\B^{n+1}\|^2}{8 \pi} u^n$ in $\f^b$ allows to retrieve a linear sub-system. Indeed, defining
\begin{equation}
	p^{n+1} = (\rho E^{n+1} - \rho k^n - \rho m^{n+1})(\gamma - 1),
	\label{eqn.pnew}
\end{equation}
and inserting the $x-$momentum equation \eqref{eq:MomentumX_implicit_E} into \eqref{eq:E_implicit_E}, the following elliptic equation on the total energy is obtained
\begin{eqnarray}
    \label{eq:Energy_implicit}
    && \displaystyle (\rho E)^{n+1} = (\rho E)^\star + (\gamma - 1)\Delta t^2\pdd{}{x}  \left(\frac{\rho E^n + p^n}{\rho^{n+1}} \, \pdd{}{x} \left( (\rho E)^{n+1} \right) \right),
\end{eqnarray}
where
\begin{align}
    \rho E^\star &= \rho E^n - \Delta t \pdd{}{x} \left(\frac{\rho E^n + p^n}{\rho^{n+1}}(\rho u)^{\star\star} + \frac{\|\B^{n+1}\|^2}{8 \pi} u^n - \frac{1}{4 \pi} B_x (\vv^n \cdot \B^{n+1})\right),\\
    (\rho u)^{\star\star} &= (\rho u)^\star - \Delta t\pdd{}{x} \left( - (\gamma - 1)\left(\rho^n k^n + \rho^{n+1} \frac{\|\B^{n+1}\|^2}{8\pi}\right) + \frac{\|\B^{n+1}\|^2}{8 \pi} - \frac{1}{4\pi} B_x^2\right). \label{eqn.ustar2}
\end{align}
The implicit energy equation \eqref{eq:Energy_implicit} is thus a linear problem for the scalar unknown $(\rho E)^{n+1}$. This is because of: i) the choice of the ideal gas EOS \eqref{eqn.EOS}; ii) the semi-implicit time discretization in the pressure and magnetic fluxes $\f^p$ and $\f^b$ of \eqref{eqn.fluxsplit}.

Once, the total energy at the new time level is obtained, the $x$-component of the momentum can be updated according to \eqref{eq:MomentumX_implicit_E} and the $y$- and $z$-components according to \eqref{eq:MomentumY_implicit} and \eqref{eq:MomentumZ_implicit}, respectively.

\subsection{Time semi-discrete scheme in multiple space dimensions}

In multiple space dimension, the preservation of the divergence-free property of the magnetic field at the discrete level is extremely important to ensure physically admissible solutions. According to \cite{Brackbill1980}, violating the $\nabla \cdot \B = 0$ constraint leads to non-physical phenomena, e.g spurious numerical monopoles may be generated or magnetic
field lines may break into open lines. To guarantee that condition \eqref{eqn.divBFree} is obeyed, instead of solving for $\B$ directly, as done in the one-dimensional case, we use the magnetic vector potential $\A$ which generates the magnetic field by its curl. Thus, since $\nabla \cdot \B = \nabla \cdot (\nabla \times \A) = 0$, the objective of the numerical scheme is to preserve the div-curl property exactly up to machine precision.
The details will be discussed later in Section \ref{ssec.discreteOp}.
Here, we provide the semi-discrete scheme in time for the ideal MHD system in multiple space dimensions equipped with $\A$ instead of $\B$.

The system of equations is non-conservative for the time evolution of $\A$, and it is given by
\begin{equation}
	\label{sys:MHD}
	\pdd{}{t} \left( \begin{array}{c}
		\displaystyle \rho \\[2mm] \displaystyle \rho \vv \\[2mm] \displaystyle \rho E \\[2mm] \displaystyle \A
	\end{array} \right) + \nabla \cdot \left( \begin{array}{c}
		\displaystyle \rho \vv \\[2mm] \displaystyle \rho \vv \otimes \vv + \left(p + \frac{\| \nabla \times \A\|^2}{8 \pi}\right) \mathbb{I} - \frac{1}{4\pi}  (\nabla \times \A) \otimes ( \nabla \times \A) \\[2mm] \displaystyle \left(\rho E + p + \frac{\| \nabla \times \A\|^2}{8 \pi}\right) \vv - \frac{1}{4 \pi} ( \nabla \times \A) (\vv \cdot ( \nabla \times \A)) \\[2mm] \displaystyle \mathbf{0}
	\end{array} \right) + \left( \begin{array}{c} \displaystyle 0 \\[2mm] \displaystyle \mathbf{0} \\[2mm] \displaystyle 0 \\[2mm] \displaystyle \left(\nabla \times \A\right) \times \vv \end{array} \right) = \mathbf{0},
\end{equation}
with the identity
\begin{equation}
	\B - \nabla \times \A = 0.
	\label{eqn.Bdef}
\end{equation}

The structure of the time semi-discrete scheme is the same as in the one-dimensional case, but the magnetic sub-system is now written in terms of $\A$. Consequently, we focus only on the magnetic sub-system where $(\rho \vv)^\star$ and $\rho^{n+1}$ are known from the convective step \eqref{eqn.qstar1D}.
Then, assuming that $\B^n = \nabla \times \A^n$ is fulfilled at the current time step, the magnetic sub-system writes
\begin{subequations}
    \label{eqn:MHD_magnetic_A}
    \begin{eqnarray}
        && \displaystyle (\rho \vv)^{n+1} = (\rho \vv)^\star - \Delta t\nabla \cdot \left(\left(p^n + \frac{\B^n\cdot (\nabla \times \A^{n+1})}{8 \pi}\right) \mathbb{I} - \frac{\B^n \otimes ( \nabla \times \A^{n+1})}{4\pi}  \right), \label{eqn.momA} \\
        && \A^{n+1} = \A^{n}- \Delta t ~\B^n \times \frac{(\rho \vv)^{n+1}}{\rho^{n+1}}. \label{eq:Implicit_A}
    \end{eqnarray}
\end{subequations}
Analogously to the one-dimensional case, plugging the momentum equations \eqref{eqn.momA} into \eqref{eq:Implicit_A} yields
\begin{eqnarray}
    && \A^{n+1} = \A^\star+ \Delta t~\B^n \times \frac{\Delta t}{\rho^{n+1}}\nabla \cdot \left(\left(\frac{\B^n\cdot (\nabla \times \A^{n+1})}{8 \pi}\right) \mathbb{I} - \frac{\B^n \otimes ( \nabla \times \A^{n+1})}{4\pi}  \right)
    \label{eqn.A_implicit}
\end{eqnarray}
where
\begin{eqnarray}
    && \A^\star = \A^{n}- \Delta t~ \B^n \times \frac{1}{\rho^{n+1}}\left((\rho\vv)^\star - \Delta t\nabla p^n\right).
    \label{eqn.Astar_implicit}
\end{eqnarray}
We remark that a semi-implicit discretization has been employed for the terms $\| \nabla \times \A\|^2 \approx (\nabla \times \A^n) \cdot (\nabla \times \A^{n+1})=\B^n \cdot (\nabla \times \A^{n+1})$ and $(\nabla \times \A) \otimes ( \nabla \times \A) \approx (\nabla \times \A^n) \otimes ( \nabla \times \A^{n+1})=\B^n \otimes ( \nabla \times \A^{n+1})$. The magnetic potential at the new time level permits to easily evaluate the new magnetic field as $\B^{n+1} = \nabla \times \A^{n+1}$.

Analogously to the one-dimensional case, we obtain the implicit energy update \eqref{eq:Energy_implicit} in a multi-dimensional setting, which remains a linear system with the scalar unknown $(\rho E)^{n+1}$. After solving the energy sub-system, the momentum equations can be updated according to
\begin{equation}
    (\rho \vv)^{n+1} = (\rho\vv)^\star - \Delta t\nabla \cdot \left(\left(p^{n+1} + \frac{\|\B^{n+1}\|^2}{8 \pi}\right) \mathbb{I} - \frac{\B^{n+1} \otimes \B^{n+1}}{4\pi}  \right),
    \label{eqn.momNew}
\end{equation}
where the new magnetic field $\B^{n+1}$ is used, and the new pressure $p^{n+1}$ is computed by \eqref{eqn.pnew}.

We remark at this point, that the scheme is asymptotic preserving (AP) for $\Ma =1$ and $\Mc$ tending to zero for well-prepared initial data towards the incompressible Euler equations.
For a proof see \cite{LukPesTho2024}, where the AP property is proven in the context of an IMEX finite volume scheme for single-temperature two-fluid all speed flows. The numerical scheme for the Euler sub-system coincides with the scheme presented in this work, thus the proof can be applied directly on the time semi-discrete scheme presented in this section for $\Ma = \mathcal{O}(1)$.

\subsection{Second order time discretization}
The semi-discrete scheme presented in the previous section is first order accurate in time. To achieve second order of accuracy, we rely on the class of semi-implicit IMplicit-EXplicit (IMEX) time integrators proposed in \cite{BosFil2016}. These Runge-Kutta methods perfectly fit to our semi-discrete scheme since we have applied a time linearization to the 3-split form of the fluxes \eqref{eqn.fluxsplit}. Such linearization is typically referred to as semi-implicit technique.

Without loss of generality, let us consider the one-dimensional system \eqref{eqn.MHD1D}-\eqref{eqn.MHD1Dterms} and the associated time semi-discrete scheme presented in Section \ref{sec.numscheme.1d}. According to \cite{BosFil2016}, the governing equations can be written in the form of an autonomous system
\begin{equation}
	\frac{\partial \q}{\partial t} = \mathcal{H}(\q_E(t),\q_I(t)),
	\label{eqn.qSI}
\end{equation}
with initial condition $\q_0=\q(t=0)$. The function $\mathcal{H}$ represents the spatial approximation of the flux terms \eqref{eqn.MHD1Dterms} that will be detailed later, and it is a function of two arguments, namely the explicit and the implicit state vector $\q_E$ and $\q_I$, respectively. To integrate the system \eqref{eqn.qSI} in time, Implicit-Explicit Runge-Kutta schemes are used, see e.g. \cite{PR_IMEX} for Stiffly Accurate (SA) IMEX-RK methods widely used for the construction of multi-scale schemes.
 An IMEX-RK scheme is described by a double Butcher tableau
\begin{equation}
	\begin{array}{c|c}
		\tilde{c} & \tilde{A} \\ \hline & \tilde{b}^\top
	\end{array} \qquad
	\begin{array}{c|c}
		c & A \\ \hline & b^\top
	\end{array},
\end{equation}
with the matrices $(\tilde{A},A) \in \R^{s \times s}$ and the vectors $(\tilde{c},c,\tilde{b},b) \in \R^s$. The tilde symbol refers to the explicit scheme and matrix $\tilde{A}=(\tilde{a}_{ij})$ is a lower triangular matrix with zero elements on the diagonal, while $A=({a}_{ij})$ is a triangular matrix which accounts for the implicit scheme, thus having non-zero elements on the diagonal. Here, we adopt the LSDIRK2 scheme from \cite{Avgerinos2019}, which is represented by the following tableaux
\begin{equation}
	\begin{array}{c|cc}
		0         & 0 &  0 \\
		\tilde{c} & \tilde{c} & 0 \\ \hline
			      & 1-\alpha & \alpha
	\end{array} \qquad
	\begin{array}{c|cc}
		\alpha & \alpha &  0 \\
		1      & 1-\alpha & \alpha \\ \hline
		& 1-\alpha & \alpha
	\end{array},
	\label{eqn.LSDIRK2}
\end{equation}
with $\alpha=1-1/\sqrt{2}$ and $\tilde{c}=1/(2\alpha)$. Since $\tilde{b}=b$ and $b$ coincides with the last stage in $A$ in \eqref{eqn.LSDIRK2}, the LSDIRK2 scheme is SA which is a crucial property for assuring asymptotic consistency and accuracy of the scheme, see also \cite{BosRus2019}. To obtain a semi-implicit IMEX-RK method, let us first set $\Q_E^n=\Q^n$ and $\Q_I^n = \Q^n$, then the stage fluxes for $i = 1, \ldots, s$ are calculated as
\begin{subequations}
	\begin{align}
		\Q_E^i &= \Q_E^n + \dt \sum \limits_{j=1}^{i-1} \tilde{a}_{ij} k_j, \qquad 2 \leq i \leq s, \label{eq.QE} \\[0.3pt]
		\tilde{\Q}_I^i &= \Q_E^n + \dt \sum \limits_{j=1}^{i-1} a_{ij} k_j, \qquad 2 \leq i \leq s, \label{eq.QI}  \\[0.3pt]
		k_i &= \mathcal{H} \left( \Q_E^i, \tilde{\Q}_I^i + \dt \, a_{ii} \, k_i \right), \qquad 1 \leq i \leq s. \label{eq.k}
	\end{align}
\end{subequations}
Note that the evaluation of the stage fluxes \eqref{eq.k} requires the solution of an implicit system, that corresponds to the implicit magnetic and energy sub-systems given by \eqref{sys:MHD_1D_Bfield_elliptic} and \eqref{eq:Energy_implicit}, respectively.

Finally, the numerical solution is updated with $\Q^{n+1}=\Q^{s}$, i.e. the solution coincides with the last Runge-Kutta stage because of the Stiffly Accurate property of the chosen IMEX scheme \eqref{eqn.LSDIRK2}.

\subsection{Discrete spatial operators} \label{ssec.discreteOp}
A triple index $\ijk$ referred to each space direction allows a cell to be uniquely identified, whereas half indexes are adopted to label cell interfaces in $x-$, $y-$ and $z-$direction, namely $(\ihp)$, $(\jhp)$ and $(\khp)$, respectively. The cell barycenter has coordinates $(x_i,y_j,z_k)$, and $|\omega^{\ijk}|=\dx \dy \dz$ is the cell volume. Collocated grids are adopted for the spatial discretization of the numerical solution, thus the discrete representation of a generic scalar quantity $q(\xx)$ at a fixed time is expressed in each cell in the finite volume sense as the cell average
\begin{equation}
	\label{eqn.qFV}
	q_{\ijk} = \frac{1}{|\omega_{\ijk}|} \, \int \limits_{\omega_{\ijk}} q(\xx) \, d\xx,
\end{equation}
Likewise, we denote by $\mathbf{q}(\xx)=(q_x,q_y,q_z)$ a generic vector, whose components undergo the same cell centered discretization \eqref{eqn.qFV}, hence implying $\mathbf{q}_{\ijk}=(q_{x,\ijk},q_{y,\ijk},q_{z,\ijk})$.

We introduce a set of discrete spatial operators that will be employed in the numerical scheme.

\paragraph{Numerical flux operator $\Fop(f(q))$.} To compute numerical fluxes, we rely on the numerical flux operator $\Fop(f(q))$ which is defined in the $x-$direction by
\begin{equation}
	\Fop_x(f(q)) = \frac{\mathcal{F}_{\ihp}(f(q))-\mathcal{F}_{\ihm}(f(q))}{\dx},
	\label{eqn.Fop}
\end{equation}
where $\mathcal{F}(f(q))$ is a Rusanov--type numerical flux function which reads
\begin{equation}
	\mathcal{F}_{\ihp}(q) = \frac{1}{2} \left( f(q_{\ihp}^+) + f(q_{\ihp}^-) \right) -\frac{1}{2} \alpha_{\ihp} \left( q_{\ihp}^+ - q_{\ihp}^- \right),
	\label{eqn.Fop_flux}
\end{equation}
with $f(q)$ representing the physical flux related to variable $q$, and $q_{\ihp}^\pm$ denoting the boundary extrapolated values of the quantity $q$ at the cell interface $(\ihp)$. The numerical dissipation is computed as
\begin{equation}
	\alpha_{\ihp} = \max \left( |u_{\ip}|,|u_{\ijk}| \right),
\end{equation}
thus it is only proportional to the material speed $u$ in $x-$direction. We expect that for low Mach number flows this amount of numerical viscosity will be enough to stabilize the scheme, while for high Mach numbers the speed of sound is bounded by the material speed. The above operator can be easily extended to the other spatial directions, namely
\begin{equation}
	\Fop(f(q))=\Fop_x(f(q))+\Fop_y(f(q))+\Fop_z(f(q)).
\end{equation}
A first order accurate numerical flux is obtained by choosing $q_{\ihp}^+=q_{\ip}$ and $q_{\ihp}^-=q_{\ijk}$ in \eqref{eqn.Fop_flux}. To attain second order of accuracy, a TVD reconstruction is carried out for the quantity $q(\xx)$ in each spatial direction, hence obtaining a piecewise polynomial representation $r(\xx)$ which, in $x-$direction for cell $\omega_{\ijk}$, writes
\begin{equation}
	r_{\ijk}(x) = c_0 \, q_{\ijk} + c_1 \,  (x-x_i),
\end{equation}
with the expansion coefficients computed using a minmod limiter \cite{ToroBook} as
\begin{equation}
	c_0 = 1, \qquad c_1 = \text{minmod} \left( \frac{q_{\ip}-q_{\ijk}}{\dx}, \, \frac{q_{\ijk}-q_{\im}}{\dx} \right).
\end{equation}
A second order numerical flux is then obtained by setting $q_{\ihp}^+=r_{\ip}(x_{i+\frac{1}{2}})$ and $q_{\ihp}^-=r_{\ijk}(x_{i-\frac{1}{2}})$ in \eqref{eqn.Fop_flux}.

\paragraph{Central flux operator $\Bop(f(q))$.} The operator $\Bop(f(q))$ is the same as $\Fop(f(q))$ but zero numerical dissipation is considered, i.e. $\alpha_{\ihp}=0$ in \eqref{eqn.Fop_flux}, hence obtaining a second order central flux operator.

\paragraph{Differential operator $\Hop(h,q)$.} This operator applies to two generic quantities $h(\xx)$ and $q(\xx)$. In $x-$direction, it computes a second derivative with second order of accuracy of the form $\partial_x(h \, \partial_x q)$ by
\begin{equation}
	\Hop_x(h,q) = \frac{1}{\dx} \left( h_{\ihp} \, \frac{q_{\ip}-q_{\ijk}}{\dx} - h_{\ihp} \, \frac{q_{\ijk}-q_{\im}}{\dx} \right),
	\label{eqn.Hop}
\end{equation}
with the interface quantities $h_{i\pm1/2 j k}$ given by the arithmetic average
\begin{equation}
	h_{\ihp} = \frac{1}{2} (h_{\ip}+h_{\ijk}), \qquad h_{\ihm} = \frac{1}{2} (h_{\ijk}+h_{\im}).
\end{equation}
The multidimensional extension writes $\Hop(h,q)=\Hop_x(h,q)+\Hop_y(h,q)+\Hop_z(h,q)$.

\paragraph{Gradient operator $\Gop(q)$.} The gradient operator is computed at the aid of a second order central finite difference discretization:
\begin{equation}
	\Gop(q) = \left( \begin{array}{c}
		\Gop_x(q)\\[8pt]
		\Gop_y(q) \\[8pt]
		\Gop_z(q)
	\end{array} \right) = \left( \begin{array}{c}
	\displaystyle	\frac{q_{\ip} - q_{\im}}{2\dx}\\[8pt]
	\displaystyle	\frac{q_{\jp} - q_{\jm}}{2\dy} \\[8pt]
	\displaystyle	\frac{q_{\kp} - q_{\km}}{2\dz}
	\end{array} \right).
	\label{eqn.Gop}
\end{equation}

\paragraph{Curl operator $\Cop(\mathbf{q})$.} The curl operator, directly applied to the vector $\mathbf{q}$, follows from the computation of the gradient operator. Indeed, it simply writes
\begin{equation}
	\Cop(\mathbf{q}) = \left( \begin{array}{c}
		\Cop_x(\mathbf{q})\\[8pt]
		\Cop_y(\mathbf{q}) \\[8pt]
		\Cop_z(\mathbf{q})
	\end{array} \right) = \left( \begin{array}{c}
		\Gop_x(q_z) - \Gop_z(q_y) \\[8pt]
		\Gop_z(q_x) - \Gop_x(q_z) \\[8pt]
		\Gop_x(q_y) - \Gop_y(q_x)
	\end{array} \right).
	\label{eqn.Cop}
\end{equation}

\paragraph{Divergence operator $\Dop(\mathbf{q})$.} The divergence operator also makes use of the gradient operator \eqref{eqn.Gop} as follows
\begin{equation}
	\Dop(\mathbf{q}) = \Gop_x(q_x) + \Gop_y(q_y) + \Gop_z(q_z).
	\label{eqn.Dop}
\end{equation}
It is easy to check that the chosen divergence and curl operators are \textit{structure preserving}, meaning that they satisfy the algebraic identity $\nabla \cdot (\nabla \times \mathbf{q})=0$ at the discrete level, namely
\begin{equation}
	\Dop(\Cop(\mathbf{q})) = \Gop_x(\Cop_x(\mathbf{q})) + \Gop_y(\Cop_y(\mathbf{q})) + \Gop_z(\Cop_z(\mathbf{q})) = 0.
	\label{eqn.divcurl}
\end{equation}
This relation is crucial for ensuring the solenoidal property of the magnetic field at the discrete level. A proof of the div-curl relation \eqref{eqn.divcurl} can be found in Appendix A of \cite{divcurlB}.

\paragraph{Laplace operator $\Lop(\kappa,q)$.} This operator approximates the Laplacian of a scalar function with a second order central finite difference discretization for a space-dependent diffusion coefficient $\kappa$ as
\begin{equation}
	\Lop(\nu,q) = \left( \begin{array}{c}
		\displaystyle \frac{1}{\dx^2} \left( \kappa_{\ihp}(q_{\ip}-q_{\ijk}) - \kappa_{\ihm}(q_{\ijk}-q_{\im}) \right) \\[8pt]
		\displaystyle\frac{1}{\dy^2} \left( \kappa_{\jhp}(q_{\jp}-q_{\ijk}) - \kappa_{\jhm}(q_{\ijk}-q_{\jm}) \right) \\[8pt]
		\displaystyle \frac{1}{\dz^2} \left( \kappa_{\khp}(q_{\kp}-q_{\ijk}) - \kappa_{\khm}(q_{\ijk}-q_{\km}) \right)
	\end{array} \right).
	\label{eqn.Lop}
\end{equation}

\subsection{Fully discrete scheme in multiple space dimensions}
A single time step of the new semi-implicit method from time $t^n$ to $t^{n+1}$ in a fully discrete set-up is made of the following steps.

\begin{enumerate}
	\item Solve the explicit convective sub-system for density and momentum as
	\begin{subequations}
		\begin{align}
		\rho^{\star}_{\ijk} &= \rho^{n}_{\ijk} - \dt \, \left( \, \Fop_x((\rho u)^{n}) + \Fop_y((\rho v)^{n}) + \Fop_z((\rho w)^{n}) \, \right), \\
		(\rho u)^{\star}_{\ijk} &= (\rho u)^{n}_{\ijk} - \dt \, \left( \, \Fop_x((\rho uu)^{n}) + \Fop_y((\rho uv)^{n}) + \Fop_z((\rho uw)^{n}) \, \right), \\
		(\rho v)^{\star}_{\ijk} &= (\rho v)^{n}_{\ijk} - \dt \, \left( \, \Fop_x((\rho vu)^{n}) + \Fop_y((\rho vv)^{n}) + \Fop_z((\rho vw)^{n}) \, \right), \\
		(\rho w)^{\star}_{\ijk} &= (\rho w)^{n}_{\ijk} - \dt \, \left( \, \Fop_x((\rho wu)^{n}) + \Fop_y((\rho wv)^{n}) + \Fop_z((\rho ww)^{n}) \, \right).
		\end{align}
	\end{subequations}
Let us notice that $\rho^{\star}_{\ijk}=\rho^{n+1}_{\ijk}$, thus the density at the next time level is already available. Furthermore, the explicit contribution of momentum is given by the vector $$(\rho \vv)^{\star}_{\ijk}=((\rho u)^{\star}_{\ijk},(\rho v)^{\star}_{\ijk},(\rho w)^{\star}_{\ijk})^\top.$$

    \item Solve the implicit magnetic sub-system \eqref{eqn.A_implicit} and obtain the vector potential $\A_{\ijk}^{n+1}$ from
    \begin{eqnarray}
	\A_{\ijk}^{n+1} &-& \dt \, \B_{\ijk}^{n} \times \frac{\dt}{\rho_{\ijk}^{n+1}} \, \Dop \left( \, \left(  \frac{\B_{\ijk}^{n} \cdot \Cop(\A_{\ijk}^{n+1})}{8\pi}  \right) \mathbb{I} - \frac{\B_{\ijk}^{n} \otimes \Cop(\A_{\ijk}^{n+1})}{4\pi} \, \right) \nonumber \\
	&+& \Lop(\lambda^b \dx,A_{x,\ijk}^{n+1}) + \Lop(\lambda^b \dy,A_{y,\ijk}^{n+1}) + \Lop(\lambda^b \dz,A_{z,\ijk}^{n+1})  \nonumber \\
	&=& \A_{\ijk}^{\star},
    \end{eqnarray}
with the right hand side
\begin{equation}
	\A_{\ijk}^{\star} = \A_{\ijk}^{n} - \dt \, \B_{\ijk}^n \times \frac{1}{\rho_{\ijk}^{n+1}} \left( (\rho \vv)^{\star}_{\ijk} - \dt \, \Gop(p^n_{\ijk})\right).
\end{equation}
    To compute the solution of this linear system, we use the GMRES method \cite{GMRES} up to a prescribed tolerance (typically set to $10^{-14}$). Note that, according to \cite{helzel2013}, we have added the numerical dissipation by means of the Laplace operator \eqref{eqn.Lop} applied to all components of $\A$ with the eigenvalues of the magnetic sub-system \eqref{eqn.fb_eig} as diffusion coefficients, e.g. $\kappa_x=\lambda^b \dx$. This operator is discretized implicitly, thus it does not affect the stability of the overall scheme.

    \item Compute the magnetic field and the magnetic energy at the next time level as
    \begin{equation}
    	\B_{\ijk}^{n+1} = \Cop(\A_{\ijk}^{n+1}), \qquad (\rho m)_{\ijk}^{n+1} = \frac{\| \B_{\ijk}^{n+1}\|^2}{8\pi}.
    \end{equation}
    Thanks to the structure preserving property \eqref{eqn.divcurl}, the magnetic field $\B_{\ijk}^{n+1}$ is divergence-free up to machine accuracy.

    \item After inserting the definition \eqref{eqn.pnew} in the momentum equations \eqref{eqn.ustar2}, add the implicit fluxes $\f^b$ of the magnetic sub-system to the momentum and energy equations as follows
    \begin{subequations}
    	\begin{align*}
    		(\rho u)^{\star\star}_{\ijk} &= (\rho u)^{\star}_{\ijk} - \dt \, \left( \, \Bop_x\left((2-\gamma)(\rho m)^{n+1} - \frac{B_x^{n+1}B_x^{n+1}}{4\pi} \right) + \Bop_y\left( - \frac{B_x^{n+1}B_y^{n+1}}{4\pi} \right) + \Bop_z\left( - \frac{B_x^{n+1}B_z^{n+1}}{4\pi} \right) \, \right), \\
    		(\rho v)^{\star\star}_{\ijk} &= (\rho v)^{\star}_{\ijk} - \dt \, \left( \, \Bop_x\left( - \frac{B_y^{n+1}B_x^{n+1}}{4\pi} \right) + \Bop_y\left((2-\gamma)(\rho m)^{n+1} - \frac{B_y^{n+1}B_y^{n+1}}{4\pi} \right) + \Bop_z\left( - \frac{B_y^{n+1}B_z^{n+1}}{4\pi} \right) \, \right), \\
    		(\rho w)^{\star\star}_{\ijk} &= (\rho w)^{\star}_{\ijk} - \dt \, \left( \, \Bop_x\left( - \frac{B_z^{n+1}B_x^{n+1}}{4\pi} \right) + \Bop_y\left( - \frac{B_z^{n+1}B_y^{n+1}}{4\pi} \right) + \Bop_z\left( (2-\gamma)(\rho m)^{n+1} - \frac{B_z^{n+1}B_z^{n+1}}{4\pi} \right) \, \right), \\
    		(\rho E)^{\star\star}_{\ijk} &= (\rho E)^{n}_{\ijk} -\dt \, \, \Bop_x\left( -(\gamma-1) (\rho k)^n + u^n(\rho m)^{n+1} - \frac{B_x^{n+1}}{4\pi} \vv^n \cdot \B^{n+1} \right) \nonumber \\
    		& \phantom{ = (\rho E)^{n}_{\ijk}} \, -\dt \, \, \Bop_y\left( -(\gamma-1) (\rho k)^n + v^n(\rho m)^{n+1} - \frac{B_y^{n+1}}{4\pi} \vv^n \cdot \B^{n+1} \right)  \nonumber \\
    		& \phantom{ = (\rho E)^{n}_{\ijk}} \, -\dt \, \, \Bop_z\left( -(\gamma-1) (\rho k)^n + w^n(\rho m)^{n+1} - \frac{B_z^{n+1}}{4\pi} \vv^n \cdot \B^{n+1} \right).
    	\end{align*}
    \end{subequations}

    \item Solve the elliptic equation for the total energy $(\rho E)_{\ijk}^{n+1}$ according to
    \begin{equation}
    	(\rho E)_{\ijk}^{n+1} - (\gamma-1) \, \dt^2 \, \left( \Hop_x(h^n,(\rho E)_{\ijk}^{n+1}) + \Hop_y(h^n,(\rho E)_{\ijk}^{n+1}) + \Hop_z(h^n,(\rho E)_{\ijk}^{n+1}) \right) = b_{\ijk}^{\star\star},
    \end{equation}
    with the right hand side
    \begin{equation}
    	b_{\ijk}^{\star\star} = (\rho E)_{\ijk}^{\star\star} - \dt \, \left( \, \Bop_x(h^n (\rho u)^{\star\star}) + \Bop_y(h^n (\rho v)^{\star\star}) + \Bop_z(h^n (\rho w)^{\star\star}) \, \right).
    \end{equation}
	Even in this case, the linear system is solved with the GMRES method up to a residual of the order of $10^{-14}$.

	\item Update the momentum at the next time level $(\rho \vv)_{\ijk}^{n+1}$ according to
		\begin{subequations}
			\begin{align}
				(\rho u)^{n+1}_{\ijk} &= (\rho u)^{\star\star}_{\ijk} - (\gamma-1)\dt \, \Bop_x\left((\rho E)^{n+1}\right), \\
				(\rho v)^{n+1}_{\ijk} &= (\rho v)^{\star\star}_{\ijk} - (\gamma-1)\dt \, \Bop_y\left((\rho E)^{n+1}\right), \\
				(\rho w)^{n+1}_{\ijk} &= (\rho w)^{\star\star}_{\ijk} - (\gamma-1)\dt \, \Bop_z\left((\rho E)^{n+1}\right),
		\end{align}
	\end{subequations}
    which is equivalent to \eqref{eqn.momNew} after the introduction of the new pressure \eqref{eqn.pnew}.
\end{enumerate}

\section{Numerical results} \label{sec.numtest}
In the sequel, we present a large suite of test cases for ideal MHD aiming at demonstrating the accuracy, robustness and effectiveness of our new scheme in different regimes of both the Mach and the Alfv{\'e}n number. If not otherwise specified, we set $\gamma = 3/5$ and the time step $\dt$ is computed according to the stability condition \eqref{eq:CFL} with $\text{CFL} = 0.9$, which only depends on the material waves. The novel scheme is labeled with \scheme, that stands for Semi-Implicit Finite Volume method with implicit discretization for the energy (E) and the magnetic field (B). If not stated otherwise, it is always used the second order version of the \scheme scheme in space and time.

\paragraph{Boundary conditions for the magnetic potential $\A$} In all the following test problems we always perform a linear extrapolation to determine the boundary conditions for the magnetic potential $\A$. Without loss of generality, let us consider the component $A_{x,i}$ for a cell attached to the left boundary $x_{i-1/2}=x_{\min}$ of the computational domain, i.e. $i=1$. The sought boundary state $A_{x,b}$ in $x$-direction is formally located at $x_b=x_{\min}-\dx/2$ and is computed as
\begin{equation}
	A_{x,b} = A_{x,i} - \alpha \dx \qquad \text{ with } \qquad \alpha = \frac{A_{x,i+1}-A_{x,i}}{\dx}.
	\label{eqn.Abnd}
\end{equation}
This procedure applies for all the components of $\A$ and all the spatial directions $\xx$.

\subsection{Accuracy of the numerical scheme}
We test the accuracy of the numerical scheme on a slight modification of the analytical vortex solution forwarded in \cite{Balsara2004} on the two-dimensional computational domain $\Omega=[-5;5]\times[-5;5]$.
The vortex is constructed as a perturbation $\delta \q$ on top of an unperturbed magnetohydrodynamics background flow given by the constant states $\q_0 = (\rho_0, \vv_0,p_0, \mathbf{0})$.
The perturbations of the velocity and magnetic field are given by
\begin{eqnarray}
     &&(\delta u, \delta v) = \frac{\tilde{v}}{2 \pi} \, \exp\left( \frac{1- r^2}{2}\right) \cdot (- r \sin(\theta), r \cos(\theta)), \\
     &&(\delta B_x, \delta B_y) = \frac{\tilde{B}}{2 \pi}  \, \exp\left( \frac{1- r^2}{2}\right) \cdot (- r \sin(\theta), r \cos(\theta)),
\end{eqnarray}
with the constants $\tilde{v} = \sqrt{2\pi}$ and $\tilde B = \sqrt{4 \pi}$.
The perturbation of the magnetic field is generated by the magnetic potential perturbation in $z$-direction $\delta A_z = \frac{\tilde B}{2\pi} \, \exp(1-r^2)$.
The pressure is set such that it balances the centrifugal force generated by the circular motion against the centripetal force stemming from the tension in the magnetic field lines. In radial coordinates, with the generic radial position $r = \sqrt{x^2 + y^2}$, the equation for the pressure perturbation writes
\begin{equation}
    \label{eq:ODE_vortex_dp}
    \frac{\partial \delta p}{\partial r} = \left( \frac{\rho v_\theta^2}{r} - \frac{\delta B_\theta^2}{4 \pi r} - \frac{\partial (\delta B_\theta)^2}{8 \pi r}\right),
\end{equation}
with the definitions
\begin{equation}
	\delta v_\theta = \frac{\tilde{v}}{2 \pi} \, r \, \exp\left( \frac{1- r^2}{2}\right), \qquad \delta B_\theta = \frac{\tilde{B}}{2 \pi} \, r \, \exp\left( \frac{1- r^2}{2}\right).
\end{equation}
Integration of \eqref{eq:ODE_vortex_dp} yields
\begin{equation}
    \label{eq:Vortex_dp}
    \delta p = \frac{1}{2} e^{1-r^2} \left(\frac{1}{8 \pi}\frac{\tilde{B}^2}{(2 \pi)^2} \left(1-r^2\right)-\frac{1}{2}\rho \left(\frac{\tilde{v}}{2\pi}\right)^2\right).
\end{equation}
We set the background velocity field to $\vv_0 = (1,1,0)$, thus the vortex is moving diagonally through the computational domain, and the background pressure is $p_0 = 1$. Since the novelty of the scheme is the stability under a material time step which is independent of the sound and Alfv\'en speeds, we study the experimental order of convergence (EOC) under a variation of the Alfv\'en speeds. To this end, let us note that simply adding a constant background magnetic field does not yield an analytic solution of the ideal MHD equations in this setup. Thus, since the Alfv\'en speeds depend on the inverse of the density, the Alfv\'en Mach number can be decreased by lowering the background density $\rho_0$. We underline that decreasing the density will increase the sound speed as well, and the kinetic influence on the pressure perturbation \eqref{eq:Vortex_dp} will decrease. However, the contribution of the magnetic field remains the same for all the cases considered. The maximal acoustic Mach number and and Alfv\'en speed are reported in Table \ref{tab:MachBeta_vortex} for different background densities $\rho_0 = 10^{-k}$, $k = 0, \ldots, 5$.

\begin{table}[!htbp]
	\renewcommand{\arraystretch}{1.25}
	\begin{center}
		\begin{tabular}{l|cccccc}
			& \multicolumn{6}{c}{Background density $\rho_0=10^{-k}$} \\
			& $k=0$ & $k={1}$ & $k={2}$ & $k=3$ & $k=4$ & $k=5$ \\
			\hline
			$\Mc$   & 1.606E+01 & 4.489E-01 & 1.529E-01 & 4.832E-02 & 1.529E-02 & 4.832E-03 \\
			$c_{A}$ & 1.549E-01 & 4.900E-01 & 1.549E+00 & 4.900E+00 & 1.549E+01 & 4.900E+01 \\
		\end{tabular}
	\end{center}
	\caption{MHD traveling vortex. Maximum acoustic Mach number $\Mc$ and Alfv\'en speed $c_{A}$ at the initial time.}
	\label{tab:MachBeta_vortex}
\end{table}


The simulations are run on a sequence of successively refined meshes until the final time $t_f=1$. In Table \ref{tab:Convergence}, the $L^2$ error in $\rho, u, p, B_x$ and $A_z$ are summarized. We observe that, for all regimes, the expected EOC of two is achieved.
Furthermore, the scheme performs well for small values of the density. We recall, however, that the scheme is not naturally positivity preserving.

\begin{table}[!htbp]
    \renewcommand{\arraystretch}{1.25}
    \begin{center}
    \begin{tabular}{c|c|cccccccccc}
       $\rho_0$ & $N_x=N_y$ & $\rho$ & & $u$ & & $p$ & & $B_x$ & & $A_z$ & \\\hline\hline
        \multirow{4}{*}{$10^{0}$}
        &32	 &  7.13e-02	&	   & 4.40e-02&		   & 5.86e-02	&	 &  8.26e-02&		&   4.94e-02    &  \\
        &64	 &1.63e-02	&2.13&	9.16e-03&	2.26&	1.65e-02	&1.83&	2.20e-02&	1.91&	1.25e-02	&1.99 \\
        &128	 &3.79e-03	&2.10&	2.18e-03&	2.07&	4.19e-03	&1.98&	5.66e-03&	1.96&	3.14e-03	&1.99 \\
        &256  &	9.26e-04	&2.03&	6.79e-04&	1.69&	1.05e-03	&1.99&	1.54e-03&	1.88&	8.12e-04	&1.95 \\\hline
        \multirow{4}{*}{$10^{-1}$}
        &32	  &  2.33e-03	&	   & 3.50e-02&		   & 6.84e-03	&	   & 8.21e-02	&	   & 4.84e-02	&    \\
        &64	 & 6.39e-04	&1.87&	8.55e-03&	2.04&	1.91e-03	&1.84&	2.17e-02	&1.92&	1.23e-02	&1.98\\
        &128	&1.62e-04	&1.98&	2.14e-03&	2.00&	4.90e-04	&1.96&	5.53e-03	&1.97&	3.09e-03	&1.99\\
        &256&	4.08e-05	&2.00&	6.71e-04&	1.67&	1.24e-04	&1.99&	1.51e-03	&1.87&	7.99e-04	&1.95\\\hline
        \multirow{4}{*}{$10^{-2}$}
        &32	    &4.12e-04&		   & 9.14e-02	&    	&3.90e-03 &   	&	9.05e-02 &   		&5.02e-02&     \\
        &64	&7.64e-05&	2.43&	2.45e-02	&1.90	&9.82e-04	&1.99&	2.32e-02	&1.96	&1.24e-02&	2.02\\
        &128	&1.63e-05&	2.23&	6.27e-03	&1.97	&2.69e-04	&1.87&	5.85e-03	&1.99	&3.08e-03&	2.00\\
        &256	&3.92e-06&	2.06&	1.63e-03	&1.95	&6.83e-05	&1.98&	1.58e-03	&1.89	&7.95e-04&	1.95\\\hline
        \multirow{4}{*}{$10^{-3}$}&32	  &  2.80e-04	&    &	2.07e-01&    	&	3.45e-03&    	&	9.44e-02&    	&	4.80e-02 &    \\
        &64  &4.50e-05	    &2.64&	5.88e-02&	1.81&	9.63e-04&	1.84&	2.06e-02&	2.20&	8.76e-03	&2.46\\
        &128 &5.95e-06	    &2.92&	1.72e-02&	1.77&	1.98e-04&	2.28&	4.45e-03&	2.21&	1.79e-03	&2.29\\
        &256&	8.10e-07	&2.88&	4.69e-03&	1.88&	4.74e-05&	2.06&	1.19e-03&	1.90&	4.61e-04	&1.96\\\hline
         \multirow{4}{*}{$10^{-4}$}
        &32  	&1.55e-04 &   		&1.67e-00 &   		&1.28e-02 &   		&2.60e-01 &   		&2.15e-01 &       \\
        &64 	&4.98e-05	&1.64	&2.25e-01	&2.89	&4.49e-03	&1.51	&8.01e-02	&1.70	&4.92e-02	&2.12 \\
        &128	&5.87e-06	&3.09	&3.44e-02	&2.71	&8.21e-04	&2.45	&1.55e-02	&2.37	&7.77e-03	&2.66 \\
        &256&	7.22e-07&	3.02&	9.02e-03&	1.93&	1.30e-04&	2.66&	2.45e-03&	2.66&	1.14e-03&	2.76  \\\hline
         \multirow{4}{*}{$10^{-5}$}
        &32  	&1.85e-05 &   		&8.14e-00 &   		&1.32e-02 &   		&2.85e-01 &   		&2.35e-01 &       \\
        &64 	&9.85e-06	&0.91	&8.62e-01	&3.24	&6.75e-03	&0.96	&1.24e-01	&1.20	&8.47e-02	&1.47 \\
        &128	&1.41e-06	&2.81	&8.42e-02	&3.36	&1.64e-03	&2.04	&3.04e-02	&2.03	&1.70e-02	&2.32 \\
        &256&	1.61e-07&	3.12&	2.05e-02&	2.04&	2.71e-04&	2.60&	4.85e-03&	2.65&	2.54e-03&	2.74
    \end{tabular}
    \end{center}
\caption{MHD traveling vortex. $L^2$ error and EOC for the second order \scheme scheme with varying background density $\rho_0$ at final time $t_f = 1$. }
\label{tab:Convergence}
\end{table}

\begin{figure}[!htbp]
	\begin{center}
		\begin{tabular}{c}
			\includegraphics[width=0.7\textwidth,keepaspectratio=true]{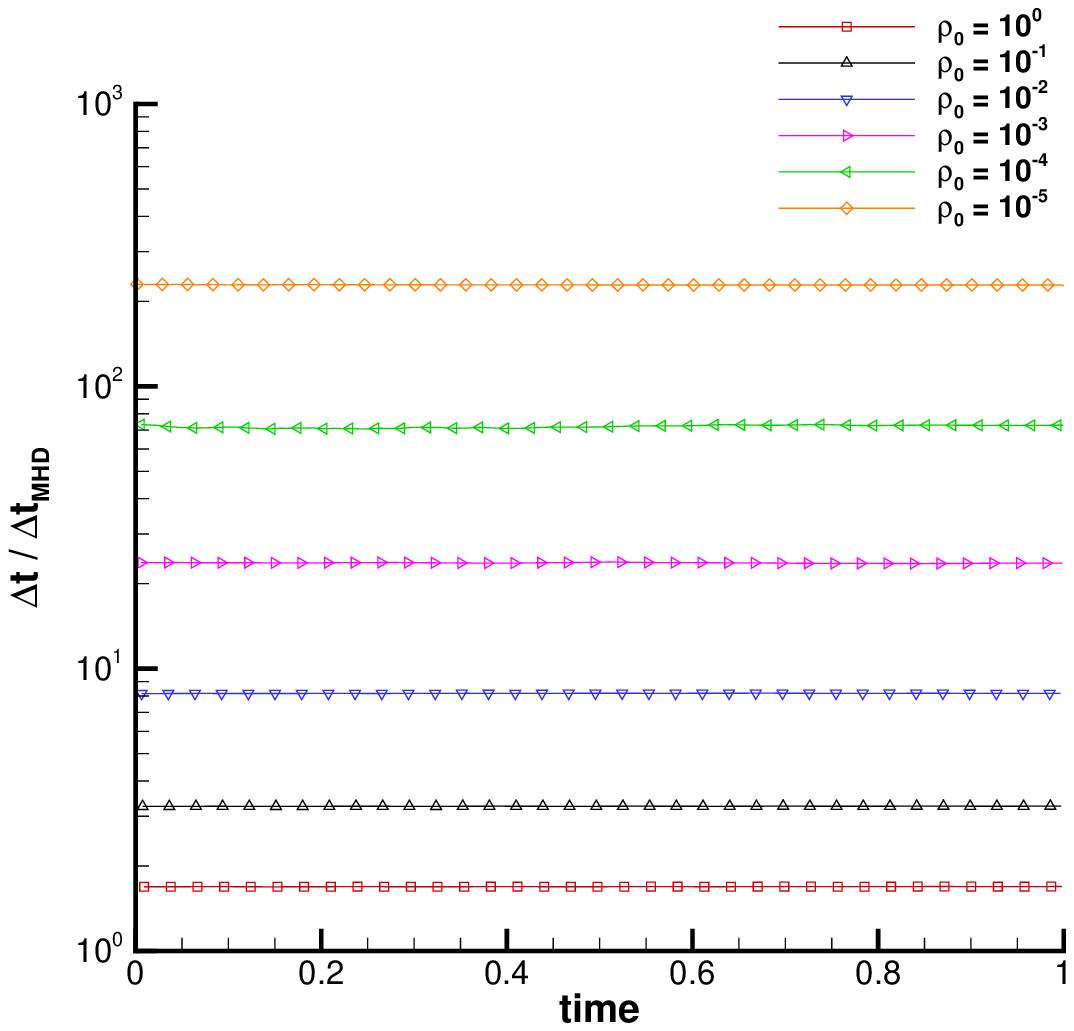}
		\end{tabular}
		\caption{MHD traveling vortex. Time evolution of the ratio between the material time step of the \scheme scheme and the magneto-sonic time step of a fully explicit finite volume scheme for different values of the background density $\rho_0$.}
		\label{fig.vortexdt}
	\end{center}
\end{figure}

To demonstrate the advantage of having a scale independent CFL stability condition, we plot in Figure \ref{fig.vortexdt} the ratio between a fully explicit time step restricted by the fastest characteristic speeds against the time step given by the CFL condition of our scheme \eqref{eq:CFL}. The \scheme scheme allows, especially for low acoustic and Alfv\'en Mach numbers, to use a time step which is several orders of magnitude larger than the one of an explicit scheme while maintaining second order accuracy in space and time.

\subsection{1D Riemann Problems}
Next, we consider a series of Riemann Problems (RP) which are standard one-dimensional test cases for ideal MHD solvers. The setup is detailed in Table \ref{tab.RPinit} and it is taken from \cite{BalDum2016}.
They have already been used to validate for instance the Osher-type scheme \cite{DumToro2011} or the HLLEM Riemann solver presented in \cite{BalDum2016}.
We compare the results obtained with our new \scheme scheme against the exact solution, that has been computed with an exact Riemann solver that has kindly been provided by Falle and Komissarov \cite{Falle2002,FalKom2001}. For an alternative exact Riemann solver of the MHD equations, see also \cite{torrilhon2003}. These test cases verify the robustness of the new \scheme scheme in the case of moderate and high Mach number flows. The computational domain is given by $\Omega=[-0.5;0.5]$ and the grid is made of 2000 cells along the $x-$direction with Dirichlet boundaries.

\begin{table}[!htbp]
	\begin{center}
		\begin{small}
			\renewcommand{\arraystretch}{1.25}
			\begin{tabular}{cr|cccccccc|c}
				Case && $\rho$ & $u$ & $v$ & $w$ & $p$ & $B_x$ & $B_y$ & $B_z$ & $t_f$ \\
				\hline\hline
				\multirow{2}{*}{RP1} &L:& 1.0 & 0.0 & 0.0 & 0.0 & 1.0 & $0.75\sqrt{4\pi}$ & +$\sqrt{4\pi}$ & 0.0 & 0.10 \\
				& R:& 0.125 & 0.0 & 0.0 & 0.0 & 0.1 & $0.75\sqrt{4\pi}$ & $-\sqrt{4\pi}$ & 0.0 &  \\\hline
				\multirow{2}{*}{RP2} &L:& 1.08 &     1.2 &   0.01 &      0.5 &    0.95 & 2.0 &    3.6 &    2.0 &  0.2 \\
				& R:& 0.9891 & $-$0.0131 & 0.0269 & 0.010037 & 0.97159 & 2.0 & 4.0244 & 2.0026 &  \\\hline
				\multirow{2}{*}{RP3} &L:& 1.7 & 0.0 & 0.0 &       0.0 & 1.7 & 3.899398 & 3.544908 &      0.0 & 0.15 \\
				& R:& 0.2 & 0.0 & 0.0 & $-$1.496891 & 0.2 & 3.899398 & 2.785898 & 2.192064 &  \\\hline
				\multirow{2}{*}{RP4} &L:& 1.0 & 0.0 & 0.0 & 0.0 & 1.0 & 1.3$\sqrt{4\pi}$ & +$\sqrt{4\pi}$ & 0.0 & 0.16 \\
				&R:& 0.4 & 0.0 & 0.0 & 0.0 & 0.4 & 1.3$\sqrt{4\pi}$ & $-\sqrt{4\pi}$ & 0.0 &   \\\hline
				\multirow{2}{*}{RP5} &L:& 0.15 &  \phantom{+}21.55 & 1.0 & 1.0 & 0.28 & 0.05 & $-$2.0 & $-$1.0 & 0.04 \\
				& R:& 0.10 & $-$26.45 & 0.0 & 0.0 & 0.10 & 0.05 & +2.0 & +1.0 &  \\\hline
				\multirow{2}{*}{RP6} & L:& 1.0 &  36.87 & $-$0.115 & $-$0.0386 & 1.0 & 4.0 & 4.0 & 1.0 & 0.03 \\
				& R:& 1.0 & $-$36.87 &    0.0 &     0.0 & 1.0 & 4.0 & 4.0 & 1.0 &  \\\hline
				\multirow{2}{*}{RP7} &L:& $1/\mu_0$ & $-$1.0 & +1.0 & $-$1.0 & 1.0 & 1.0 & $-$1.0 & 1.0 & 0.25 \\
				&R:& $1/\mu_0$ & $-$1.0 & $-$1.0 & $-$1.0 & 1.0 & 1.0 & +1.0 & 1.0 & \\
			\end{tabular}
		\end{small}
	\end{center}
\caption{Initialization of Riemann problems. Initial states left (L) and right (R) are reported as well as the final time of the simulation $t_f$. In all cases $\gamma=5/3$ and the position of the initial discontinuity is $x_d=0$, except for RP2 and RP3 where it is $x_d=-0.1$.}
	\label{tab.RPinit}
\end{table}

\begin{figure}[!htbp]
	\begin{center}
		\begin{tabular}{ccc}
			\includegraphics[width=0.33\textwidth]{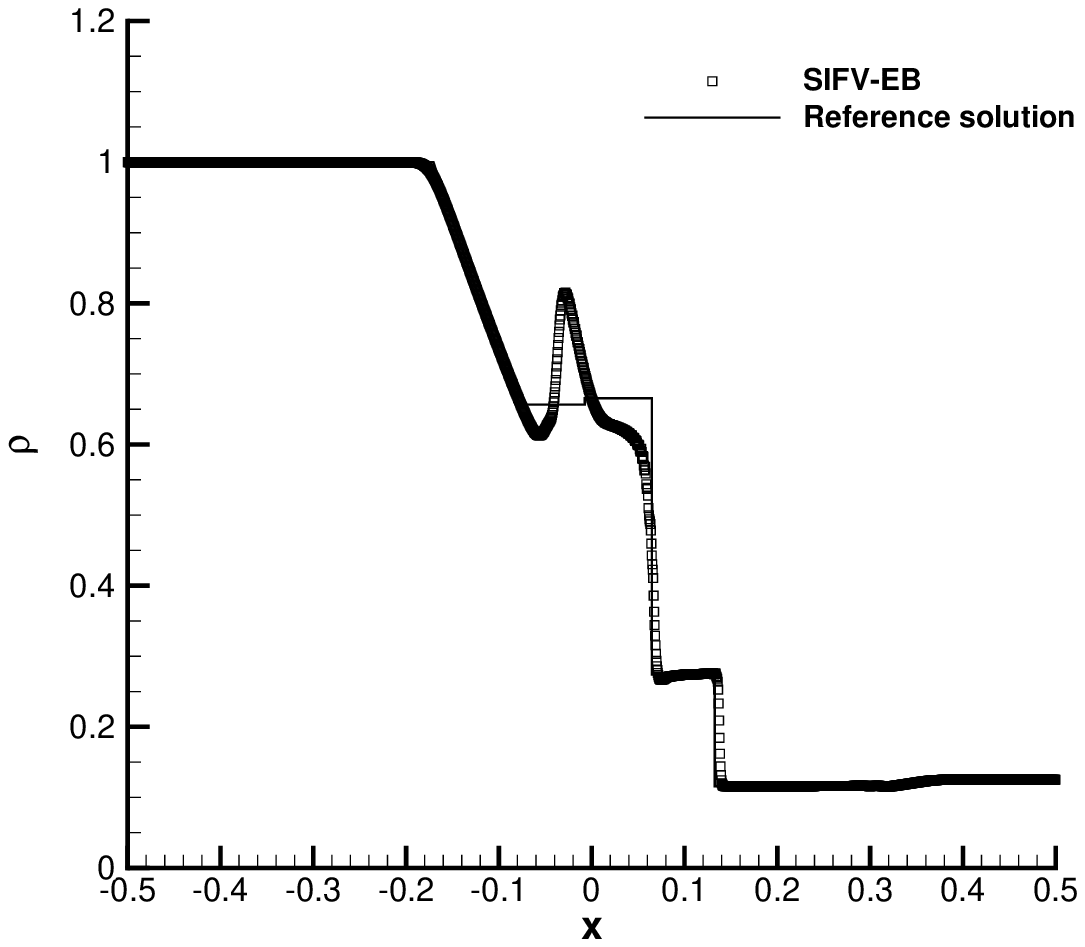}  &
			\includegraphics[width=0.33\textwidth]{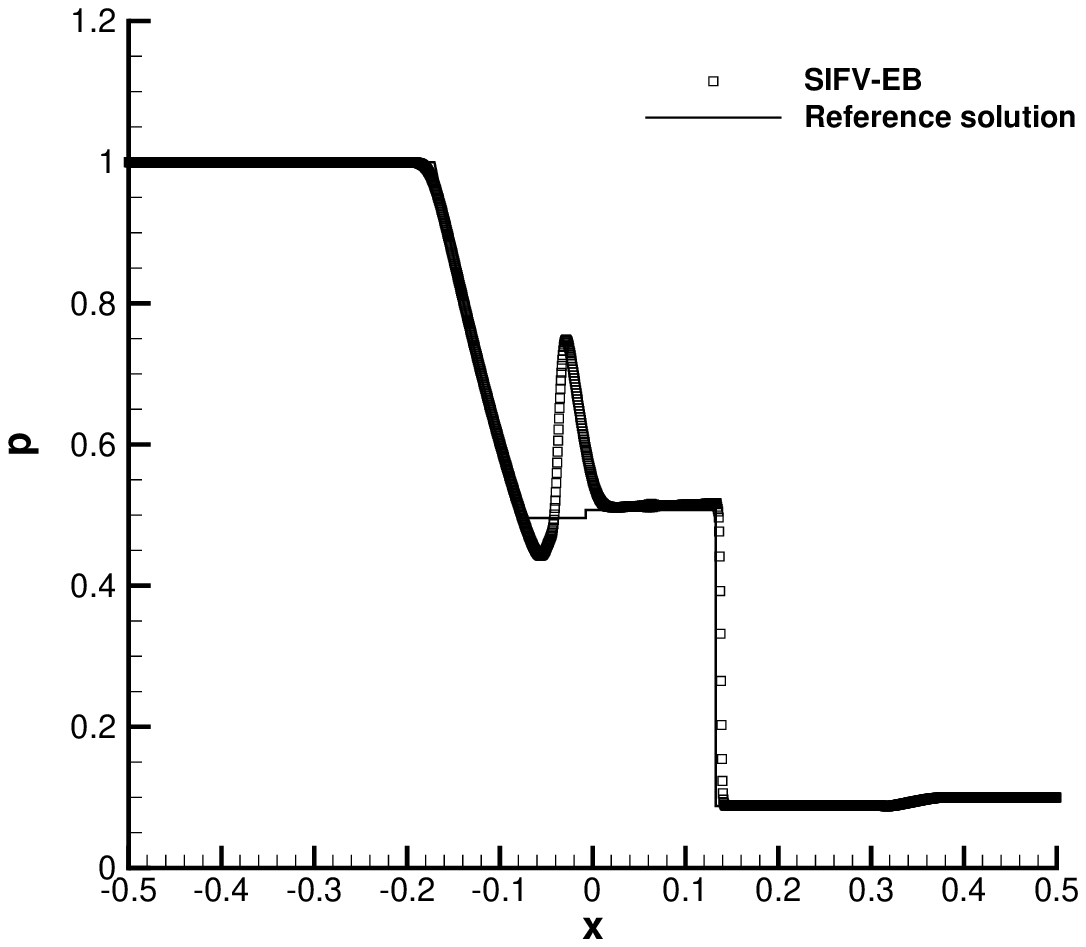} &
			\includegraphics[width=0.33\textwidth]{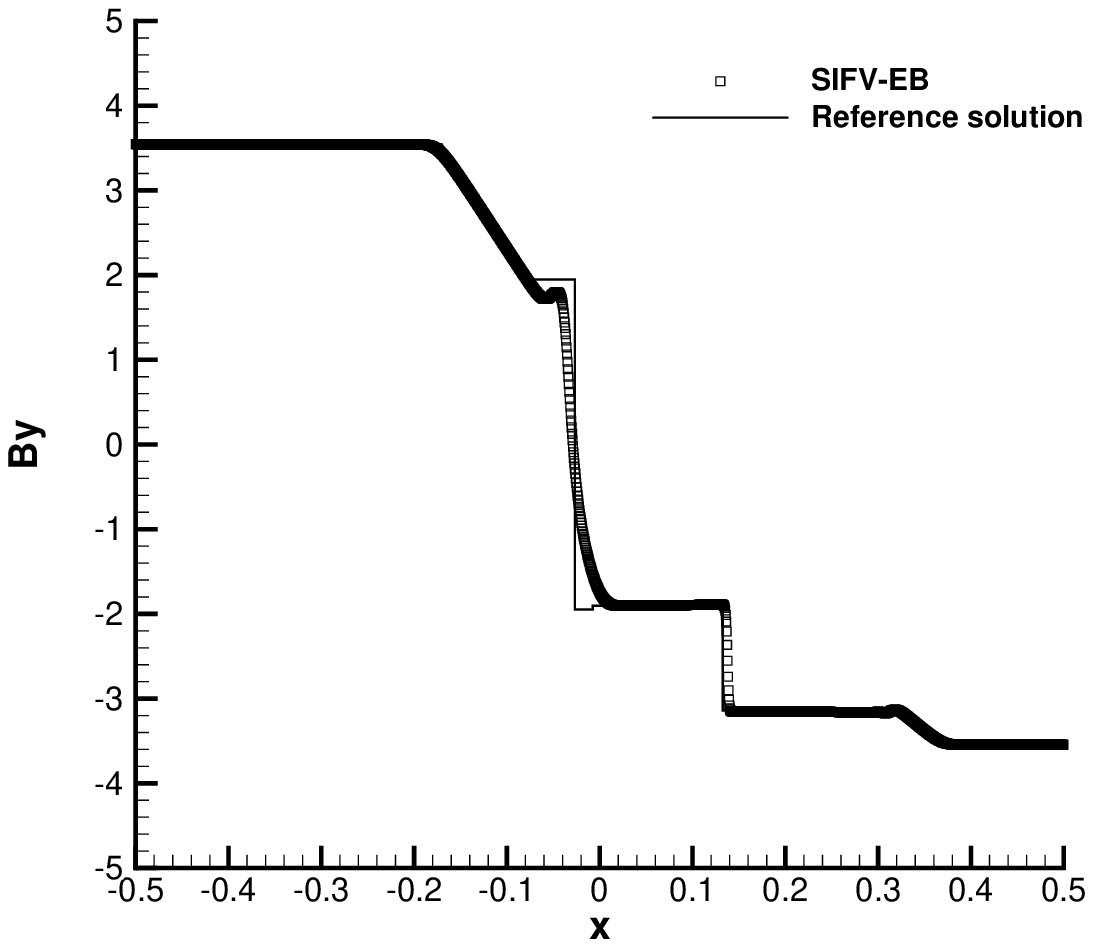}  \\
			\includegraphics[width=0.33\textwidth]{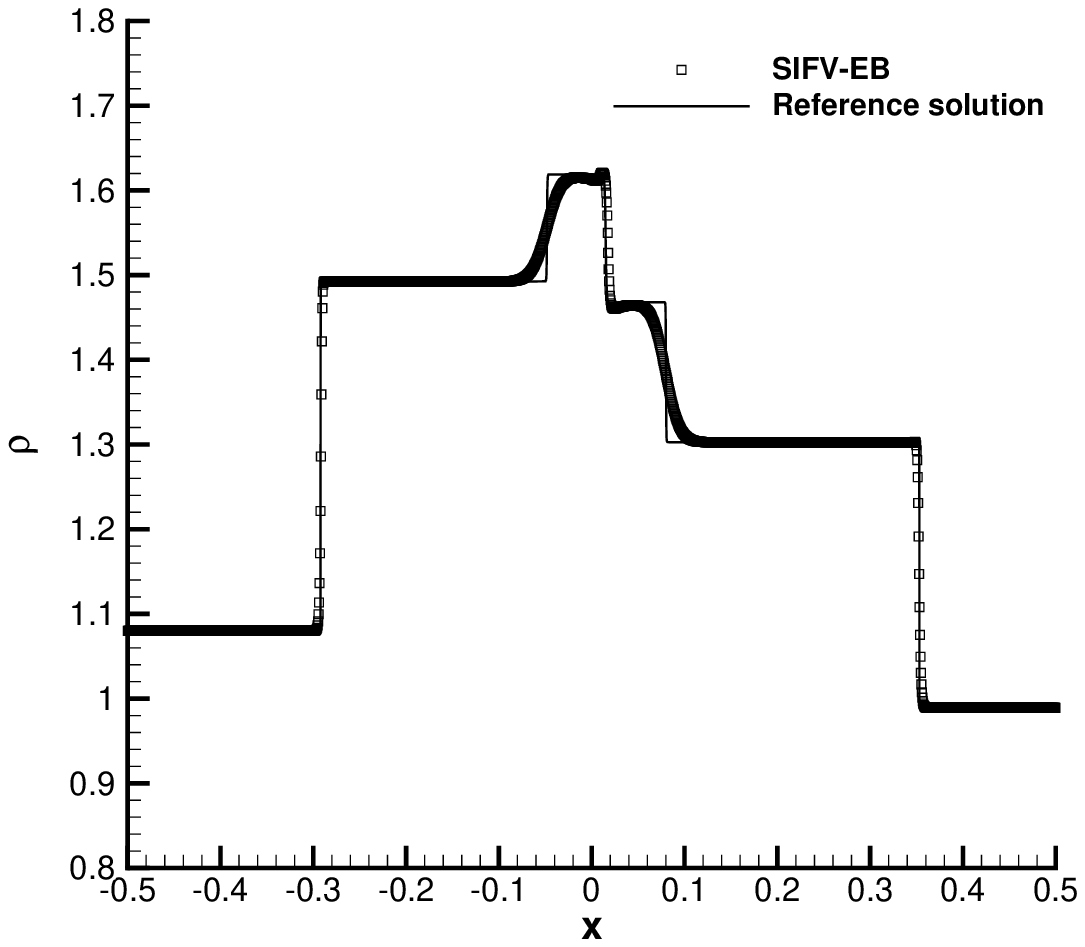}  &
			\includegraphics[width=0.33\textwidth]{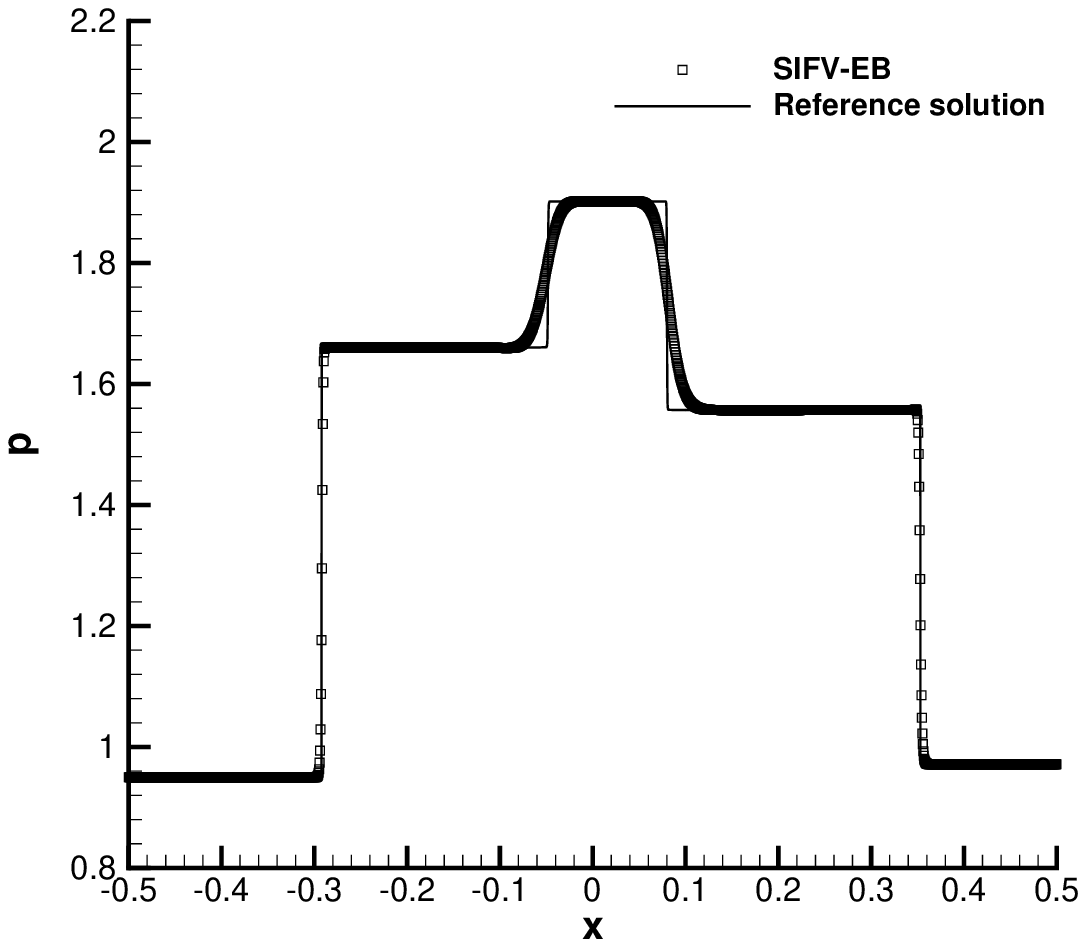} &
			\includegraphics[width=0.33\textwidth]{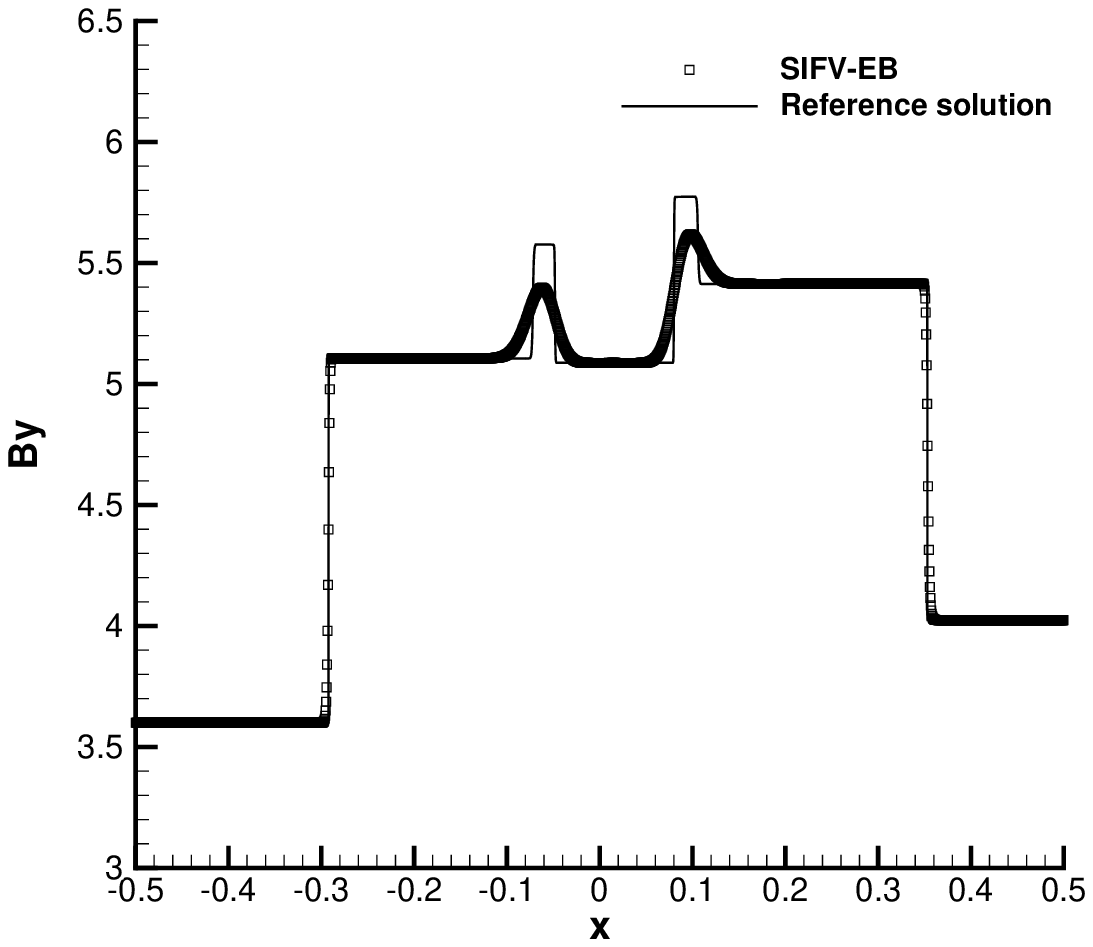}  \\
			\includegraphics[width=0.33\textwidth]{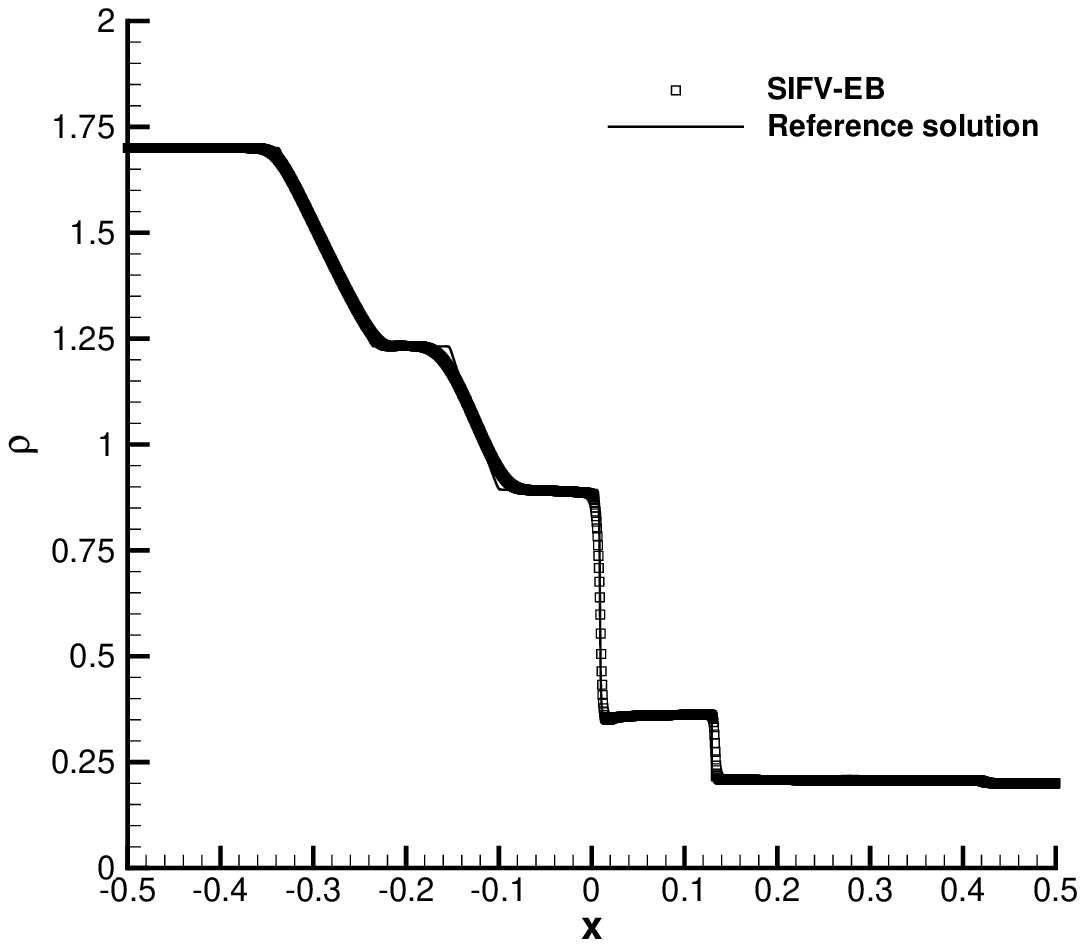}  &
			\includegraphics[width=0.33\textwidth]{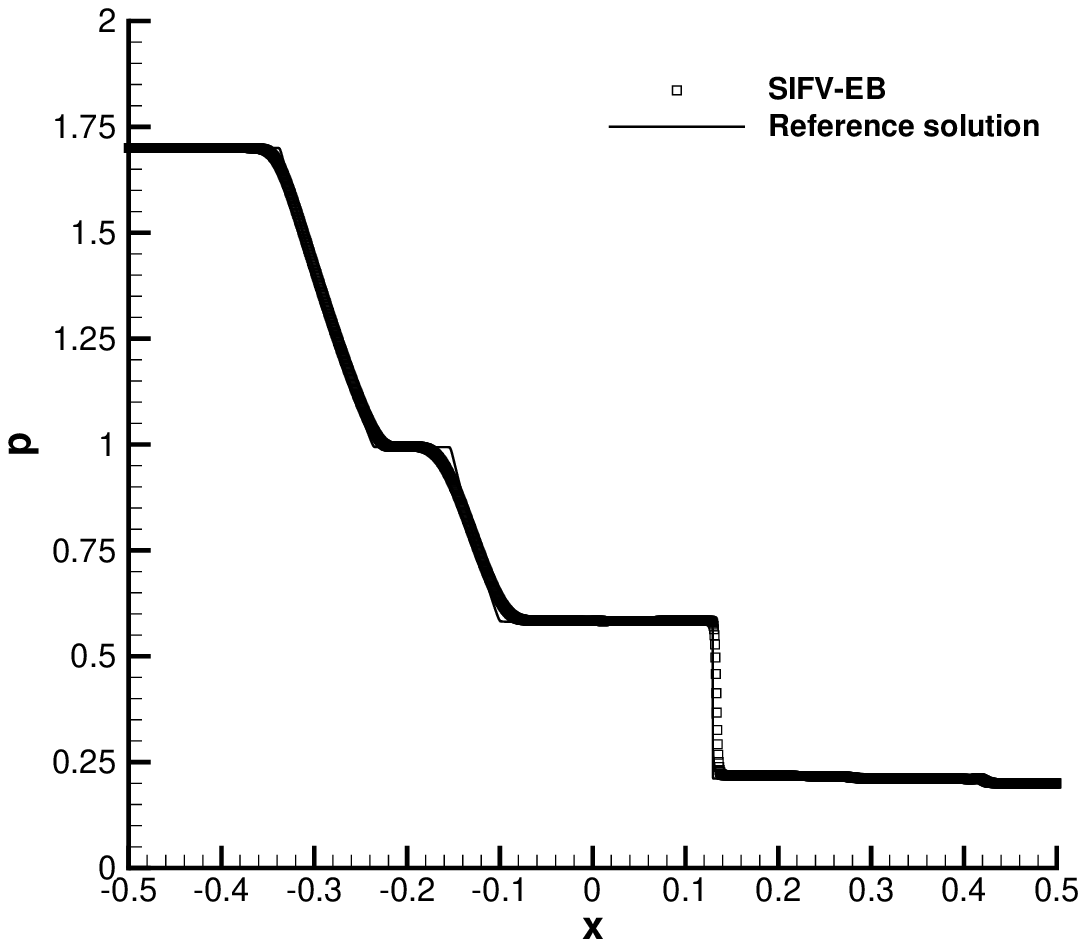} &
			\includegraphics[width=0.33\textwidth]{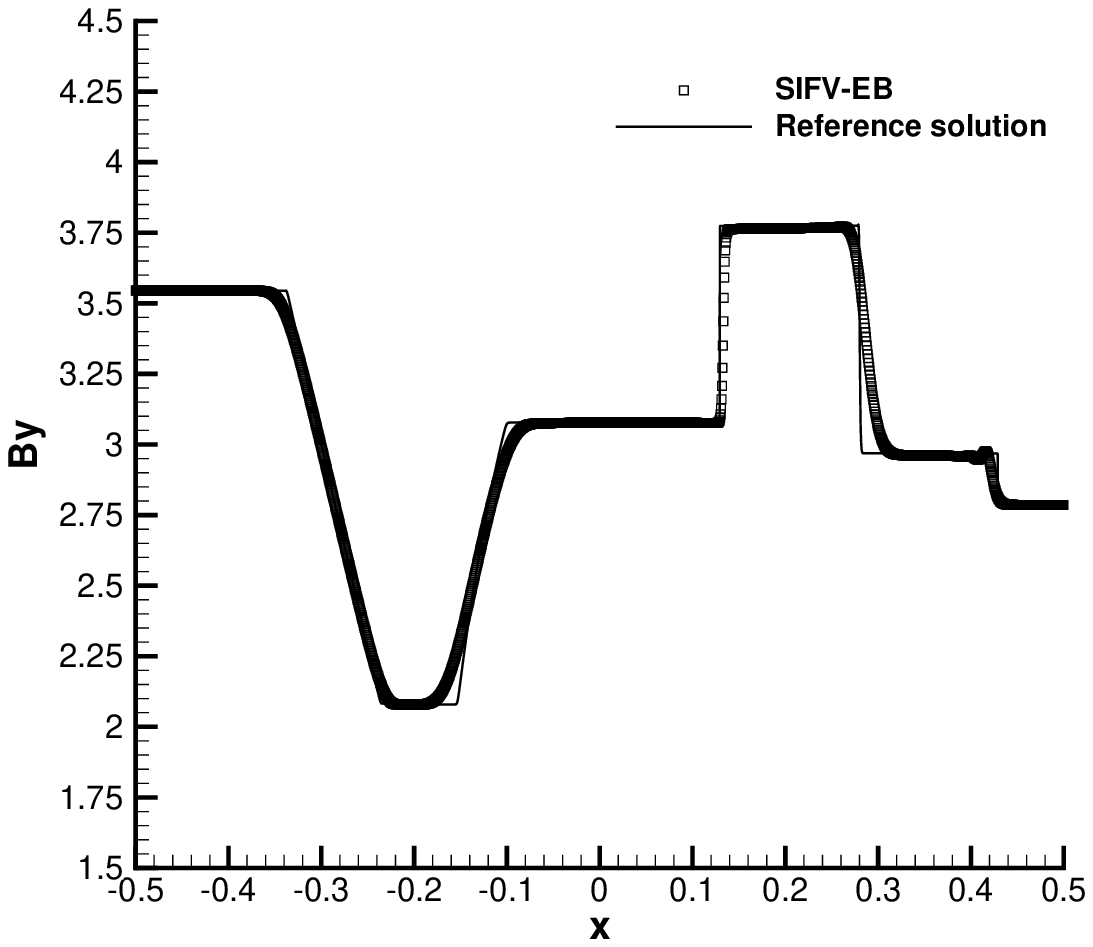}  \\
		\end{tabular}
		\caption{Riemann problems. RP1 at time $t_f=0.1$ (top), RP2 at time $t_f=0.2$ (middle) and RP3 at time $t_f=0.15$ (bottom). Left: density $\rho$. Center: magnetic field component $B_y$. Right: pressure $p$.}
		\label{fig.RP123}
	\end{center}
\end{figure}

\begin{figure}[!htbp]
	\begin{center}
		\begin{tabular}{ccc}
			\includegraphics[width=0.33\textwidth]{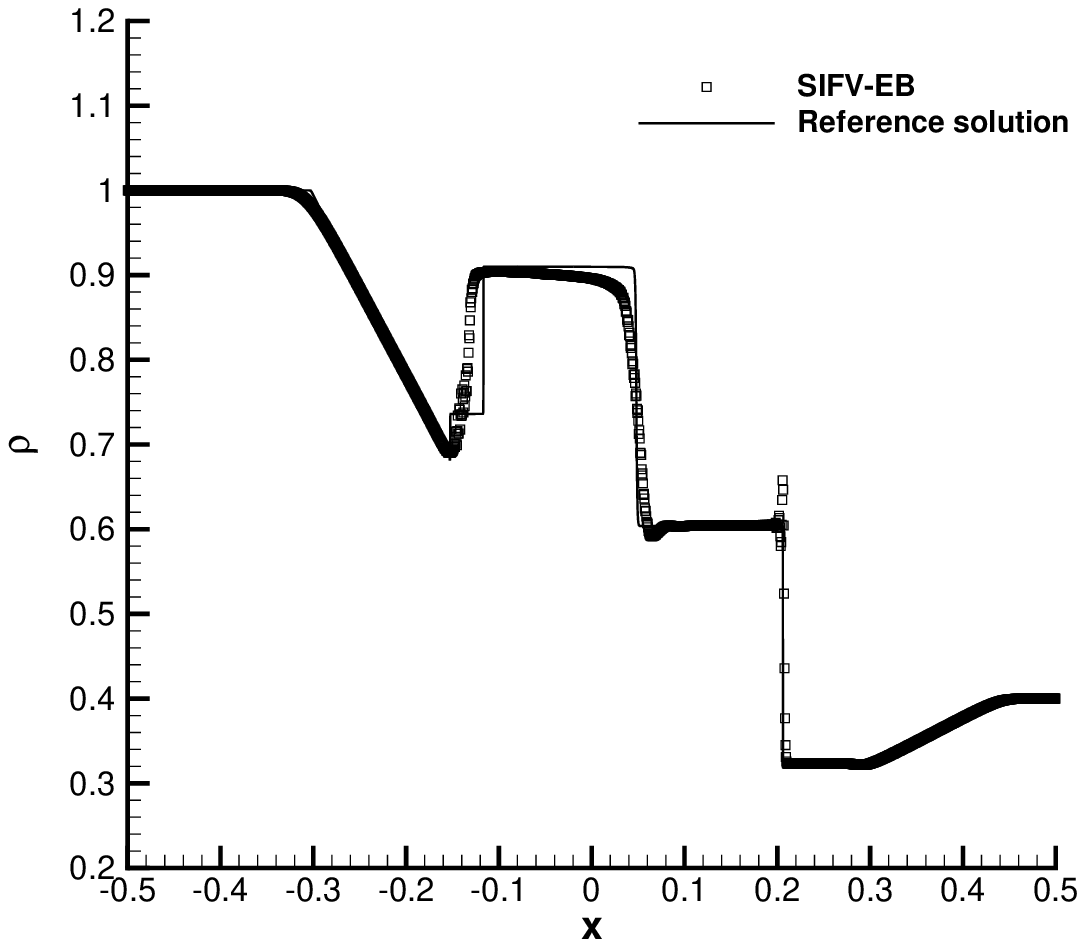}  &
			\includegraphics[width=0.33\textwidth]{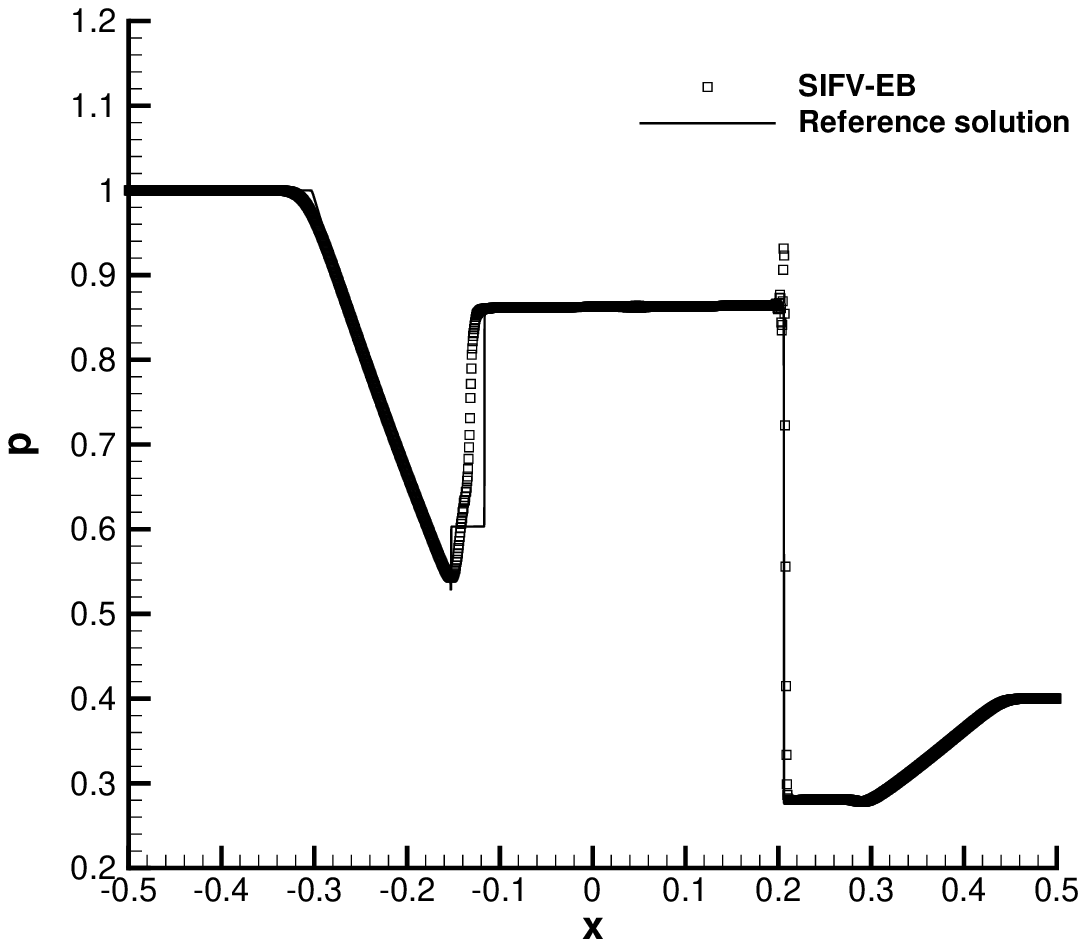} &
			\includegraphics[width=0.33\textwidth]{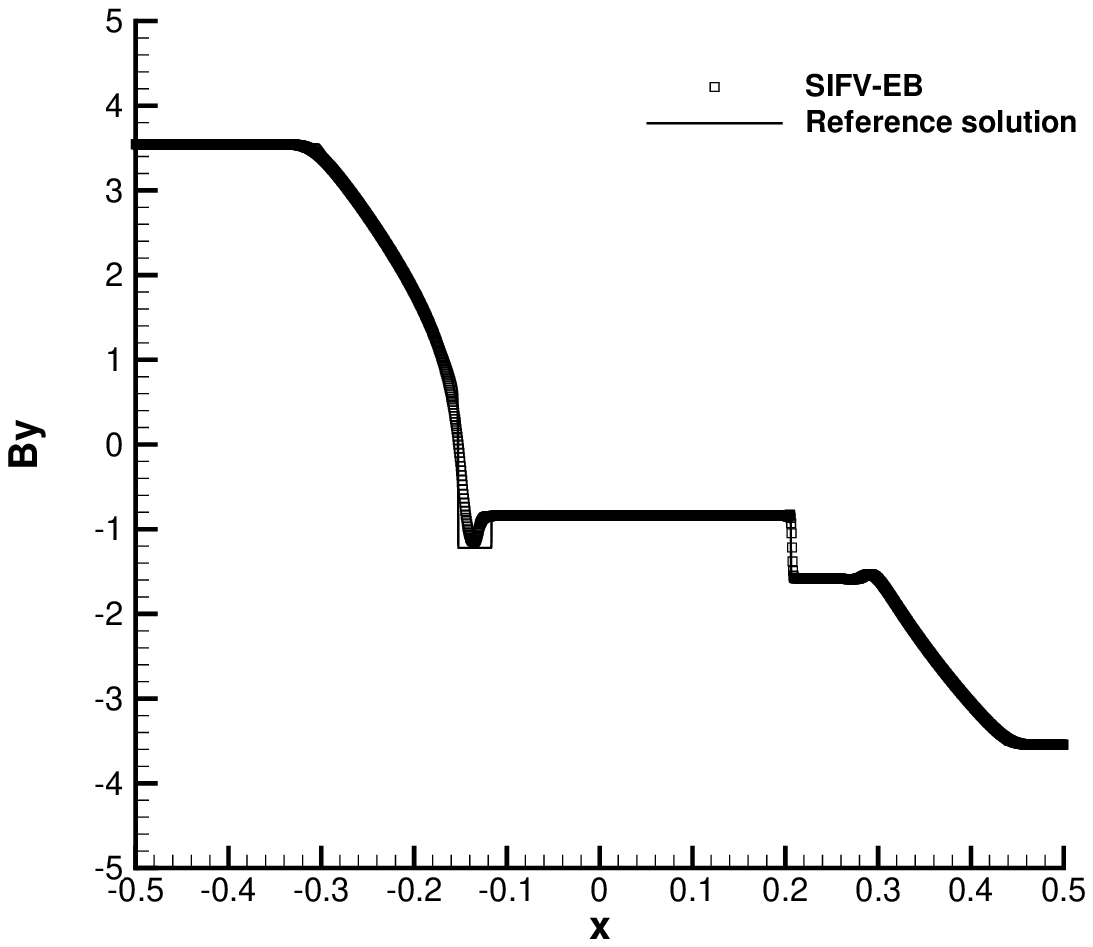}  \\
			\includegraphics[width=0.33\textwidth]{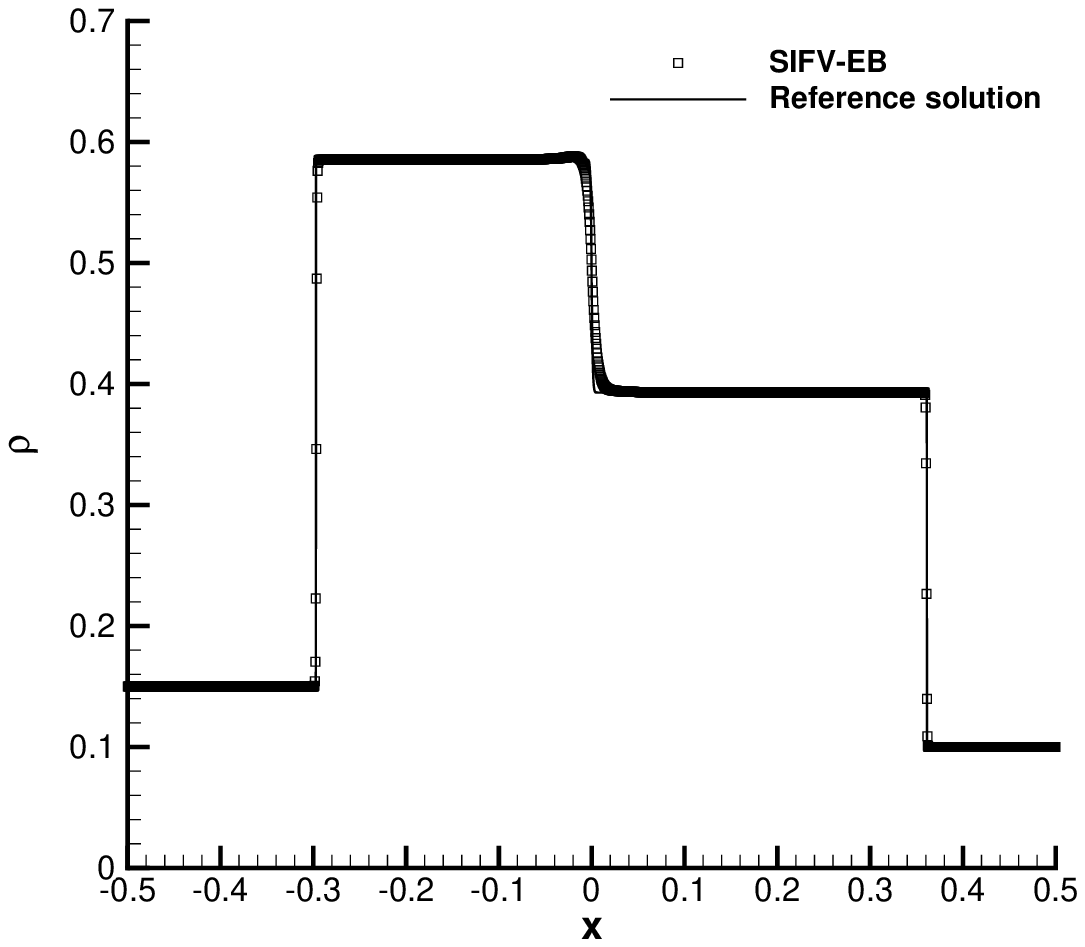}  &
			\includegraphics[width=0.33\textwidth]{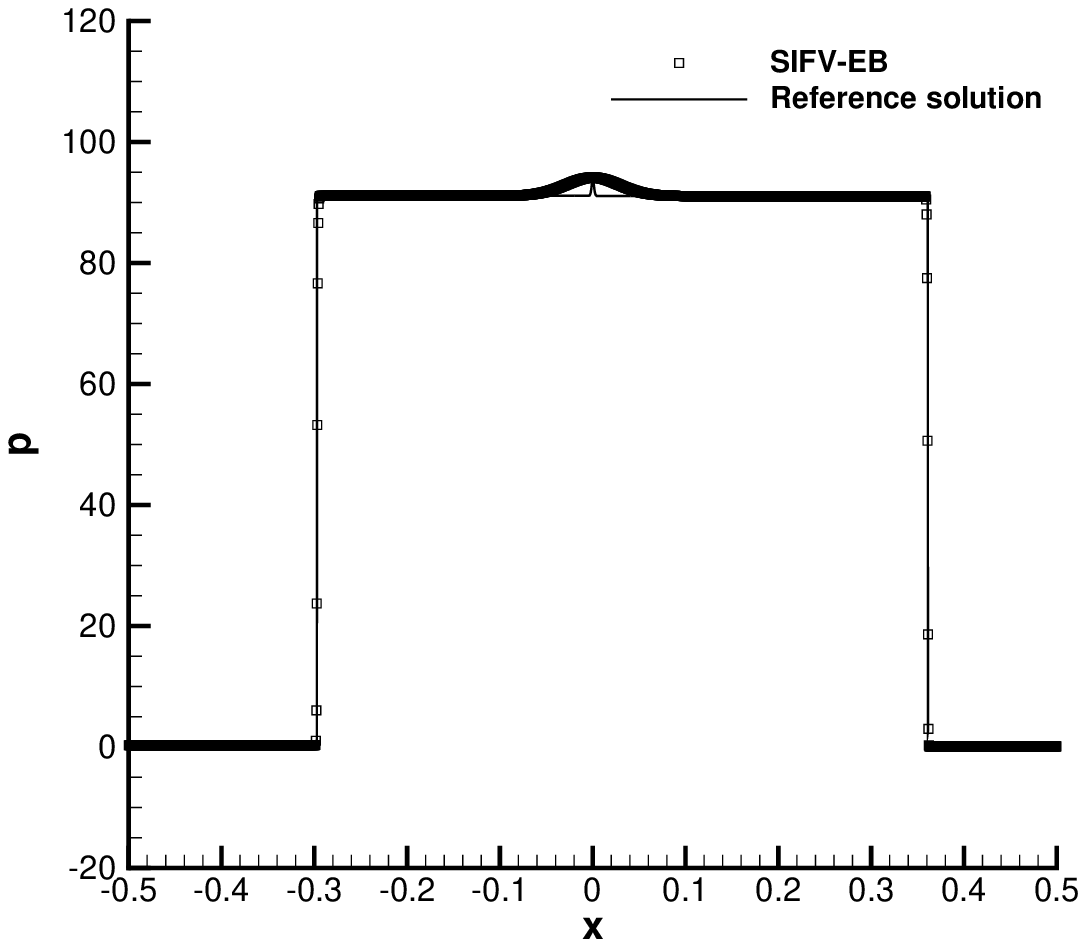} &
			\includegraphics[width=0.33\textwidth]{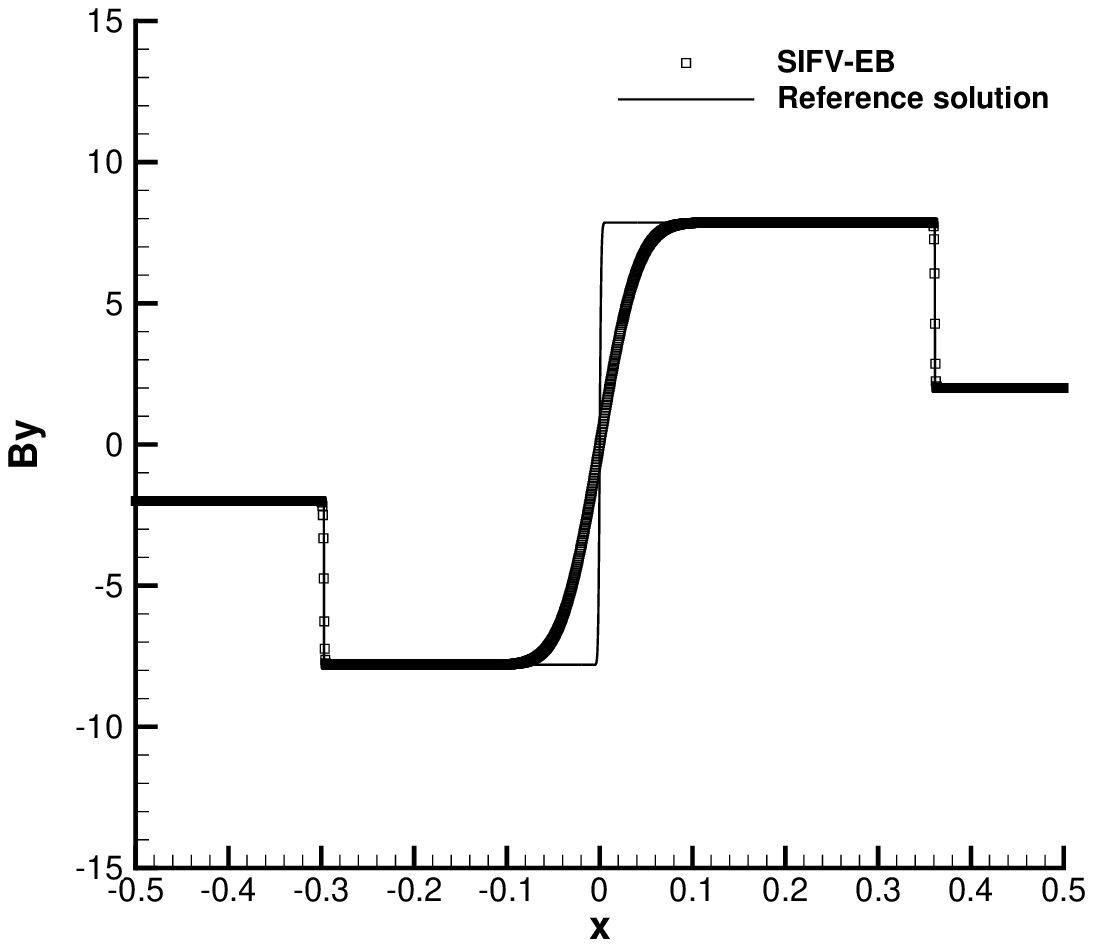}  \\
			\includegraphics[width=0.33\textwidth]{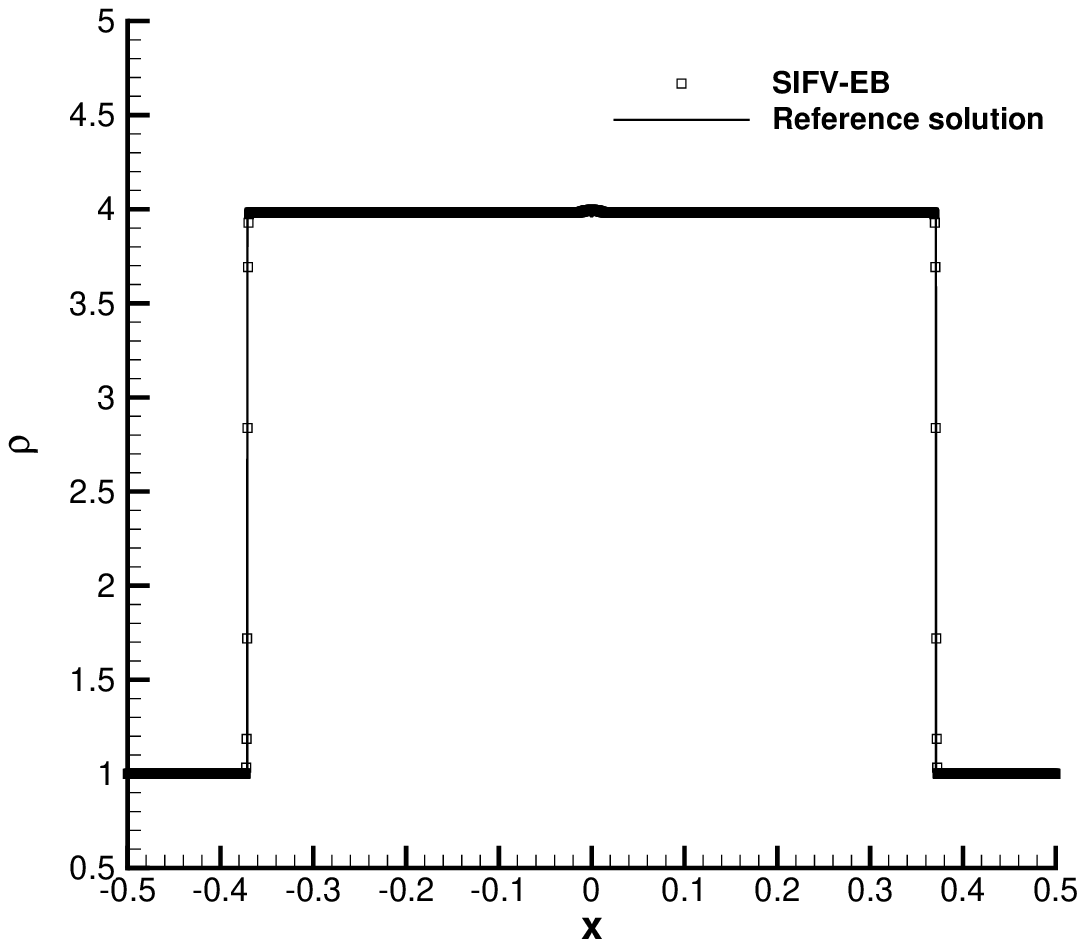}  &
			\includegraphics[width=0.33\textwidth]{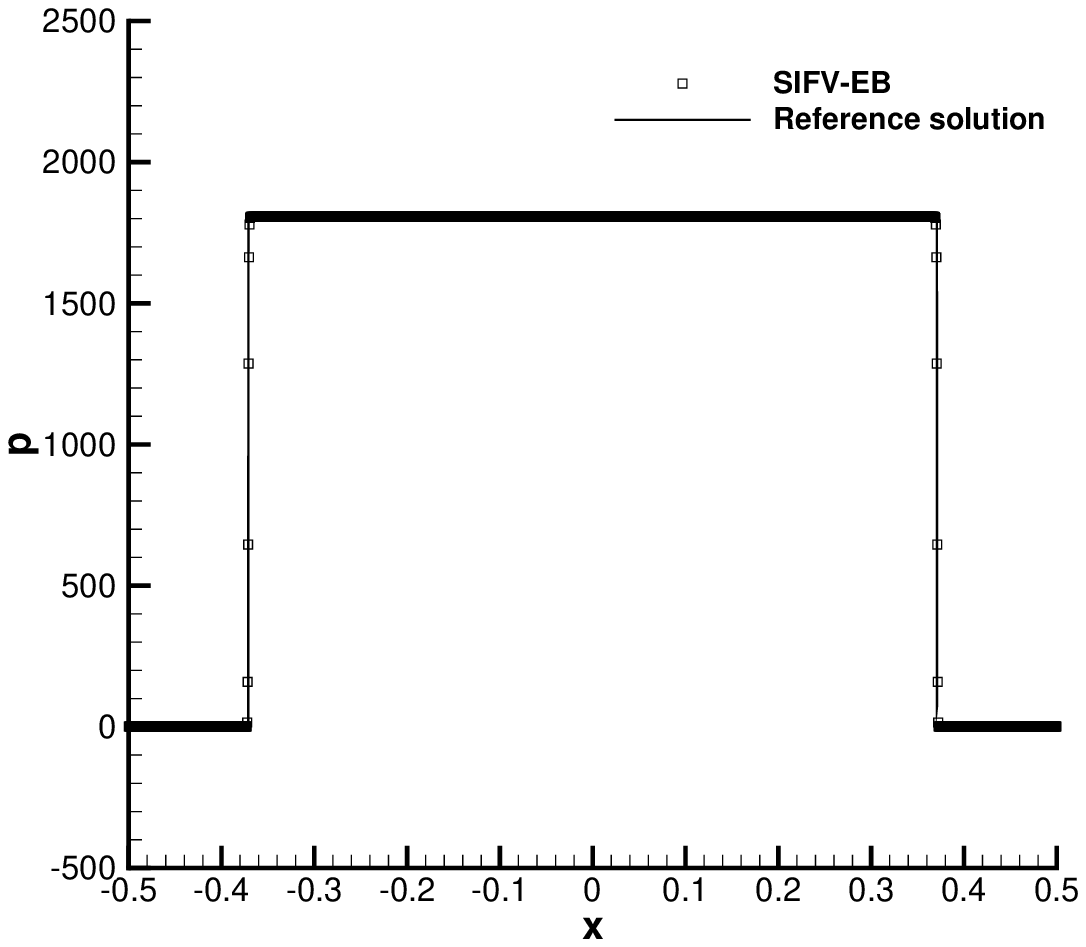} &
			\includegraphics[width=0.33\textwidth]{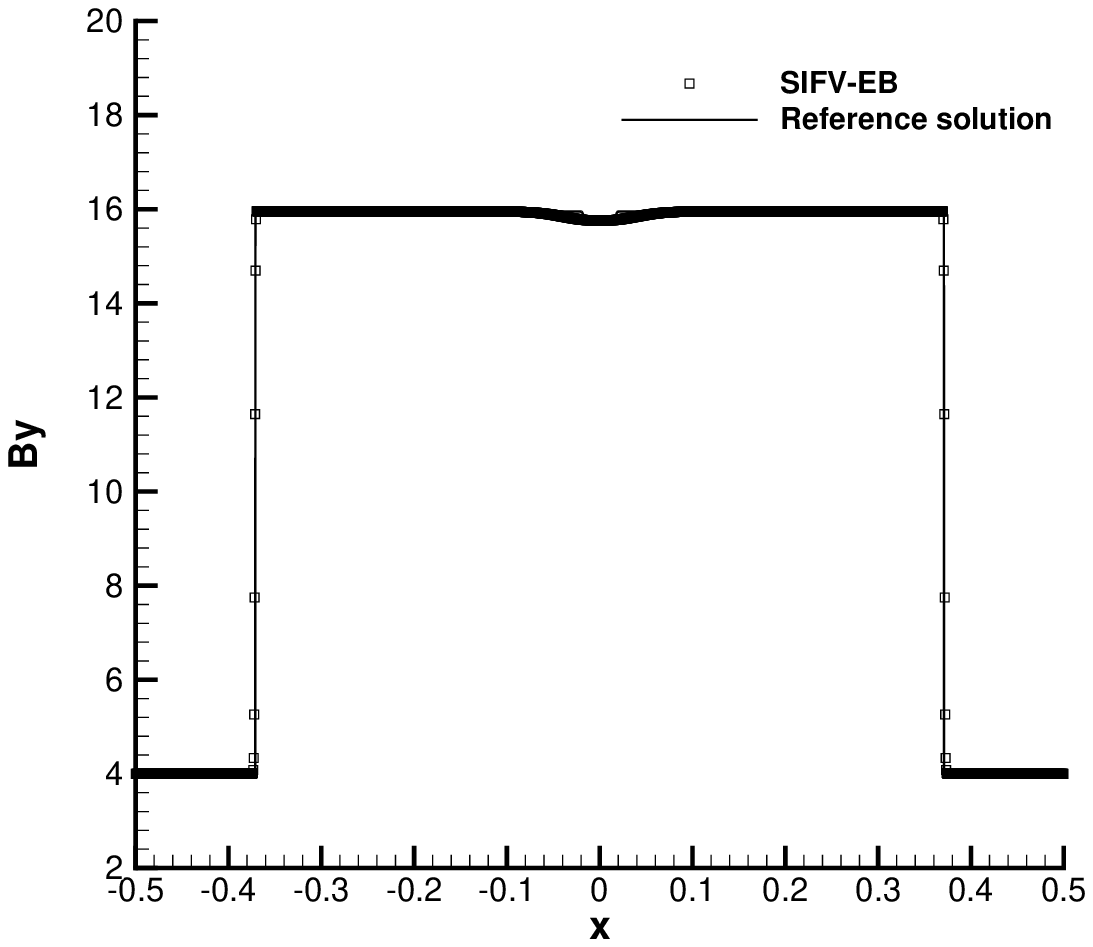}  \\
			\includegraphics[width=0.33\textwidth]{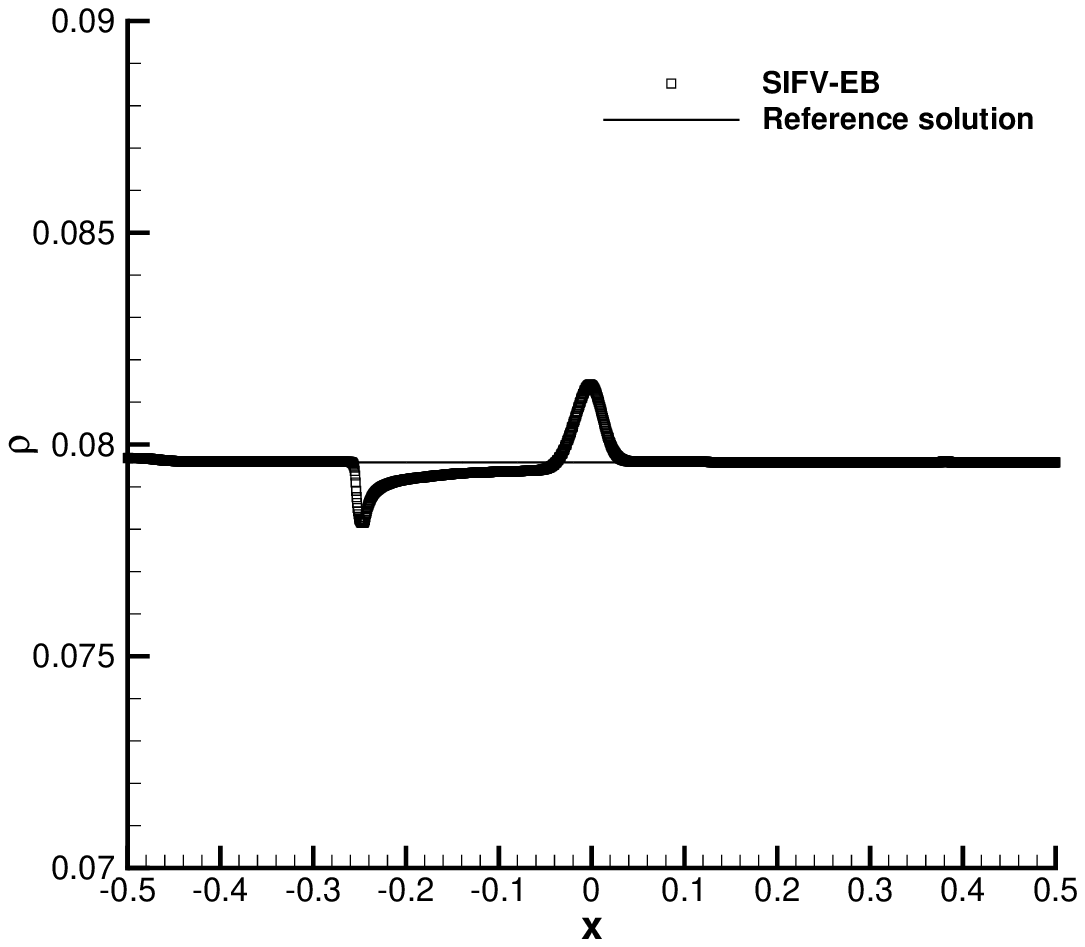}  &
			\includegraphics[width=0.33\textwidth]{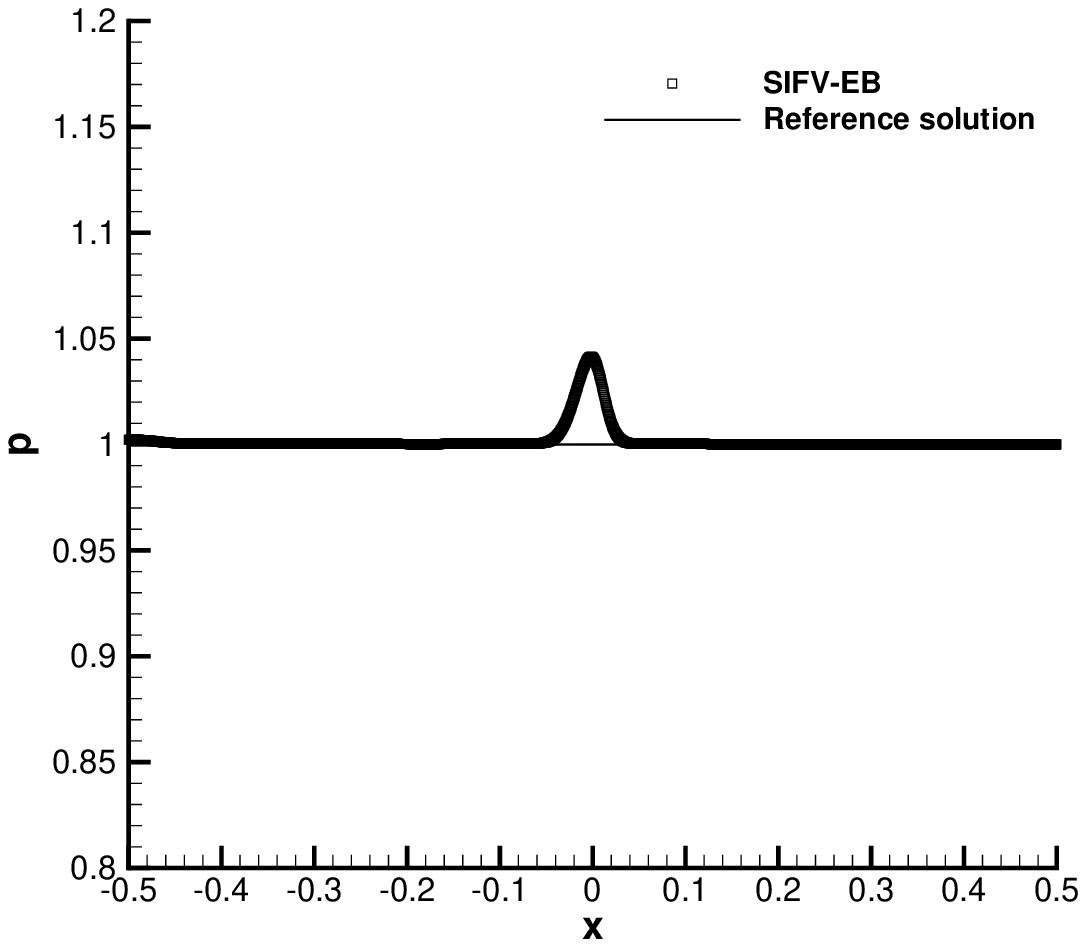} &
			\includegraphics[width=0.33\textwidth]{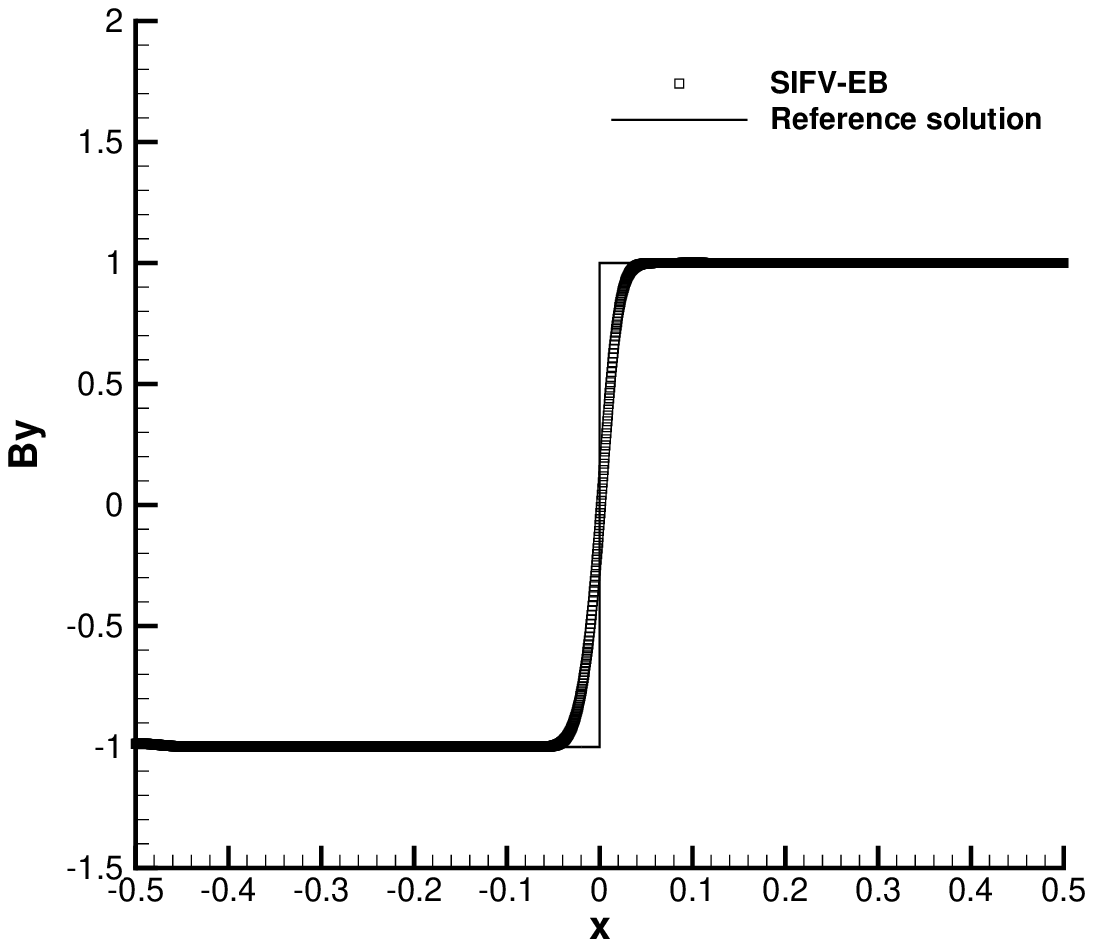}  \\
		\end{tabular}
		\caption{Riemann problems. RP4 at time $t_f=0.16$, RP5 at time $t_f=0.04$, RP6 at time $t_f=0.03$  and RP7 at time $t_f=0.2$ (from top to bottom). Left: density $\rho$. Center: magnetic field component $B_y$. Right: pressure $p$.}
		\label{fig.RP4567}
	\end{center}
\end{figure}

\begin{figure}[!htbp]
\begin{center}
	\begin{tabular}{c}
		\includegraphics[width=0.7\textwidth,keepaspectratio=true]{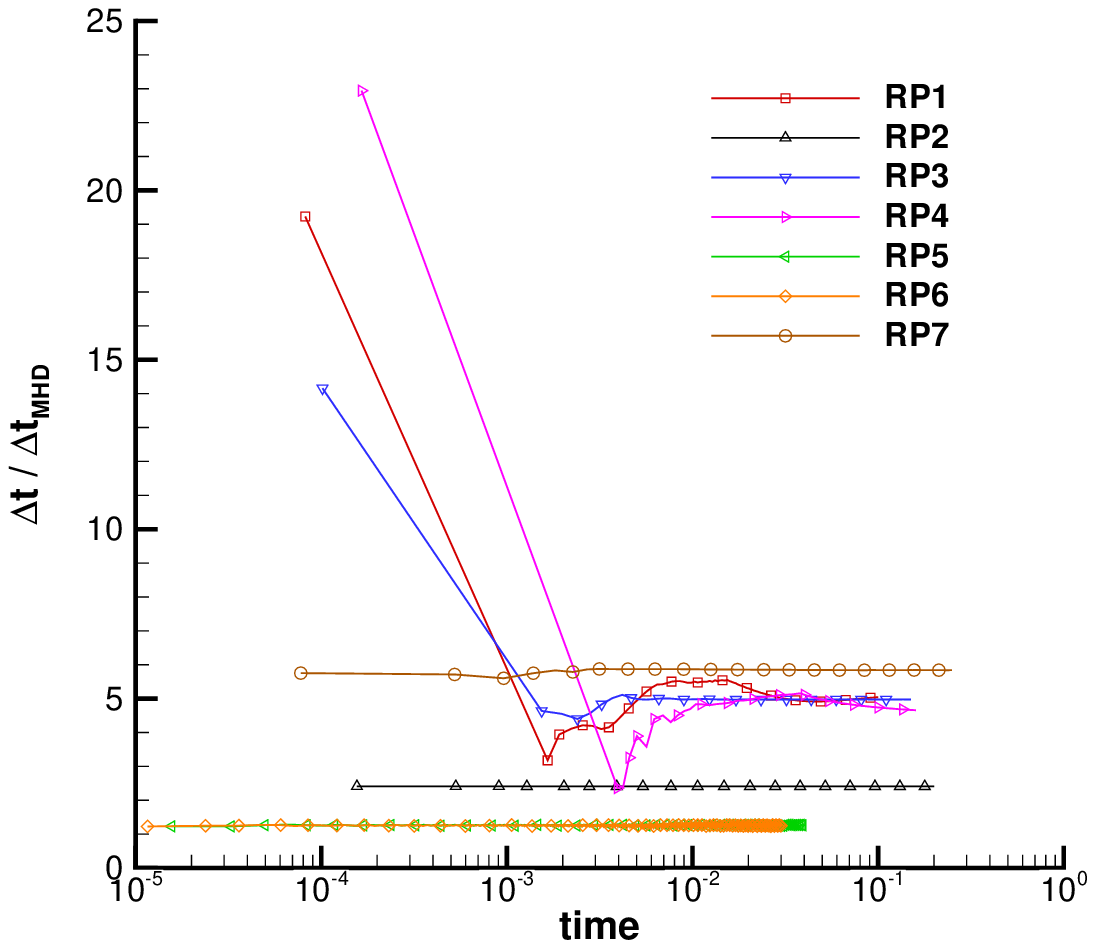}
	\end{tabular}
	\caption{Riemann problems. Time evolution of the ratio between the material time step of the \scheme scheme and the magneto-sonic time step of a fully explicit finite volume scheme for all the Riemann problems RP1--RP7.}
	\label{fig.RPdt}
\end{center}
\end{figure}

For each RP, depicted in Figures \ref{fig.RP123} and \ref{fig.RP4567}. The \scheme scheme is able to capture the correct shock position and all waves accurately for RP2--RP7.
In particular, the overshoot in density and pressure in RP1 is an artefact that also appears in approximate solutions produced with the HLLEM scheme from \cite{BalDum2016}, the Osher-type scheme \cite{DumToro2011} or the Rusanov scheme.
Moreover, RP7 contains an isolated steady Alfv\'en wave which is well resolved.
However, small perturbations of density and pressure appear around the jump position which then disappear as the mesh is refined.
The new \scheme scheme overall is clearly more diffusive than an explicit approximate Riemann solver. This is due to the semi-implicit nature and the material time stepping that allows much larger time steps, see Figure \ref{fig.RPdt} where a comparison between the material time stepping of the \scheme scheme with respect to a fully explicit time step, enforced by the fast magneto-sonic wave speeds $\lambda_{1,8}$, is given.
Even though most RP are placed in a compressible regime, our scheme allows for larger time steps that is reflected by slightly more diffusive fast traveling waves, while the material waves are captured accurately, as expected. Let us note that the first time step is always given by a magneto-sonic explicit time step, since many tests are initialized with a zero or small velocity which would lead to an infinitely large first time step and thus to inaccurate numerical results.

\subsection{Blast wave}
\begin{figure}[!htbp]
	\begin{center}
		\begin{tabular}{cc}
			\includegraphics[trim=2 2 2 2,clip,width=0.47\textwidth]{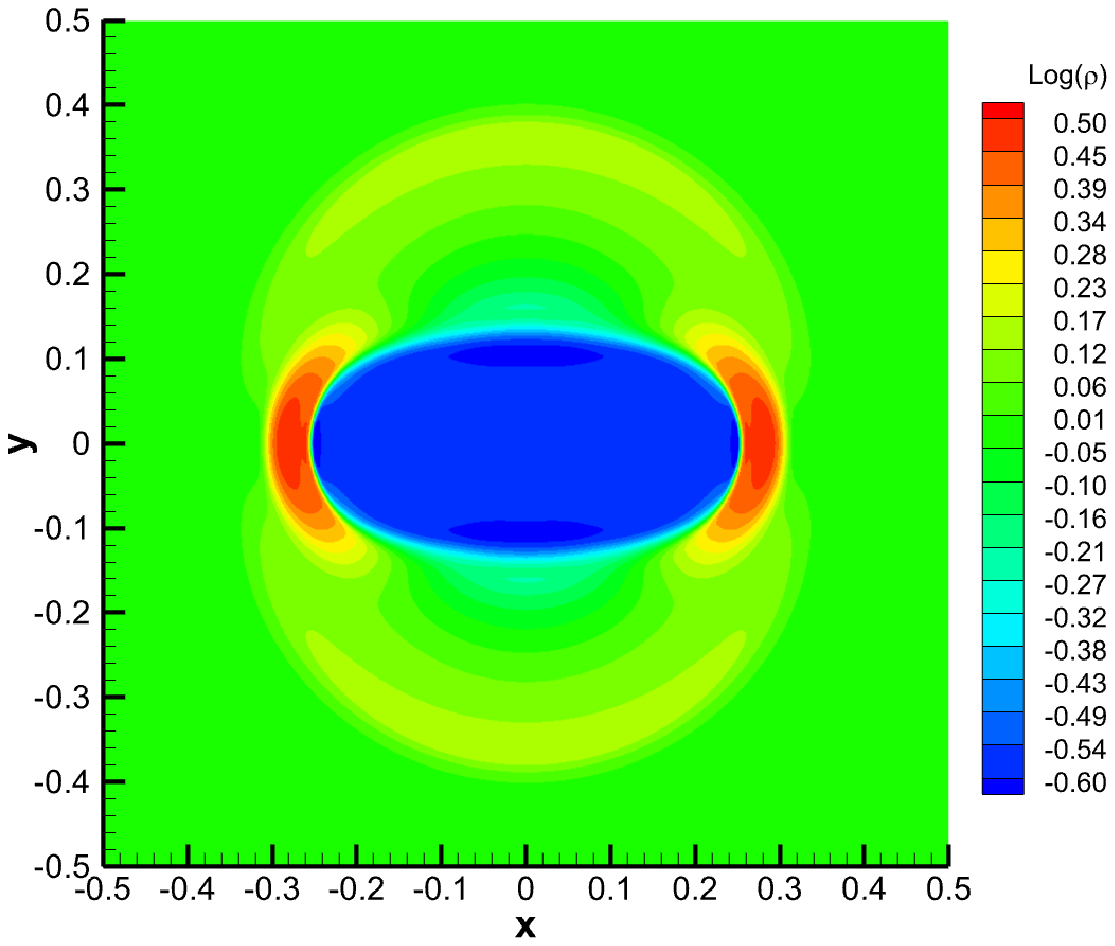} &
			\includegraphics[trim=2 2 2 2,clip,width=0.47\textwidth]{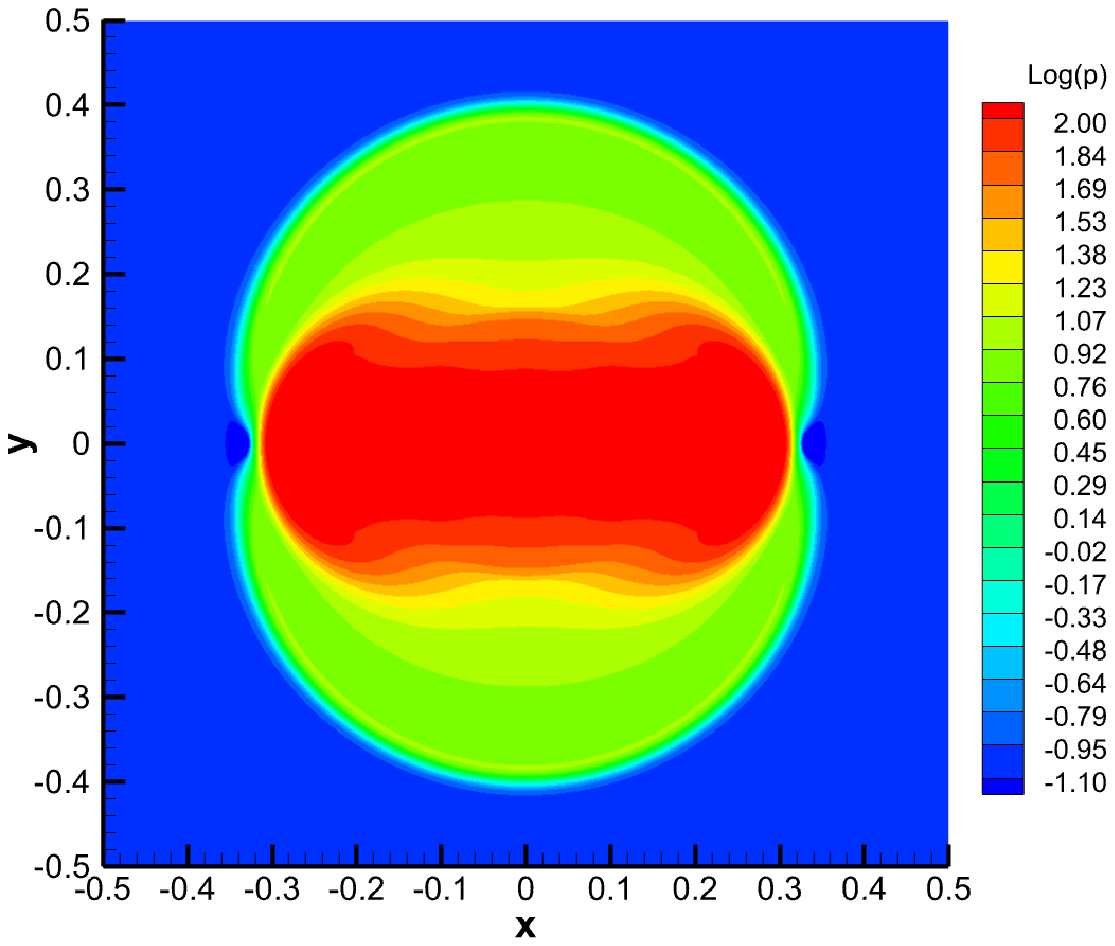} \\
			\includegraphics[trim=2 2 2 2,clip,width=0.47\textwidth]{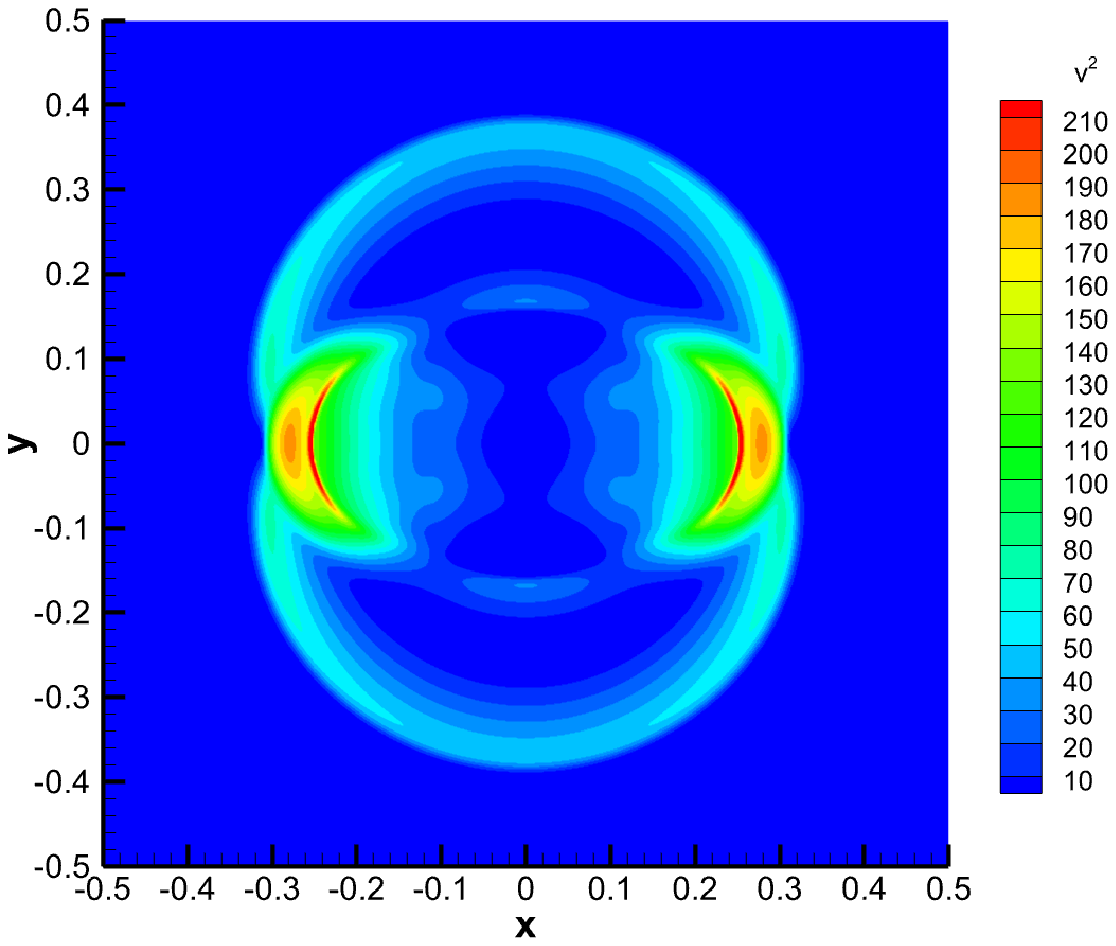} &
			\includegraphics[trim=2 2 2 2,clip,width=0.47\textwidth]{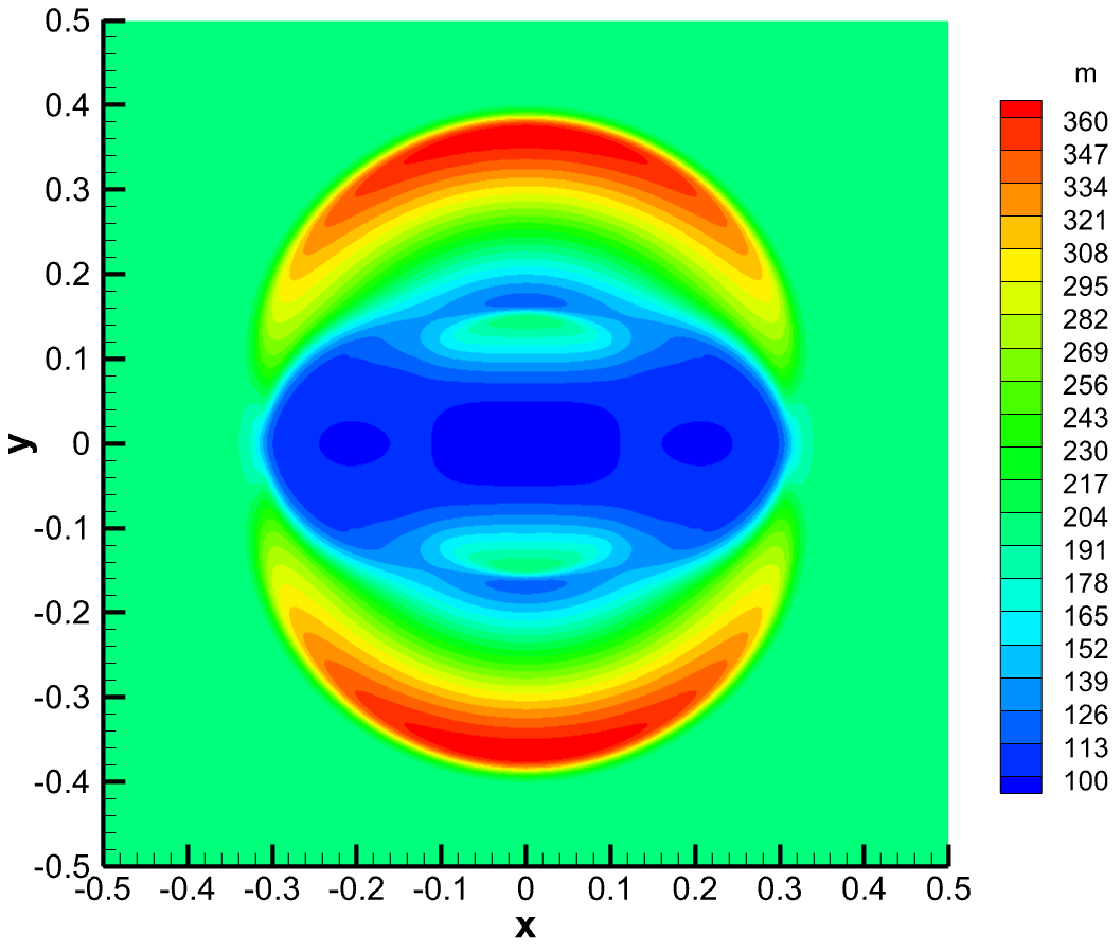} \\
		\end{tabular}
		\caption{Blast problem at time $t_f=0.01$. Results for the logarithm (base 10) of density and pressure, magnitude of the velocity and magnetic pressure (from top left to bottom right).}
		\label{fig.MHDBlast}
	\end{center}
\end{figure}

The MHD blast wave problem introduced in \cite{BalsaraSpicer1999} represents a challenging test since it contains a strong magnetic field and a jump in the pressure of four orders of magnitude.
The initial density is constant $\rho = 1$, the initial velocity field is set to zero $(u,v,w) = (0,0,0)$ and the magnetic field is $\B = (70,0,0)$. The computational domain $\Omega=[-0.5;0.5]\times[-0.5;0.5]$ is paved with $N_x \times N_y=1024 \times 1024$ cells, and Dirichlet boundary conditions are imposed everywhere. Inside the circular region $r < 0.1$ with $r = \sqrt{x^2 + y^2}$, the pressure is set to $p = 1000$ whereas outside it is assigned with $p = 0.1$.
Thus initially, the Alfv\'en waves are faster than the sound speed inside the circular domain, that are given by $c_A = 70$ and $c = 37.4166$, respectively.
Here, we set $\gamma = 1.4$ and the final time is chosen to be $t_f = 0.01$.
The results at final time for density, pressure, magnitude of the velocity computed as $\vv^2$ as well as the magnetic pressure $m$ are displayed in Figure \ref{fig.MHDBlast} for the second-order \scheme scheme.
We observe that there are no oscillations arising in the simulation and all waves are captured accurately. The results compare qualitatively well with the literature \cite{LagrangeMHD,divFreeMHD_BalDum2015,3splitMHD}.

\subsection{Ideal MHD rotor problem}
\begin{figure}[!htbp]
	\begin{center}
		\begin{tabular}{cc}
			\includegraphics[trim=2 2 2 2,clip,width=0.47\textwidth]{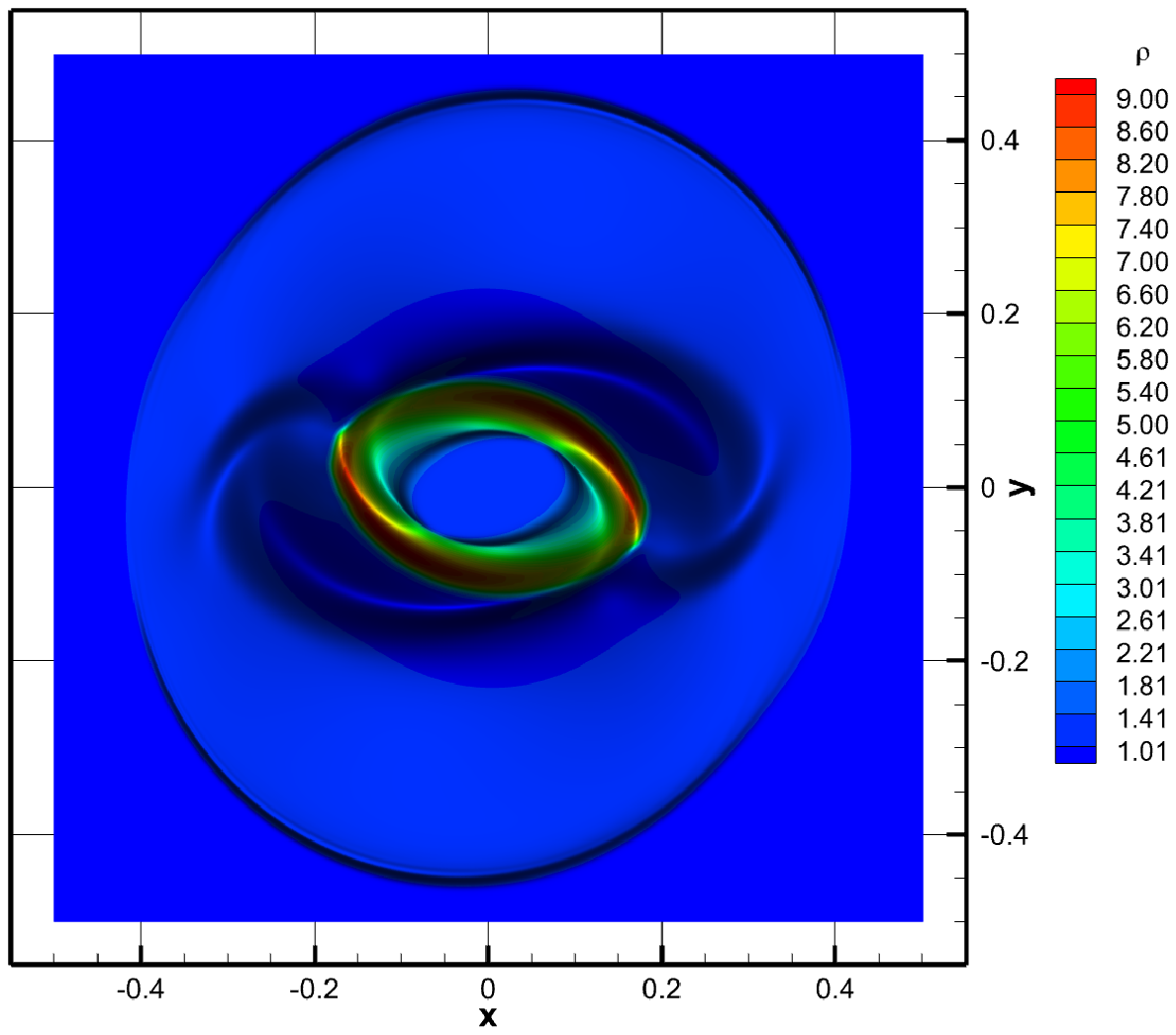} &
			\includegraphics[trim=2 2 2 2,clip,width=0.47\textwidth]{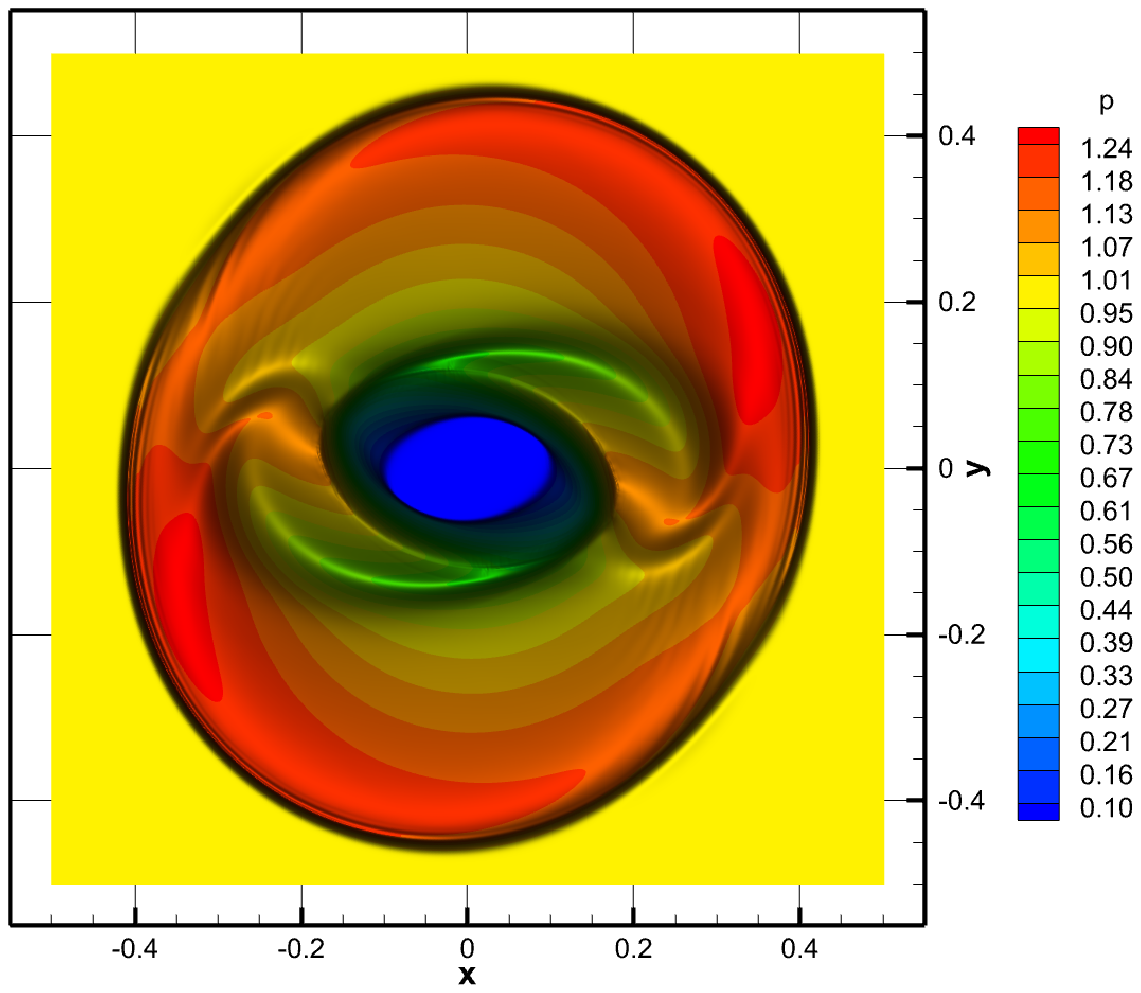} \\
			\includegraphics[trim=2 2 2 2,clip,width=0.47\textwidth]{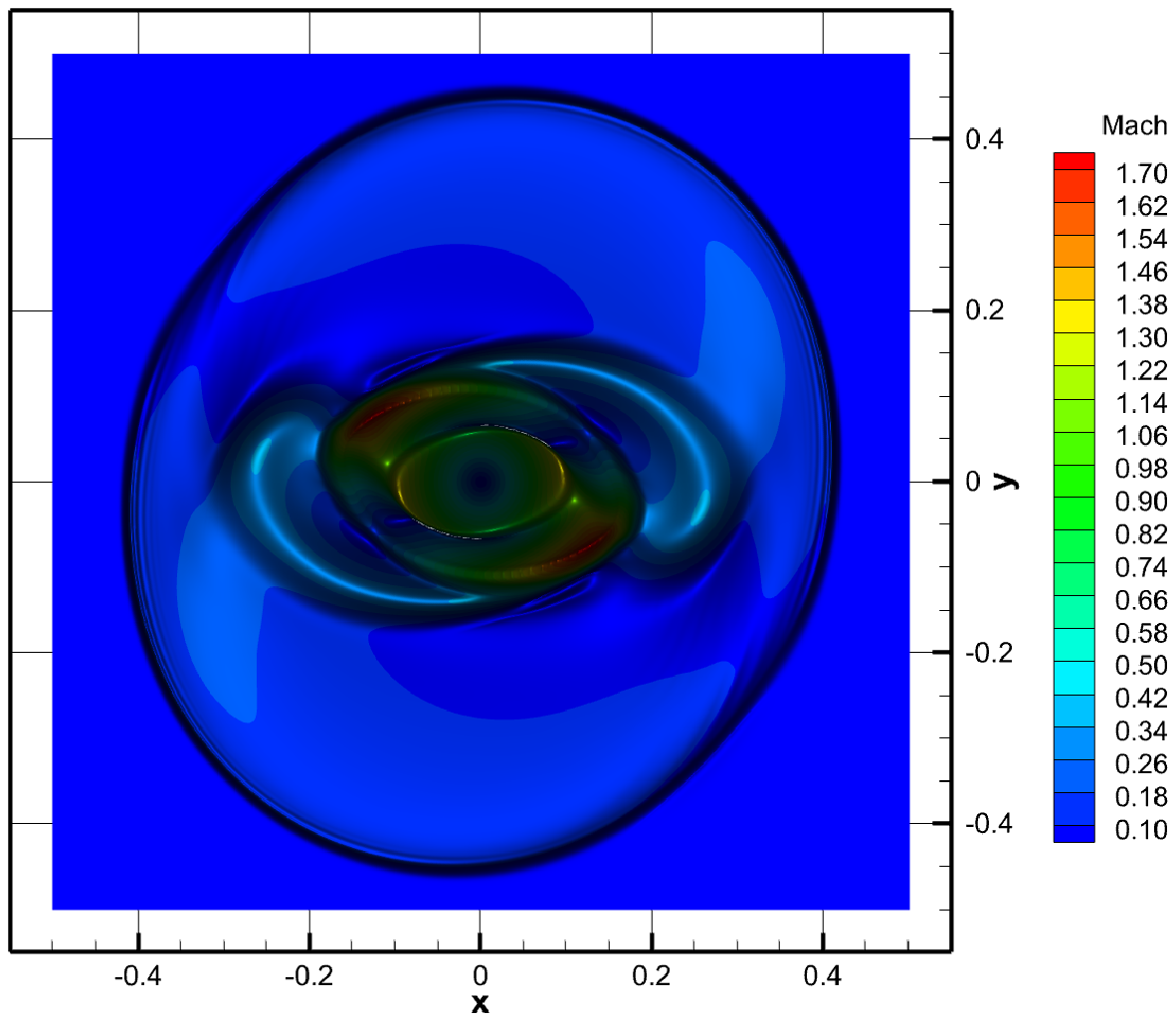} &
			\includegraphics[trim=2 2 2 2,clip,width=0.47\textwidth]{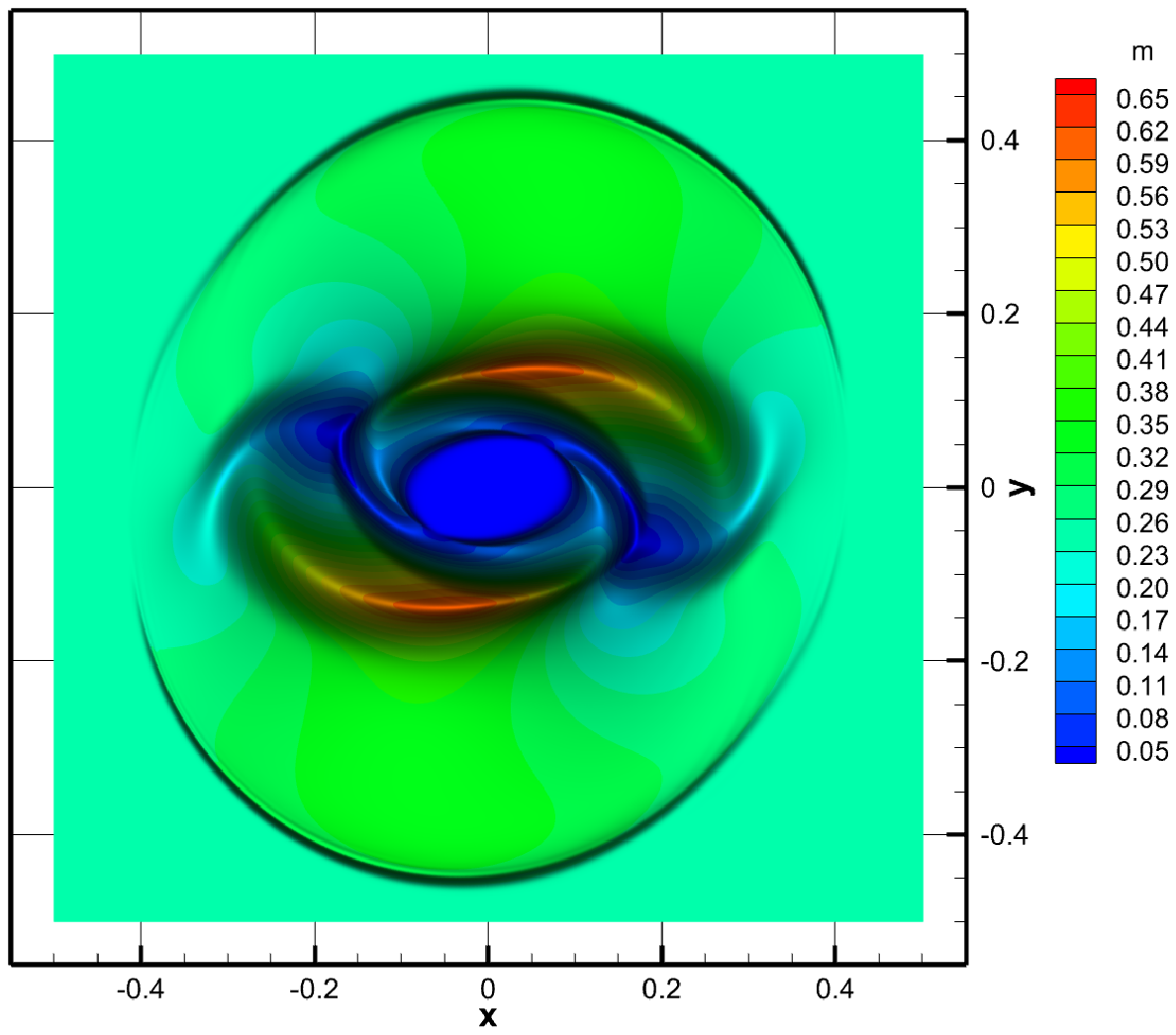} \\
		\end{tabular}
		\caption{Rotor problem at time $t_f=0.25$. Results for density, pressure, Mach number and magnetic pressure (from top left to bottom right).}
		\label{fig.MHDRotor}
	\end{center}
\end{figure}

The next test case concerns a benchmark for the ideal MHD equations, the so-called rotor problem proposed in \cite{BalsaraSpicer1999}, which is characterized by a big density difference inside and outside of the rotor. The computational domain is given by the square $\Omega=[-0.5;0.5]\times[-0.5;0.5]$ with Dirichlet boundaries, and it is discretized using $N_x \times N_y=1024 \times 1024$ control volumes. The background pressure and magnetic field are constant and set to $p = 1$ and $\B = (2.5,0,0)^T$, respectively. The corresponding magnetic vector potential in $z$-direction is $A_z = 2.5 y$. The density and velocity field are defined differently inside and outside of the rotor, i.e.
\begin{equation}
    \rho = \begin{cases}
        10 & \text{if} \quad 0 \leq r \leq 0.1, \\
        1  & \text{else},
    \end{cases}\quad \quad
    (u,v,w) =
    \begin{cases}
         (0,0,10) \times (x,y,z) & \text{if} \quad 0 \leq r \leq 0.1,\\
         (0,0,0) & \text{else.}
    \end{cases}
\end{equation}
The ratio of specific heats is set to $\gamma = 1.4$ and the simulation is run until the final time $t_f=0.25$. In Figure \ref{fig.MHDRotor}, the results for the second-order \scheme scheme are plotted for the density, pressure, acoustic Mach number $\Mc$ and the magnetic pressure $m$.
The scheme, though spanning a Mach number regime from weakly compressible to compressible with $\Mc > 1$, yields results that compare well with simulations reported in literature see e.g. \cite{SIMHD_Dumbser2019,LagrangeMHD}.

\begin{figure}[!htbp]
	\begin{center}
		\begin{tabular}{cc}
			\includegraphics[trim=2 2 2 2,clip,width=0.47\textwidth]{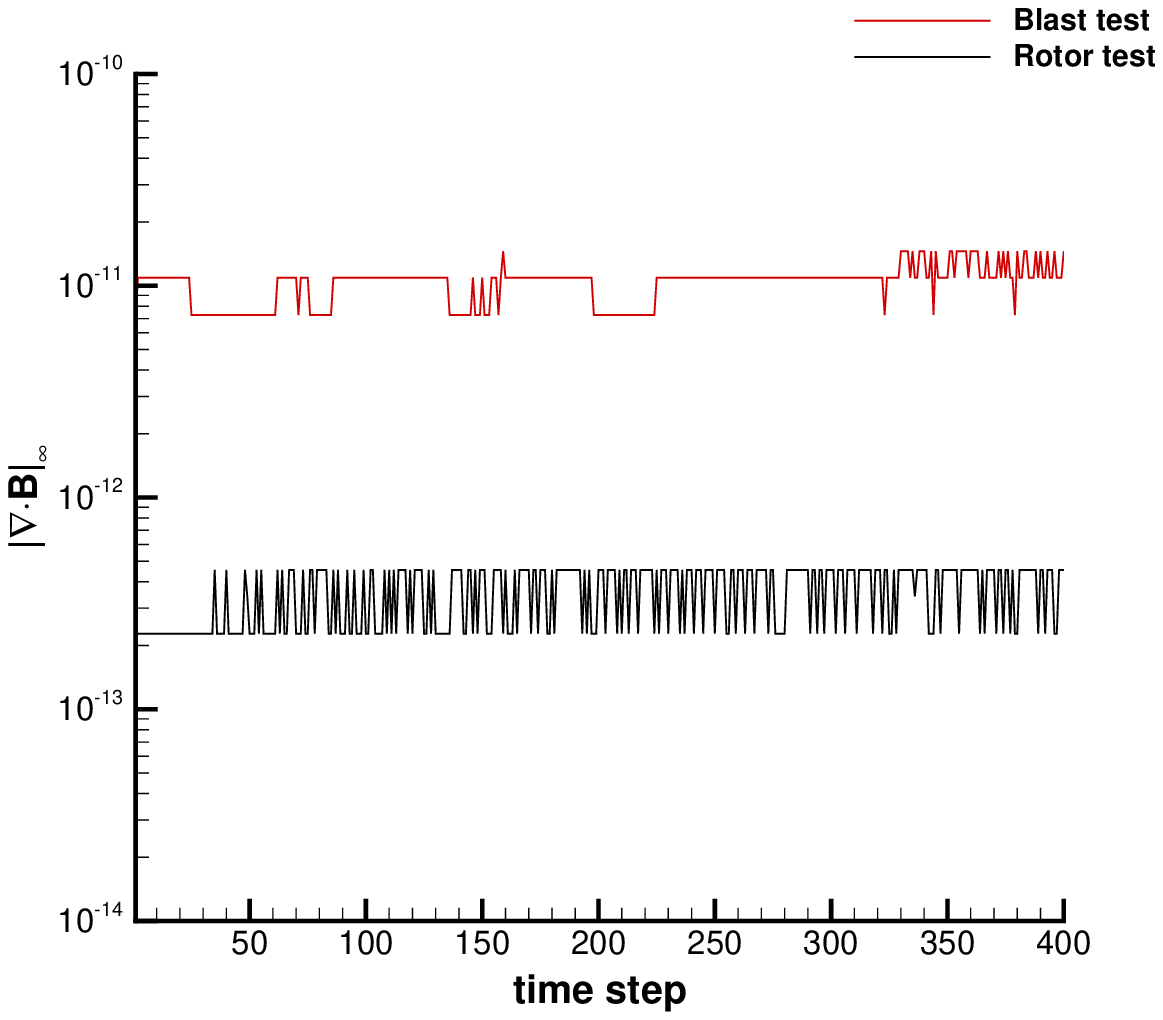} &
			\includegraphics[trim=2 2 2 2,clip,width=0.47\textwidth]{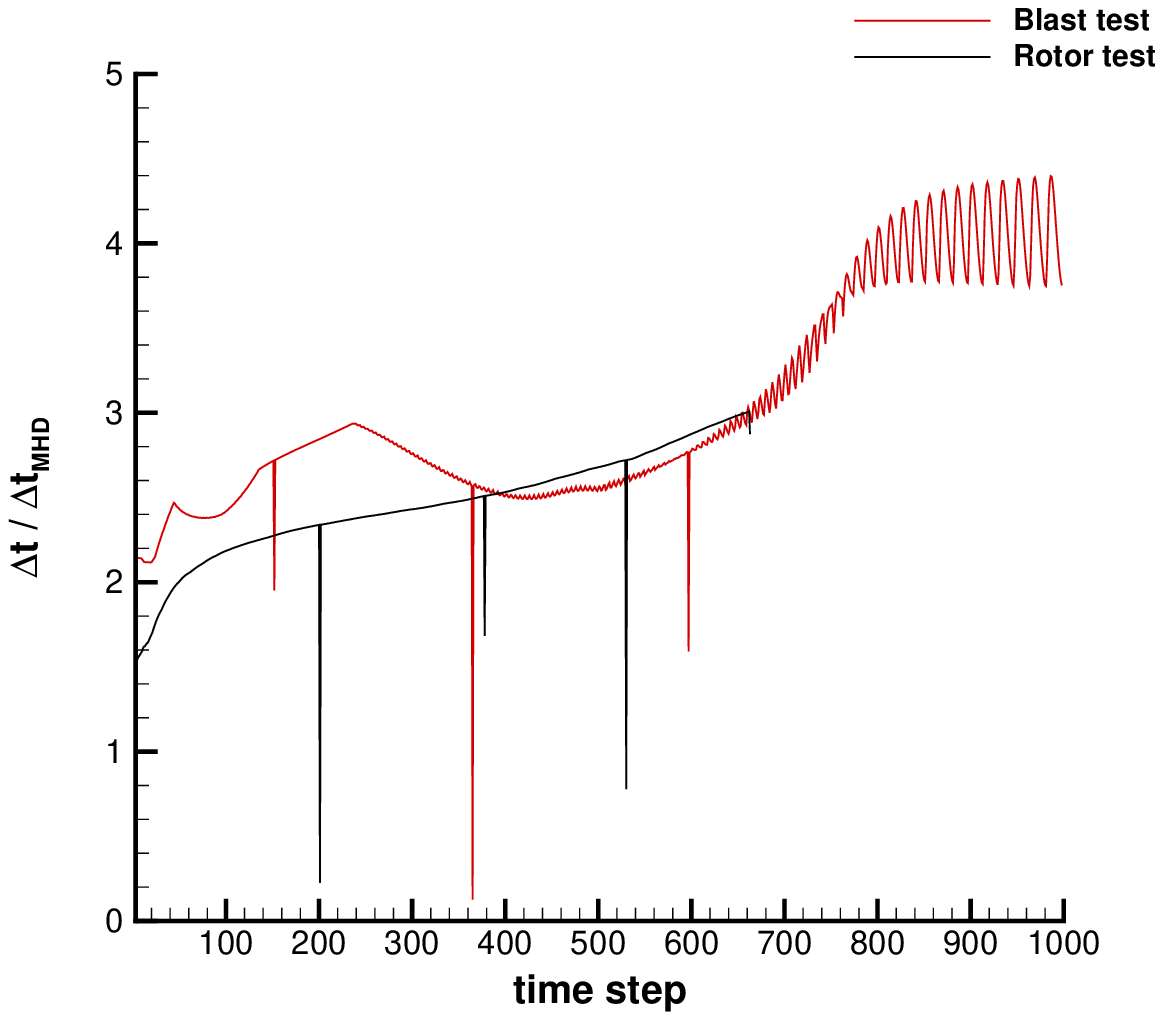} \\
		\end{tabular}
		\caption{Blast problem and rotor problem. left: time evolution of the infinity norm of the divergence of the magnetic field $|\nabla \cdot \B|_{\infty}$ for 1000 time steps. Right: Time evolution of the ratio between the material time step of the \scheme scheme and the magneto-sonic time step of a fully explicit finite volume scheme }
		\label{fig.divB_dt_Blast-Rotor}
	\end{center}
\end{figure}

Furthermore, in Figure \ref{fig.divB_dt_Blast-Rotor} we also monitor the error in the divergence of the magnetic field, which is measure in the $L_{\infty}$ norm over the whole computational domain at all time instants. For both the blast wave and the rotor problem we achieve a discrete divergence-free condition up to machine accuracy. Alongside, the new semi-implicit schemes allows the time step to be still larger than the one of a fully explicit scheme, even if strong shocks occur in these test cases. Looking at the right panel of Figure \ref{fig.divB_dt_Blast-Rotor}, we clearly see that the \scheme method can lead to a time marching scheme with a time step that is up to four times larger than a classical explicit solver.

\subsection{Orszag-Tang vortex}
\begin{figure}[!htbp]
	\begin{center}
		\begin{tabular}{cc}
			\includegraphics[trim=2 2 2 2,clip,width=0.47\textwidth]{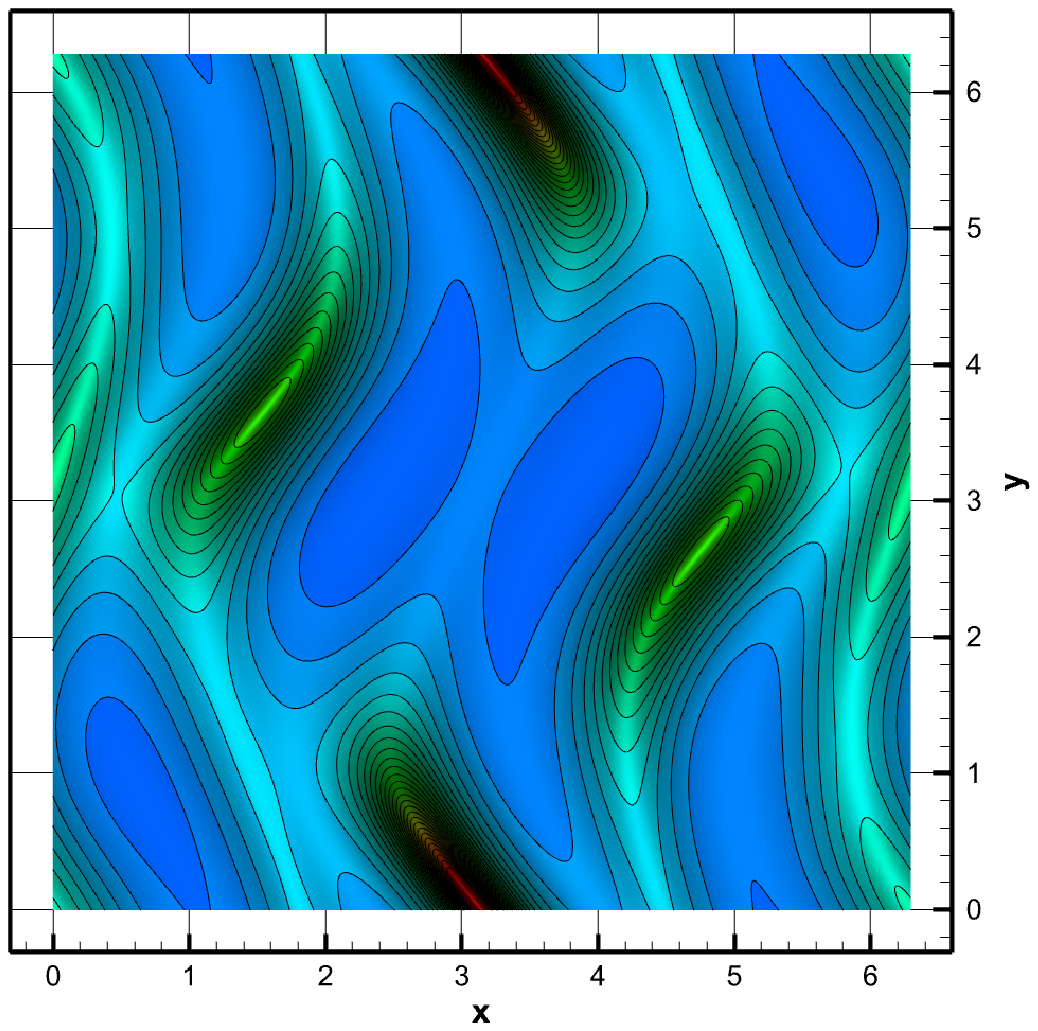} &
			\includegraphics[trim=2 2 2 2,clip,width=0.47\textwidth]{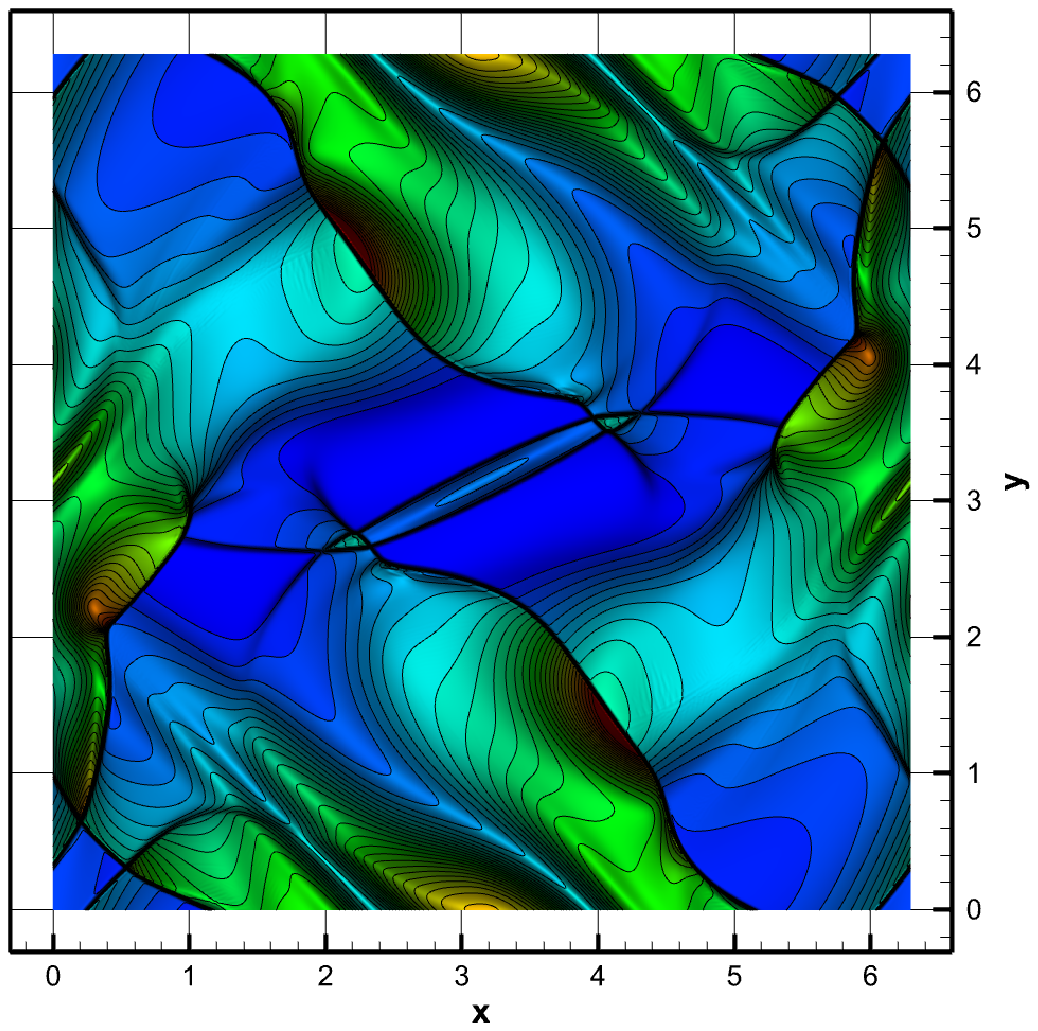} \\
			\includegraphics[trim=2 2 2 2,clip,width=0.47\textwidth]{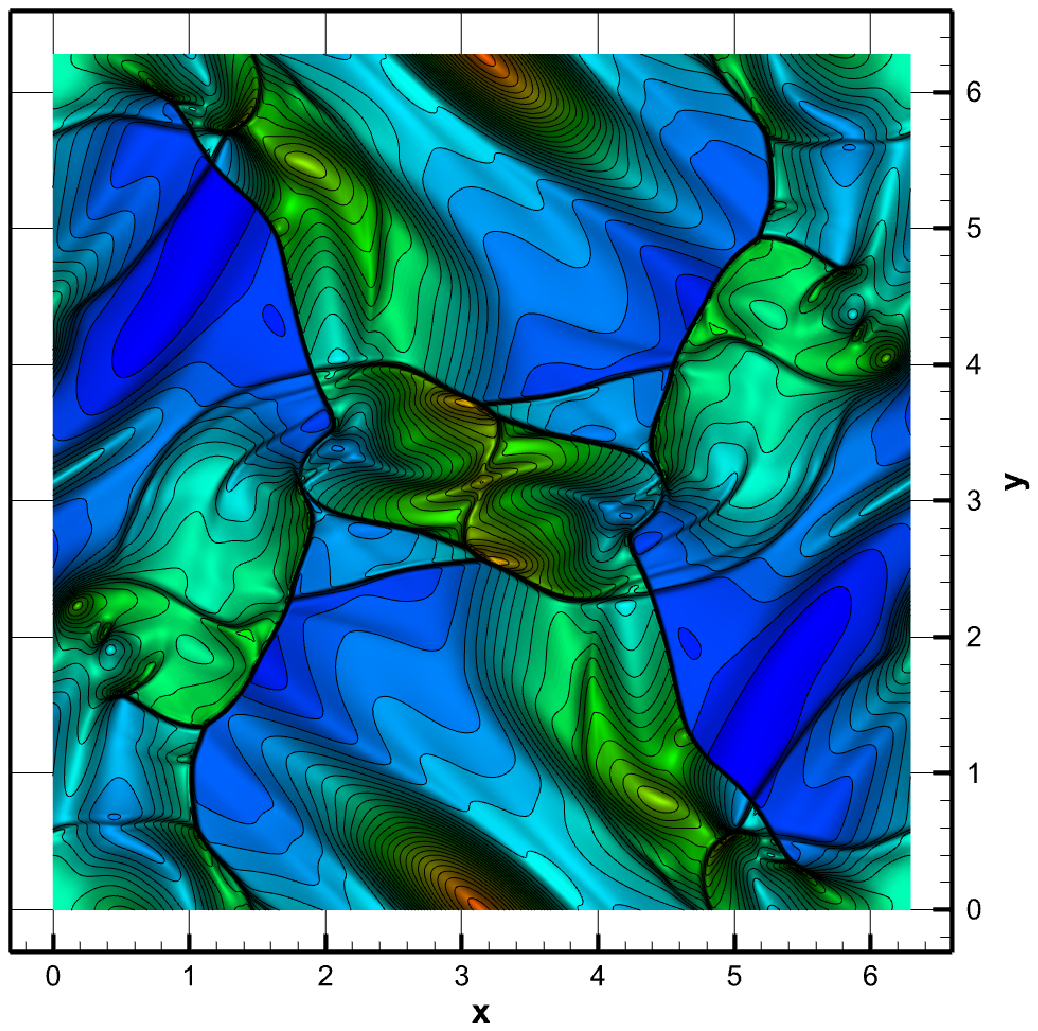} &
			\includegraphics[trim=2 2 2 2,clip,width=0.47\textwidth]{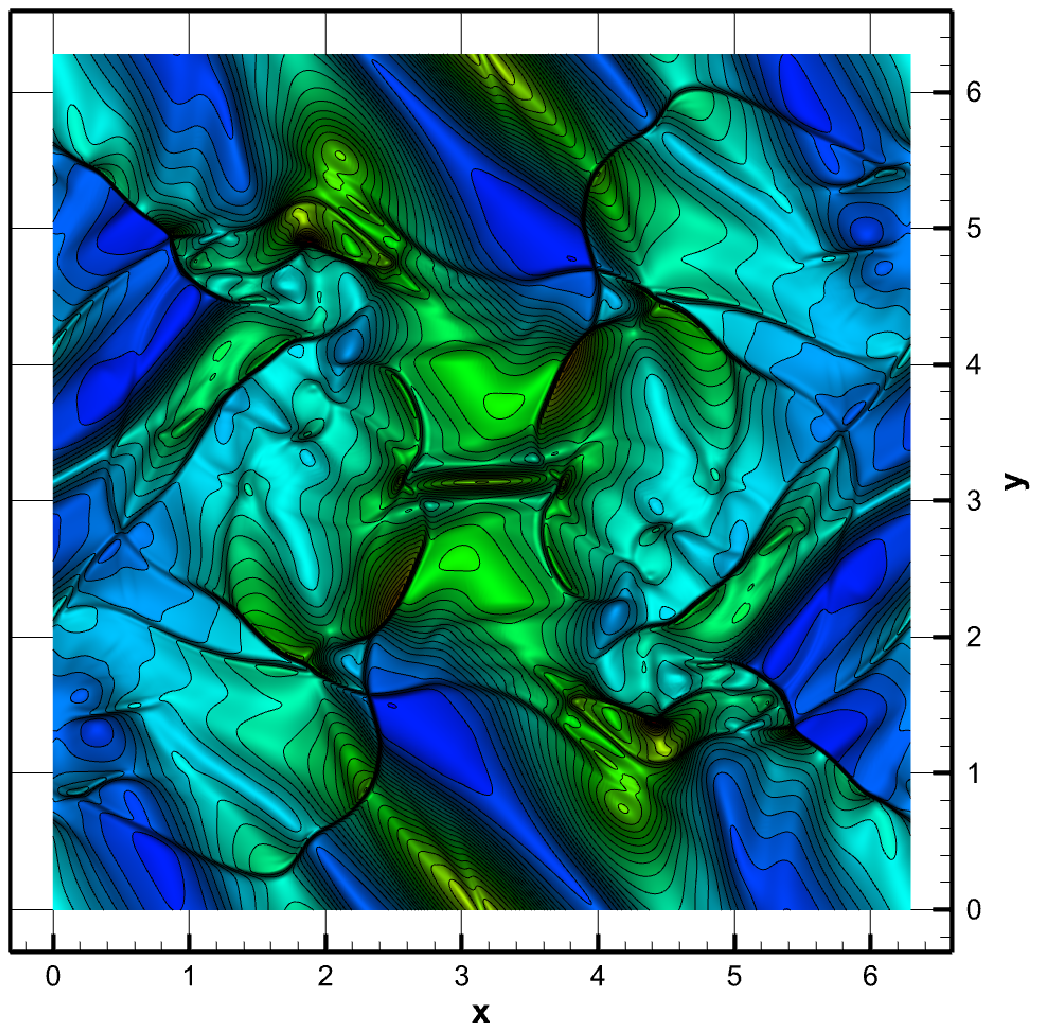} \\
		\end{tabular}
		\caption{Orszag-Tang vortex at time $t=0.5$, $t=2$, $t=3$ and $t=4$ (from top left to bottom right). Results are given for the pressure $p$ with 31 equidistant contour lines in the interval $[0.5,5.0]$.}
		\label{fig.OrszagTang}
	\end{center}
\end{figure}

We consider the well-known Orszag-Tang vortex for the two-dimensional ideal MHD equations \cite{OrszagTang1979,dahlburg1989,picone1991}.
The initial condition is smooth and reads
\begin{equation}
    (\rho, u, v, w, p, B_x, B_y, B_z) =
    \left(\frac{25}{36\pi}, - \sin(2 \pi y), \sin(2 \pi x), 0, \frac{5}{12 \pi}, - \frac{\sin (2 \pi y)}{\sqrt{4 \pi} }, \frac{\sin (4 \pi y)}{\sqrt{4 \pi} }, 0\right).
\end{equation}
The magnetic field $\B$ is generated by the $z$-component of the magnetic potential
\begin{equation}
    A_z = \frac{1}{8 \pi^{3/2}} (\cos(2 \pi x) + 2\cos(\pi y)).
\end{equation}
Though initially smooth, the dynamics develop shocks along the diagonal direction in combination with a vortex located at the center of the computational domain, which is defined by $\Omega=[0;1]\times [0;1]$ with periodic boundaries. We use a mesh made of $N_x \times N_y=1024 \times 1024$ cells. Figure \ref{fig.OrszagTang} depicts snapshots of the pressure computed with our second-order \scheme scheme at times $t = \{0.5, 1, 2, 4\}$.
The features of the flows and the development of the shock are captured well and the results are qualitatively in very good agreement with those reported in the literature \cite{SIMHD_Dumbser2019,divFreeMHD_BalDum2015,balsara2010HLLE}.

\subsection{Field loop advection}
\begin{figure}[!htbp]
	\begin{center}
		\begin{tabular}{c}
			\includegraphics[trim=2 2 2 2,clip,width=0.9\textwidth]{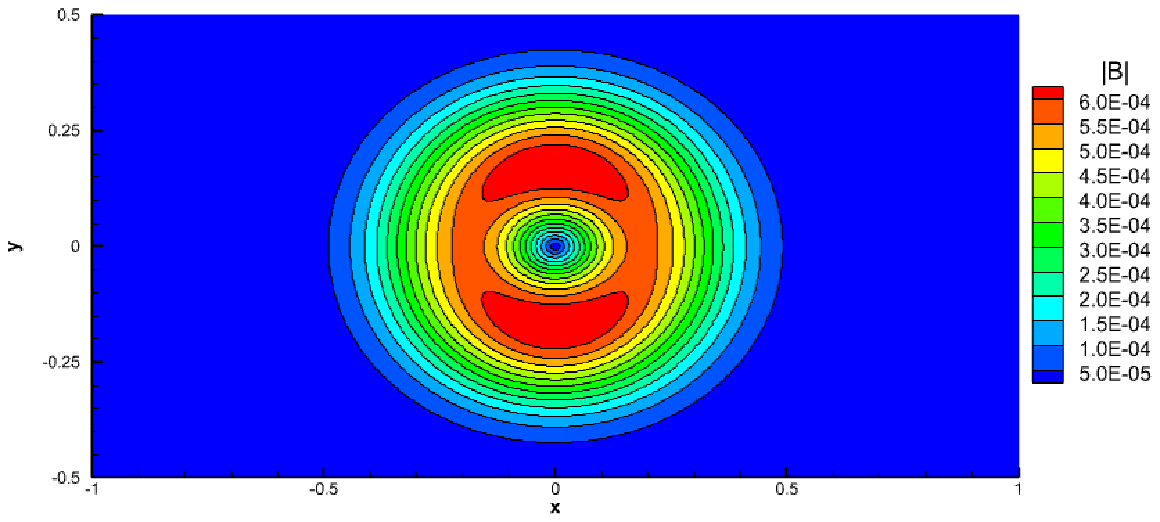} \\
			\includegraphics[trim=2 2 2 2,clip,width=0.9\textwidth]{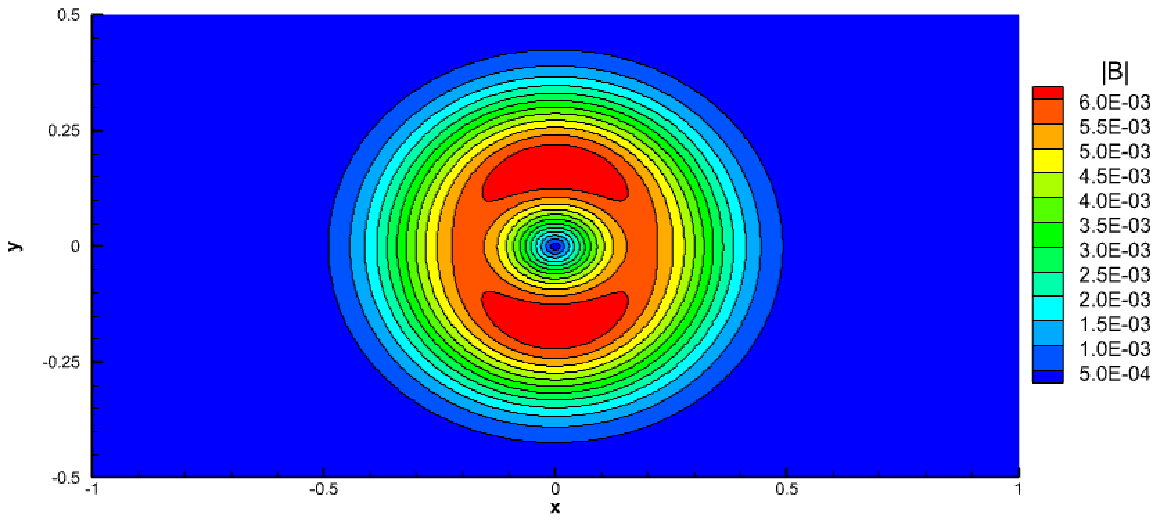} \\
			\includegraphics[trim=2 2 2 2,clip,width=0.9\textwidth]{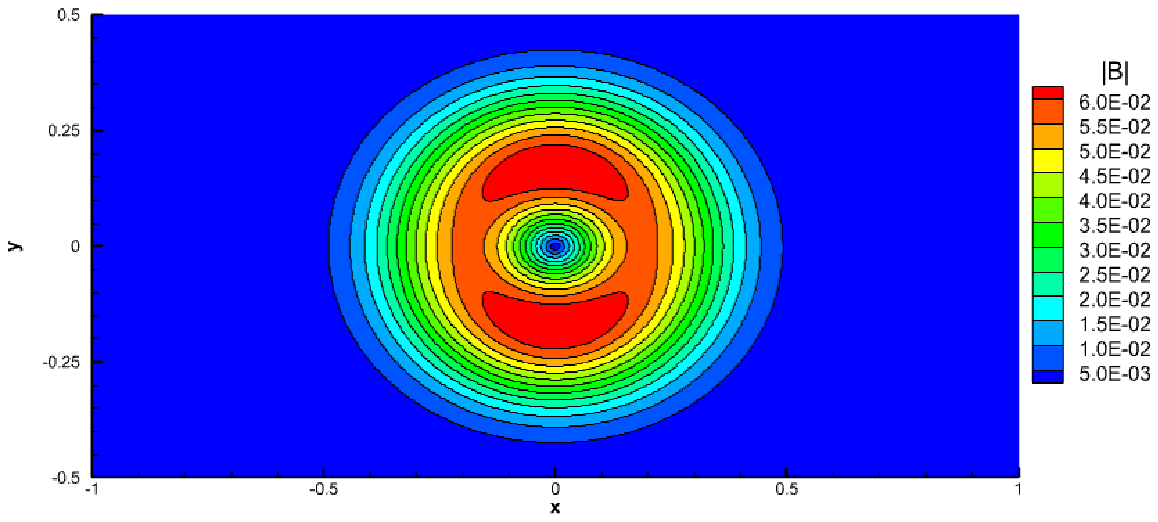} \\
			\includegraphics[trim=2 2 2 2,clip,width=0.9\textwidth]{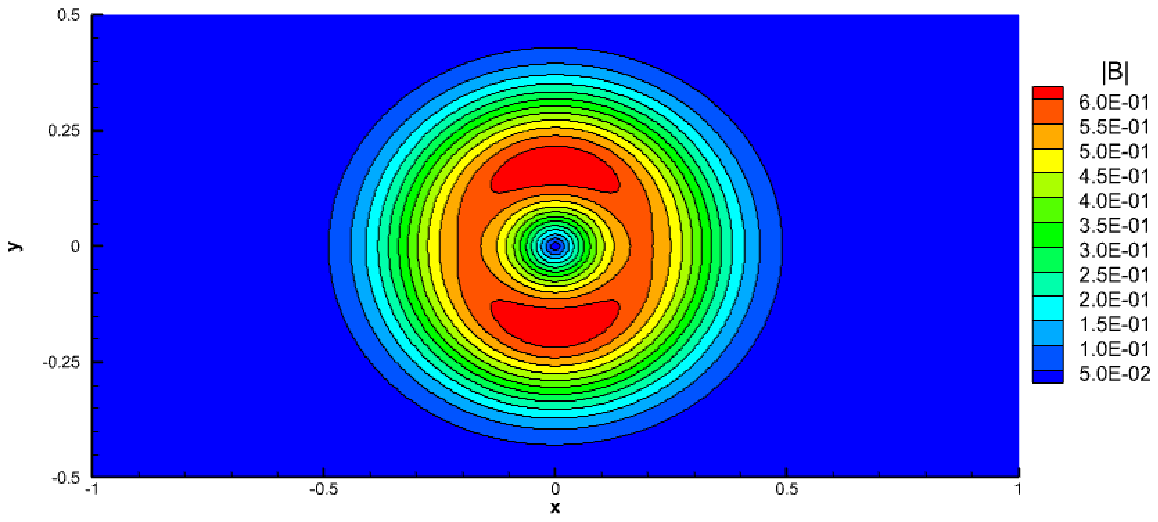} \\
		\end{tabular}
		\caption{Field loop advection problem at time $t_f=1$. Magnitude of the magnetic field with $A_0 = \{ 10^{-3}, 10^{-2}, 10^{-1}, 10^{0}\}$ (from top to bottom).}
		\label{fig.FieldLoop}
	\end{center}
\end{figure}

We solve the magnetic field loop advection problem originally proposed in \cite{Gardiner2005}, in the low Mach number version of \cite{SIMHD_Dumbser2019} with $\rho = 1$, an initial velocity field $(u,v,w) = (2,1,0)$ and the pressure $p = 10^5$. The computational domain is given by $\Omega=[-1;1]\times[-0.5;0.5]$ with periodic boundaries, and we use a mesh composed of $N_x \times N_y = 512 \times 256$ cells. The magnetic field is prescribed by the magnetic vector potential in $z$-direction given by
\begin{equation}
    A_z = \begin{cases}
        A_0 (R- r) & \text{if} \quad r \leq R,\\
        0          & \text{else},
    \end{cases}
\end{equation}
where $A_0=10^{-3}$, $R = 0.3$ and $r = \sqrt{x^2 + y^2}$. This test case is very challenging because of a singular point in the magnetic field $\B$ at the center of the computational domain. This setting places the flow in a low acoustic Mach number regime $\Mc \approx 6 \cdot 10^{-3}$ and high Alfv\'en Mach number regime $\Ma \approx 7.9 \cdot 10^3$. To make the test more challenging, we perform four different simulations with a progressively increasing value of $A_0$, reaching a configuration with an Alfv\'en Mach number $\Ma \approx 7.9$ corresponding to $A_0=1$. The final time is chosen to be $t_f = 1$, and the results are depicted in Figure \ref{fig.FieldLoop}. The shape of the magnitude of the magnetic field remains qualitatively unchanged despite the difference in the Alfv\'en Mach number, hence demonstrating that the new \scheme scheme is independent of the magnetic scales. Indeed, the time step for all the simulations is the same and is only determined by the convective velocity $(u,v)=(2,1)$ which transports the loop diagonally through the domain.

\begin{figure}[!htbp]
	\begin{center}
		\begin{tabular}{c}
			\includegraphics[trim=2 2 2 2,clip,width=0.7\textwidth]{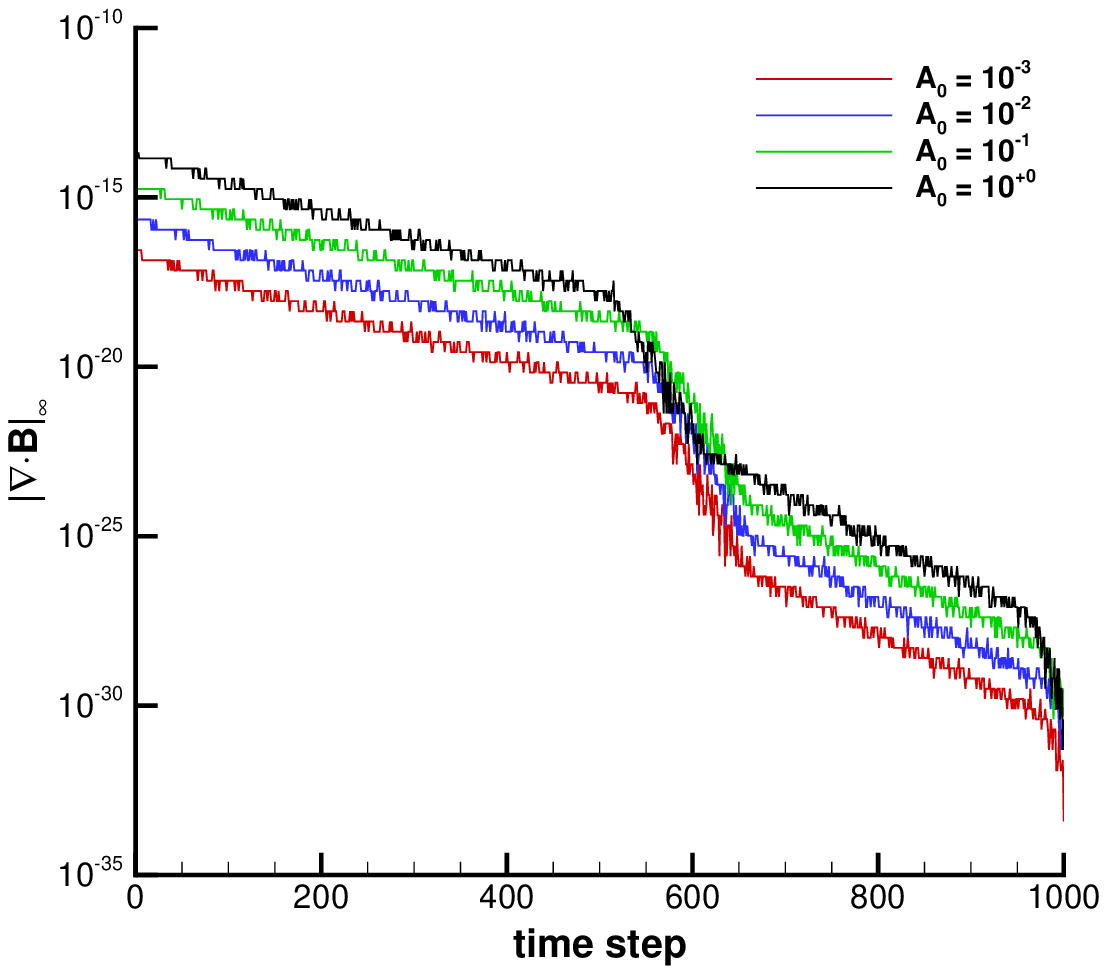}
		\end{tabular}
		\caption{Field loop advection problem with $A_0 = \cdot \{ 10^{-3}, 10^{-2}, 10^{-1}, 10^{0}\}$. Evolution of the infinity norm of the divergence of the magnetic field $|\nabla \cdot \B|_{\infty}$ for 1000 time steps of each simulation.}
		\label{fig.FieldLoop_divB}
	\end{center}
\end{figure}

In Figure \ref{fig.FieldLoop_divB} we show the errors related to the divergence of the magnetic field for all four simulations, demonstrating that the discrete divergence-free condition is preserved up to machine accuracy. The infinity norm of the divergence of $\B$ is considered over the whole domain and at each time level.

\subsection{Magnetized Kelvin Helmholtz instability}
\begin{figure}[!htbp]
	\begin{center}
		\begin{tabular}{cc}
			\includegraphics[trim=2 2 2 2,clip,width=0.48\textwidth]{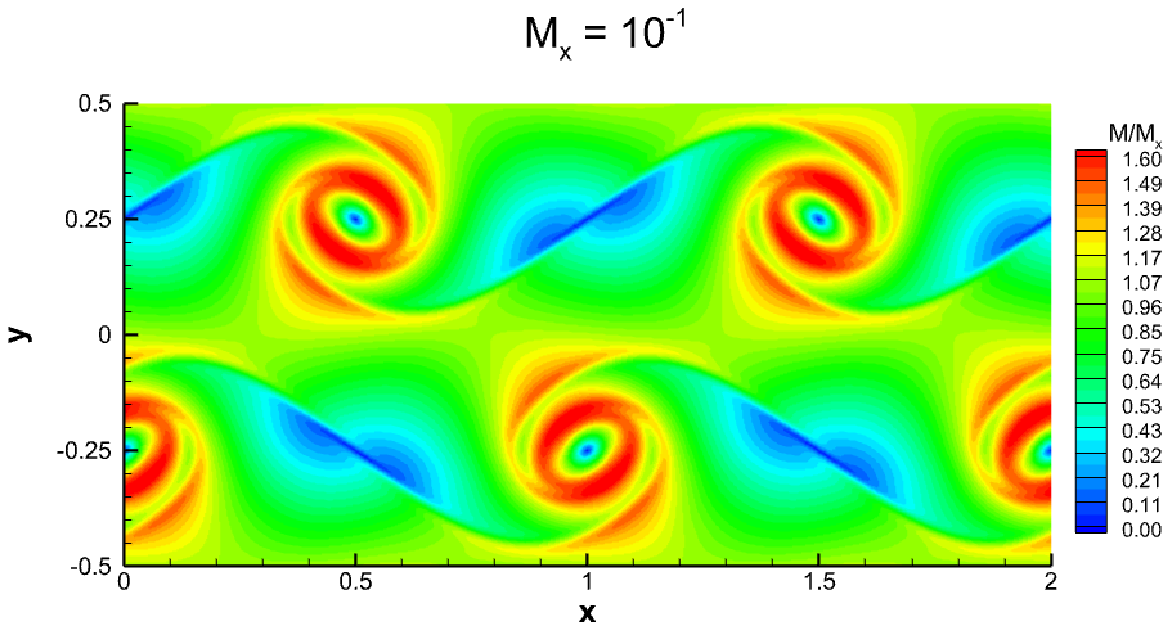} &
			\includegraphics[trim=2 2 2 2,clip,width=0.48\textwidth]{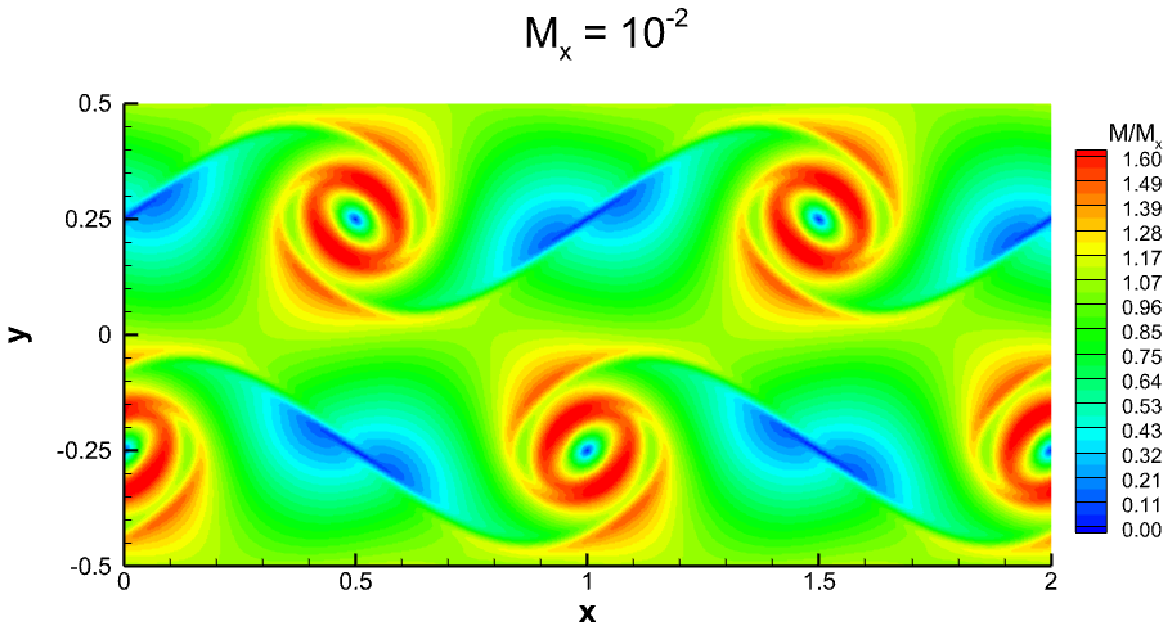} \\
			\includegraphics[trim=2 2 2 2,clip,width=0.48\textwidth]{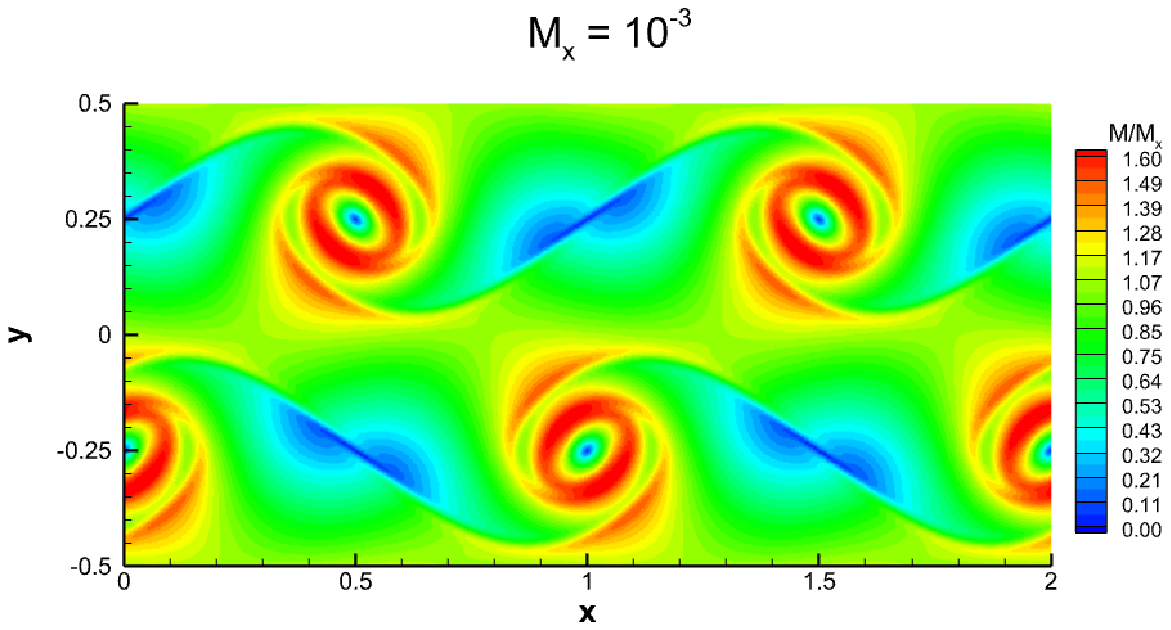} &
			\includegraphics[trim=2 2 2 2,clip,width=0.48\textwidth]{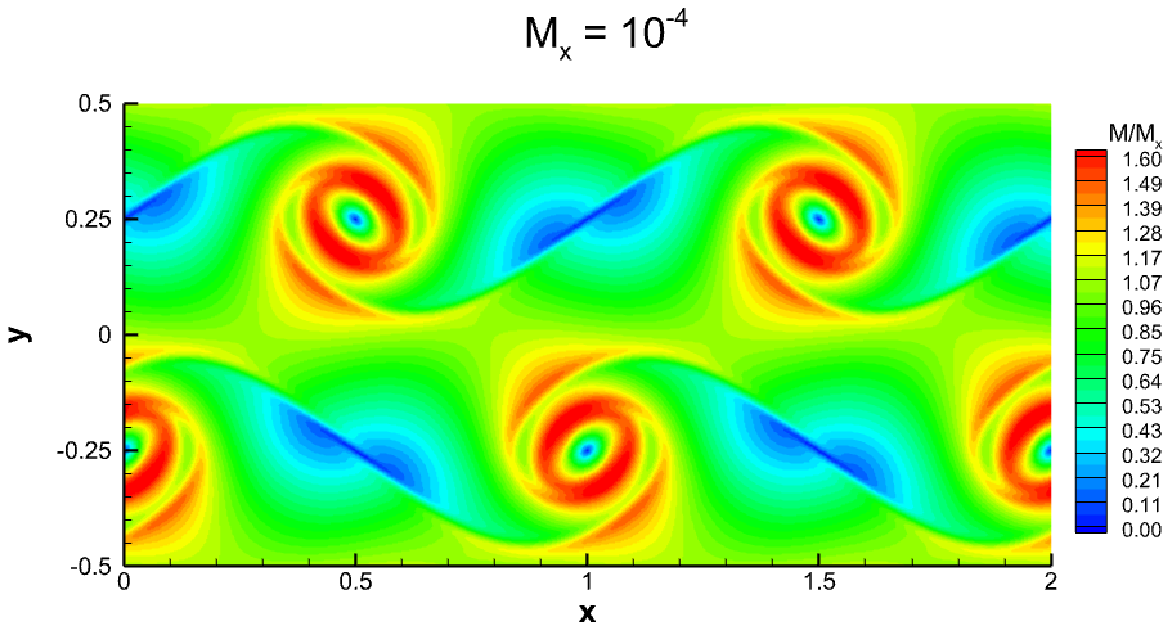} \\
		\end{tabular}
		\caption{Kelvin-Helmholtz instability test at time $t/t_{\max}=1/6$ for $M_x=\{10^{-1},10^{-2},10^{-3},10^{-4}\}$. Distribution of the acoustic Mach number rescaled by the corresponding value of $M_x$ ($\Mc/M_x$) with 31 equidistant contour lines in the interval $[0,1.6]$.}
		\label{fig.KH_M}
	\end{center}
\end{figure}

\begin{figure}[!htbp]
	\begin{center}
		\begin{tabular}{c}
			\includegraphics[width=0.7\textwidth,keepaspectratio=true]{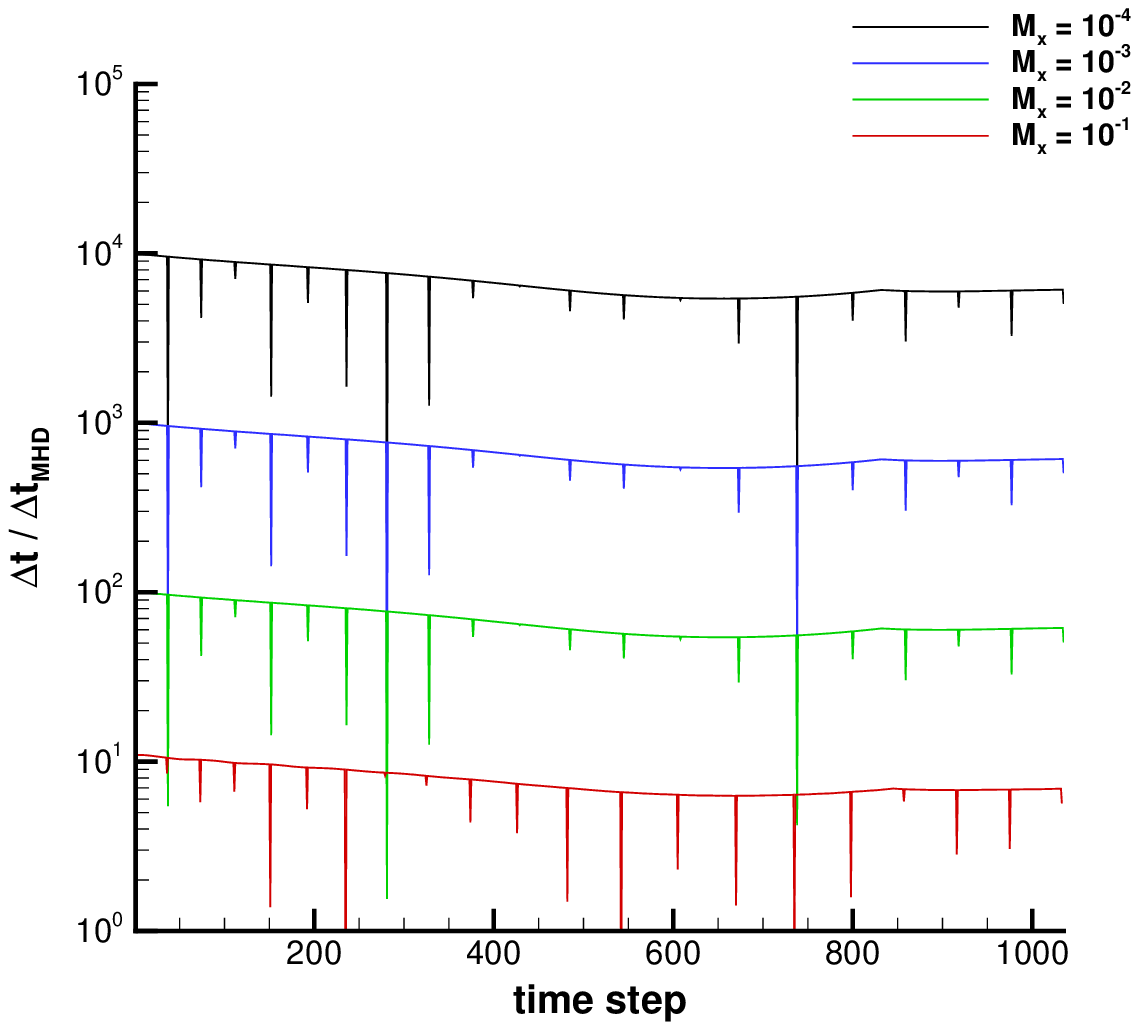}
		\end{tabular}
		\caption{Kelvin-Helmholtz instability test. Time evolution of the ratio between the material time step of the \scheme scheme and the magneto-sonic time step of a fully explicit finite volume scheme for different values of $M_x$.}
		\label{fig.KHdt}
	\end{center}
\end{figure}

We consider the setup for the Kelvin-Helmholtz instability recently proposed in \cite{BK2022} for ideal MHD. The computational domain is the rectangular box $\Omega=[0;2]\times[-0.5;0.5]$ with periodic boundaries. It is paved with a computational grid made of $N_x \times N_y = 512 \times 256$ control volumes. A uniform density $\rho=\gamma$ with $\gamma=1.4$ and a uniform pressure $p=1$ is initially assigned to the fluid. The velocity and the magnetic field are prescribed as a function of the parameter $M_x$, which represents the maximum acoustic Mach number of the horizontal flow. The velocity field is given by
\begin{equation}
	u=M_x \, (1-2\eta(x)), \qquad v = 0.1 \, M_x \, \sin(2\pi x),
\end{equation}
with the perturbation $\eta(x)$
\begin{equation}
	\eta(x) = \begin{cases}
		0.5 \left(  1 + \sin(16\pi(y+0.25)) \right) & \text{if } y \geq -\frac{9}{32} \text{ and } y < -\frac{7}{32},\\
		1 & \text{if } y \geq -\frac{7}{32} \text{ and } y < \frac{7}{32}, \\
		0.5 \left(  1 - \sin(16\pi(y-0.25)) \right) & \text{if } y \geq \frac{7}{32} \text{ and } y < \frac{9}{32},\\
		0 & \text{otherwise}.
	\end{cases}
\end{equation}
The magnetic field at the initial time is uniform and horizontal, i.e. $B_x=0.1 \, M_x$, and the associated $z-$component of the magnetic potential is given by $A_0=B_x \, y$. This configuration at time $t=0$ leads to a unit adiabatic sound speed and to a minimum Alfv\'en mach number $M_{b,\min{}}=11.82$. The simulation is run with different values of $M_x$ until the final time $t_f=t_{\max}/6$ with $t_{\max}=4.8/M_x$.
The flow develops a smooth and resolved interface across the shear flow, hence yielding to convergent results in the early stages of the dynamics. Figure \ref{fig.KH_M} plots the acoustic Mach number at the end of the simulation. It is rescaled with the prescribed value of $M_x$. Note that all the results are almost identical since excessive numerical dissipation is avoided thanks to the central implicit discretization of the pressure and the magnetic terms. Fully explicit schemes would not preserve the vortical structures as the value of $M_x$ gets smaller, as noticed in \cite{BK2022}. Furthermore, our \scheme scheme can run these test cases with a time step that is up to 10,000 times larger compared to a fully explicit method, as shown in Figure \ref{fig.KHdt}.

\subsection{3D cloud-shock interaction problem}
\begin{figure}[!htbp]
	\begin{center}
		\begin{tabular}{cc}
			\includegraphics[trim=2 2 2 2,clip,width=0.47\textwidth]{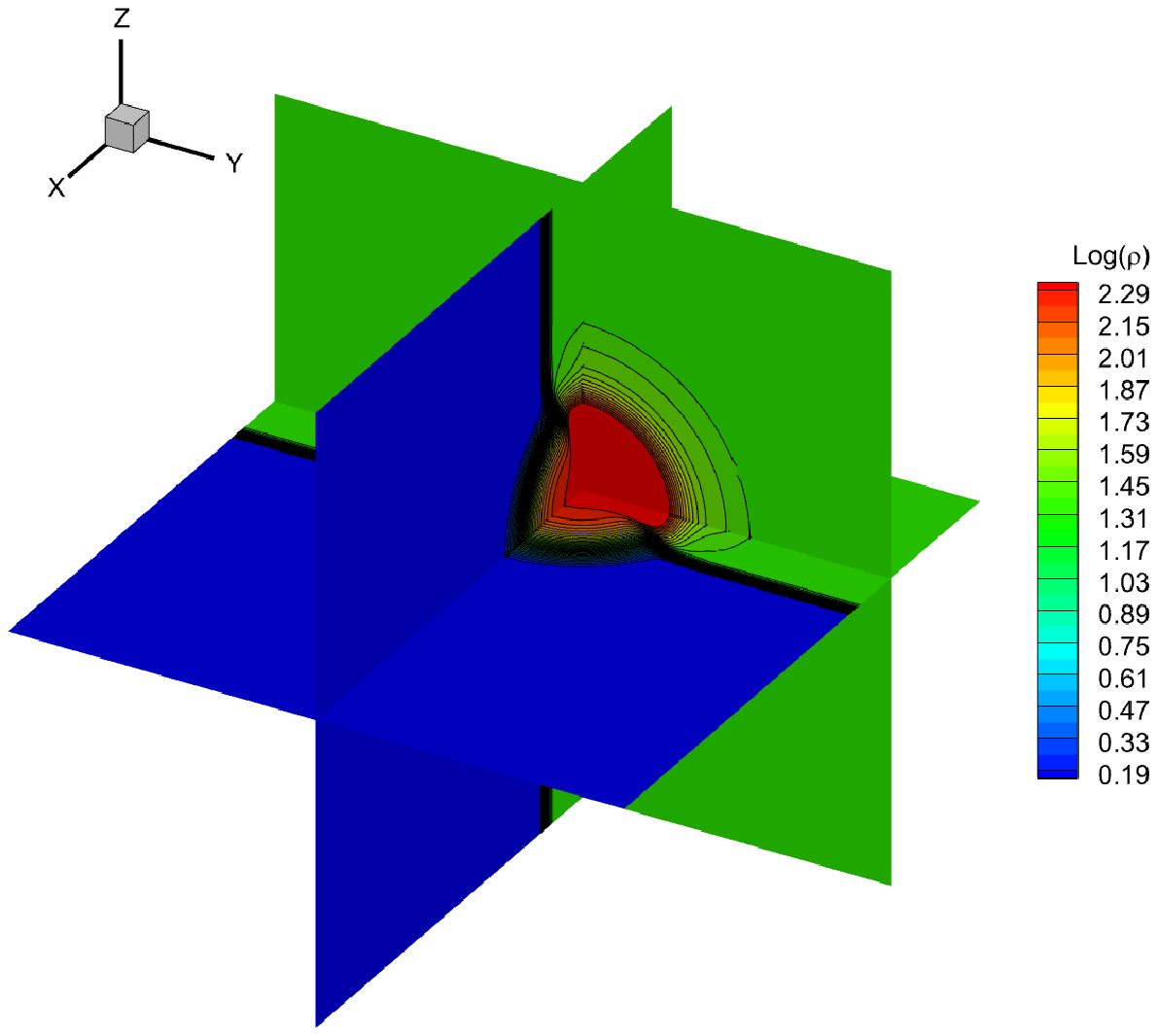} &
			\includegraphics[trim=2 2 2 2,clip,width=0.47\textwidth]{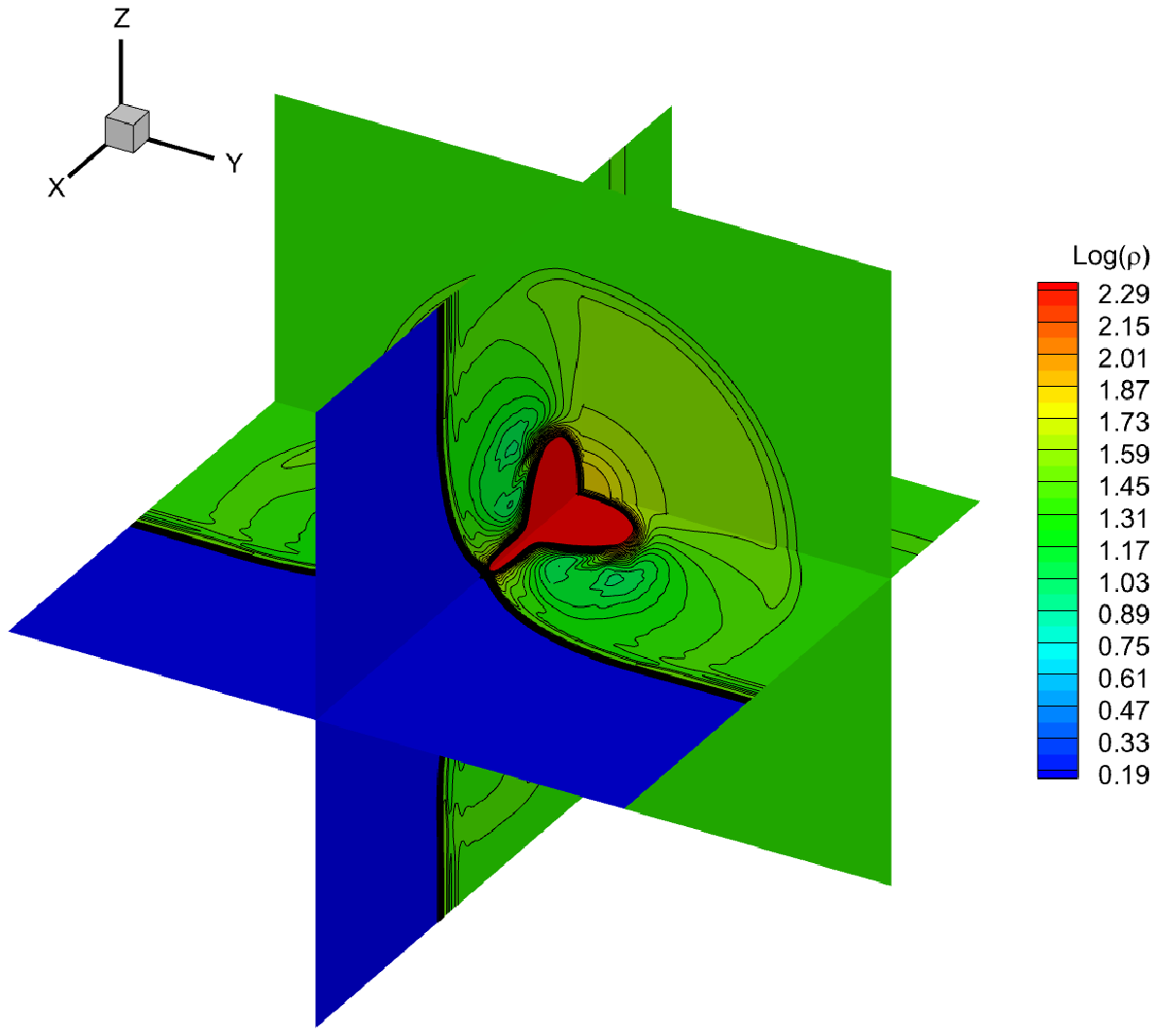} \\
			\includegraphics[trim=2 2 2 2,clip,width=0.47\textwidth]{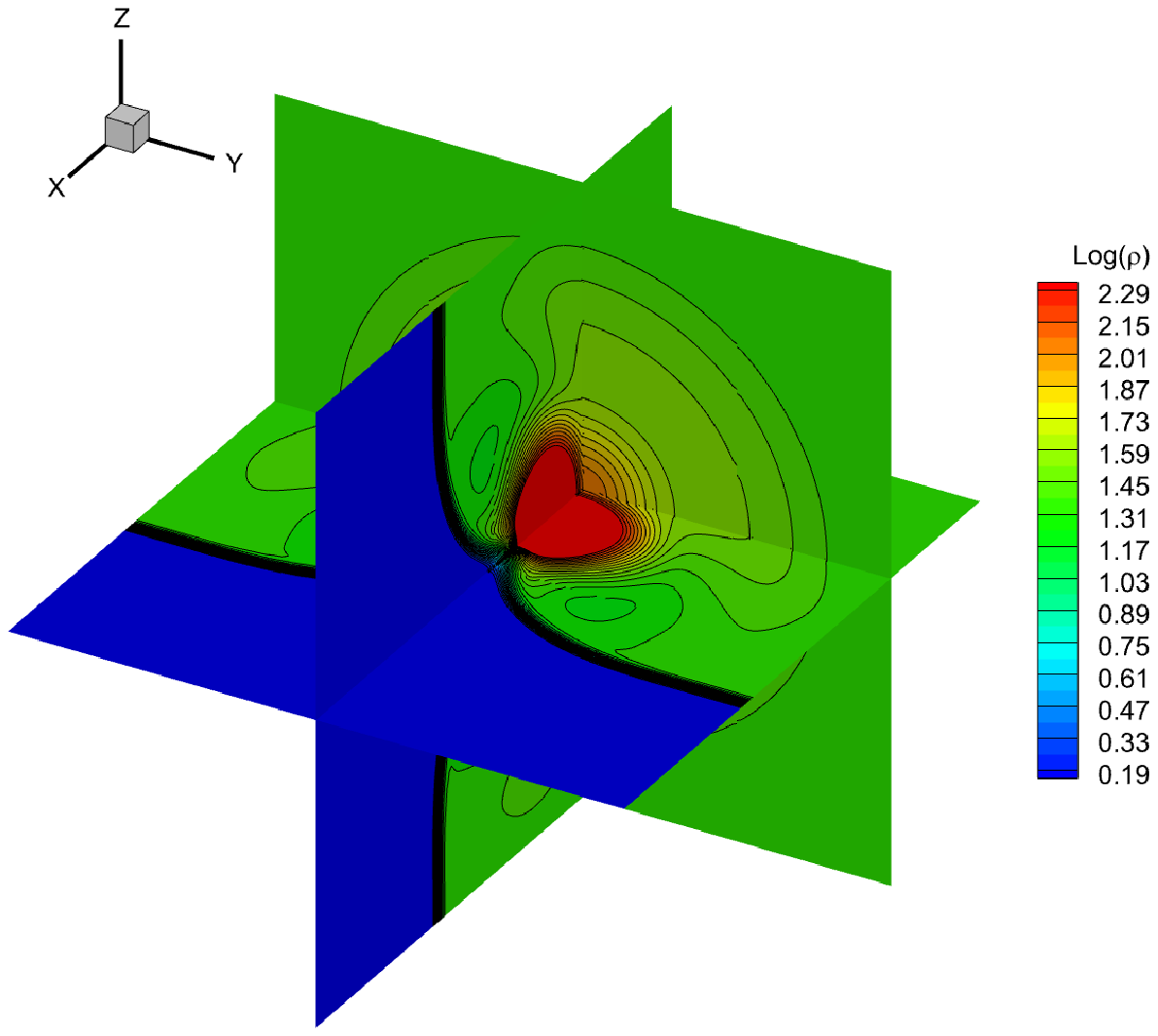} &
			\includegraphics[trim=2 2 2 2,clip,width=0.47\textwidth]{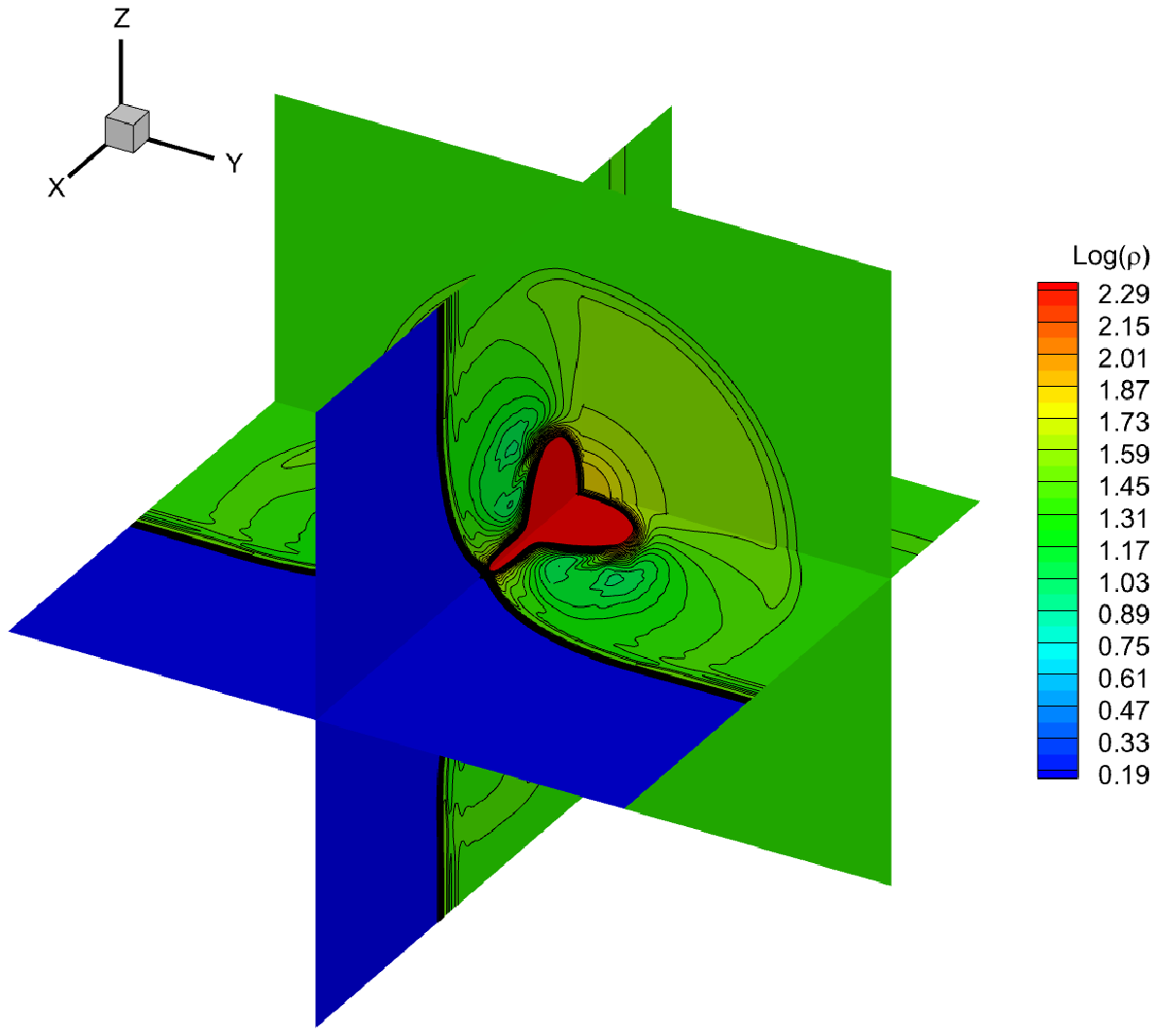} \\
			\includegraphics[trim=2 2 2 2,clip,width=0.47\textwidth]{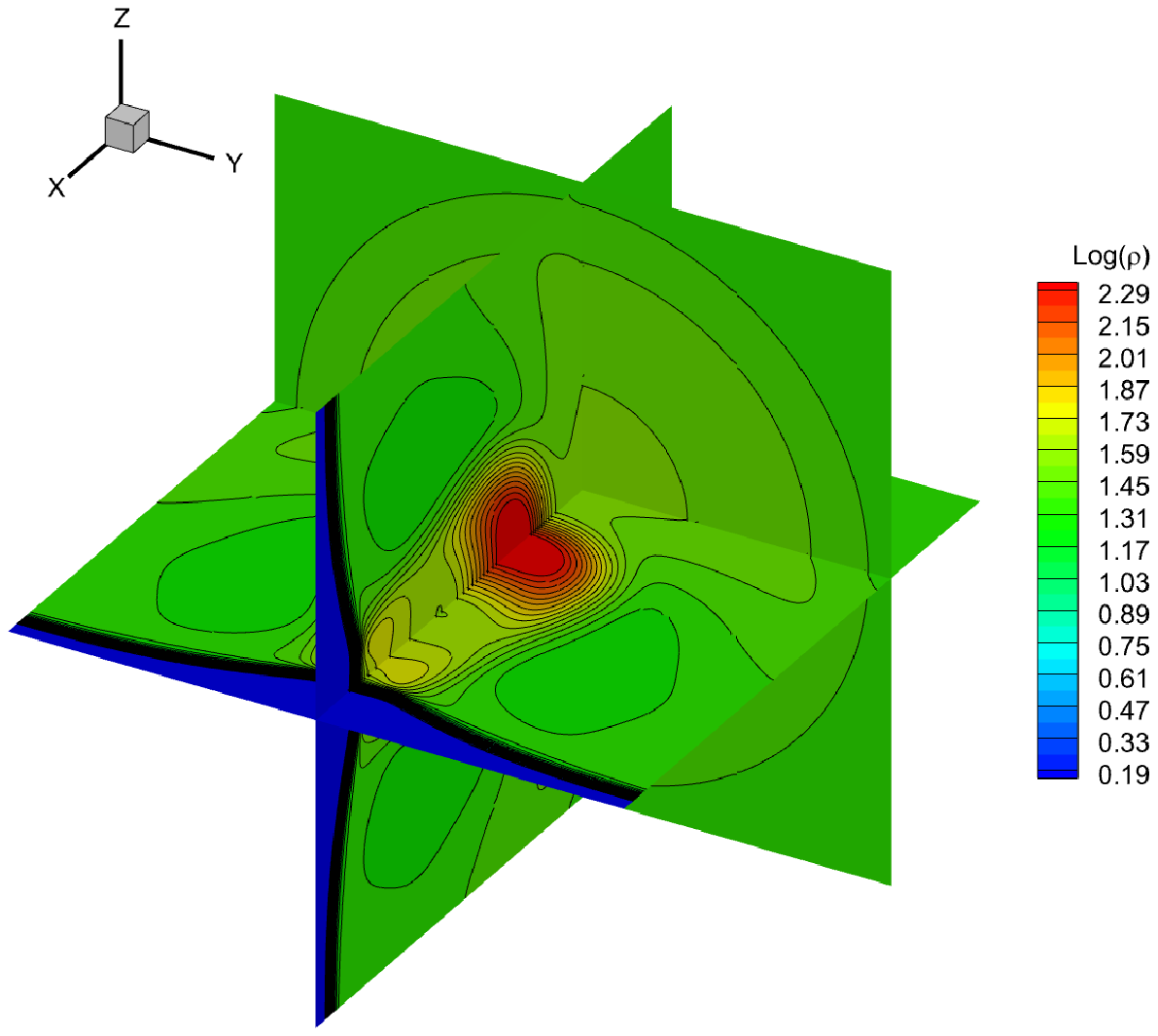} &
			\includegraphics[trim=2 2 2 2,clip,width=0.47\textwidth]{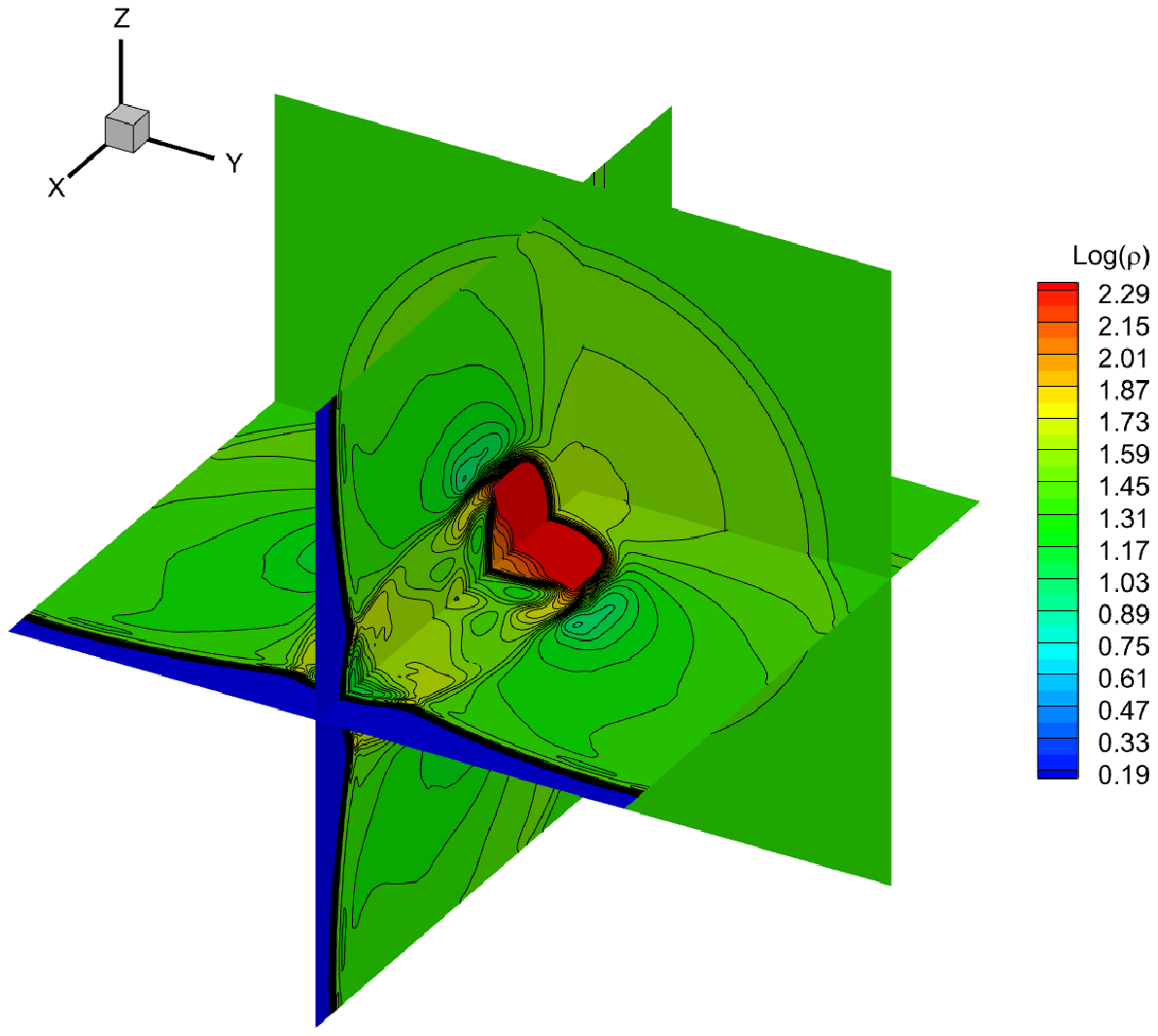} \\
		\end{tabular}
		\caption{3D cloud-shock interaction problem at time $t=0.02$ (top), $t=0.04$ (middle) and $t=0.06$ (bottom). Distribution of the logarithm of density obtained with first (left) and second (right) order \scheme schemes. Slices along the planes $x=0.25$, $y=0.5$ and $z=0.5$.}
		\label{fig.cloudshock}
	\end{center}
\end{figure}

As last test case we consider the 3D cloud-shock interaction problem. The computational domain is given by $\Omega=[0;1]^3$ with Dirichlet boundary conditions, and we employ a computational mesh with a total number of $N_x=N_y=N_z=180^3$ cells, corresponding to approximately $16.5 \cdot 10^6$ degrees of freedom. The initial condition is taken from \cite{helzel2013} where the initial shock front is located at $x_0 = 0.0$ with
\begin{multline}
    (\rho, u, v, w, p, B_x, B_y, B_z) = \\
    \begin{cases}
        \begin{aligned}
        (\hskip-5pt&&3.86859, &&11.2536, \quad &0, &0, &&167.345, \quad &0, &2.1826182, &&-2.1826182~)& \quad\text{if } x < x_0 \\
        (\hskip-5pt&&1, &&0, \quad & 0,& 0, &&1,\quad &0, &0.56418958, &&0.56418958~)& \quad \text{if } x \geq x_0
        \end{aligned}
    \end{cases}.
\end{multline}
The associated magnetic vector potential writes
\begin{equation}
    \A = \begin{cases}
        \begin{aligned}
            (\hskip-5pt&&2.1826182 y, &&0, \quad && -2.1826182 (x - 0.05) ~)^T& \quad\text{if } x < x_0 \\
            (\hskip-5pt&&-0.56418958 y, &&0, \quad &&  0.56418956 (x - 0.05)~)^T& \quad \text{if } x \geq x_0
        \end{aligned}
    \end{cases}.
\end{equation}
This shock front travels towards a spherical cloud of density $\rho = 10$ situated at $\xx = (0.25,0.5,0.5)$ with radius $r=0.15$. The final time is set to $t_f = 0.06$.
In Figure \ref{fig.cloudshock}, we compare the first against the second order \scheme scheme. The density profiles along the planes cutting through the center of the initial cloud position are displayed. The first snapshot is taken shortly after the shock front hits the density cloud, the second one after the shock front passes through the cloud, and the third one shortly before the shock front leaves the computational domain. One can observe that the second-order scheme greatly improves the resolution of the flow structures, which are qualitatively in very good agreement with the results presented in \cite{helzel2013}. This test case demonstrates the capability of the \scheme scheme to handle genuinely 3D magnetized flows at reasonable compute times due to the linear implicit structure of our novel \scheme scheme.

\section{Conclusions}
\label{sec.concl}
In this work, we have introduced a new divergence-free semi-implicit finite volume method for the simulation of ideal magnetized flows in multiple space dimensions. To tackle the time multi-scale nature of the governing equations, a flux splitting technique is adopted, which allows to separate the convective terms, the hydrodynamic pressure contribution and its magnetic counterpart. More specifically, the non-linear convective contribution is discretized explicitly with a classical finite volume scheme, while we use a finite difference semi-implicit strategy to compute the hydrodynamic pressure and the magnetic fluxes, that ultimately yields two linear systems to be solved at each time step of the algorithm. We underline that these two systems are uncoupled, thus firstly we solve for the unknown magnetic field, and secondly we determine the unknown pressure contribution. This approach makes the scheme very efficient, especially in the low acoustic or Alfv\'en regimes, since the time step is restricted by a milder CFL stability condition that is only based on the material waves. Furthermore, no numerical dissipation is added in the discretization of the implicit pressure fluxes, hence the scheme does not pollute the solution in the low acoustic Mach number regime, contrarily to standard explicit finite volume methods. To satisfy the solenoidal property of the magnetic field, the equations are written in terms of the magnetic vector potential in a multi-dimensional setting, and a set of finite difference operators are specifically designed to mimic the algebraic property of the zero div-curl at the discrete level, thus preserving the structure of the underlying model. The entire spatial discretization is cell-centered, and second-order of accuracy is achieved by means of a TVD reconstruction procedure in the finite volume scheme. The class of semi-implicit IMEX Runge-Kutta schemes is employed to reach second-order accuracy also in time. The new schemes have been thoroughly tested in terms of accuracy and robustness against a wide set of test cases, addressing different acoustic and Alfv\'en Mach number regimes, demonstrating the capability of our novel methods to span a wide range of physical settings in ideal MHD.

In the future, we plan to carry out a theoretical and numerical study of the limit models in the low Alfv\'en Mach number, which is still rather unexplored. The extension of the presented approach to unstructured meshes is also foreseen, with particular care in the development of divergence-free operators that can guarantee the solenoidal property of the magnetic field. Finally, future research lines will consider the application to viscous and resistive MHD equations in the context of Tokamak simulations.

\section*{Acknowledgments}
WB received financial support by Fondazione Cariplo and Fondazione CDP (Italy) under the project No. 2022-1895 and by the Italian Ministry of University and Research (MUR) with the PRIN Project 2022 No. 2022N9BM3N. WB is member of the GNCS-INdAM (\textit{Istituto Nazionale di Alta Matematica}) group.

	%
	\section*{Declarations}

	\paragraph{Funding.} WB received financial support by Fondazione Cariplo and Fondazione CDP (Italy) under the project No. 2022-1895 and by the Italian Ministry of University (MUR) with the PRIN Project 2022 No. 2022N9BM3N.

	\paragraph{Conflicts of interest.} The authors declare that they have no conflict of interest.

	\paragraph{Availability of data and material.} Data and material are available upon reasonable request addressed to the corresponding author.

	\paragraph{Code availability.} The code is written in Fortran programming language and is available upon reasonable request addressed to the corresponding author.

	\paragraph{Ethics approval.} Not applicable.

	\paragraph{Consent to participate.} Not applicable.

	\paragraph{Consent for publication.} Not applicable.

	\section*{Data availability}
	The datasets generated during the current study are available from the corresponding author upon reasonable request.

	\bibliographystyle{plain}
	\bibliography{biblio}   

\end{document}